%% file: main.tex
\documentclass[12pt]{amsart}
\usepackage[margin=1.25in]{geometry}

\usepackage{amsmath}
\usepackage{amsfonts}
\usepackage{amssymb}
\usepackage{amscd}
\usepackage{graphicx}
\usepackage[abbrev,alphabetic]{amsrefs}
\RequirePackage[dvipsnames,usenames]{color}
\usepackage{soul,xcolor}
\setstcolor{red}
\usepackage{stmaryrd}
\usepackage{multicol}
\usepackage{url}

\newcommand\hidden[1]{} 
\usepackage{mathtools}
\newtagform{hidden}[\hidden]{}{}

\usepackage{hyperref}

\usepackage{amsthm}
\usepackage{comment}
\usepackage[all,cmtip]{xy}
\usepackage{tikz-cd}
\usetikzlibrary{cd}

\makeatletter
\def\@tocline#1#2#3#4#5#6#7{\relax
  \ifnum #1>\c@tocdepth 
  \else
    \par \addpenalty\@secpenalty\addvspace{#2}%
    \begingroup \hyphenpenalty\@M
    \@ifempty{#4}{%
      \@tempdima\csname r@tocindent\number#1\endcsname\relax
    }{%
      \@tempdima#4\relax
    }%
    \parindent\z@ \leftskip#3\relax \advance\leftskip\@tempdima\relax
    \rightskip\@pnumwidth plus4em \parfillskip-\@pnumwidth
    #5\leavevmode\hskip-\@tempdima
      \ifcase #1
       \or\or \hskip 1em \or \hskip 2em \else \hskip 3em \fi%
      #6\nobreak\relax
    \hfill\hbox to\@pnumwidth{\@tocpagenum{#7}}\par
    \nobreak
    \endgroup
  \fi}
\makeatother

\hypersetup{
bookmarks,
bookmarksdepth=3,
bookmarksopen,
bookmarksnumbered,
pdfstartview=FitH,
colorlinks,backref,hyperindex,
linkcolor=Sepia,
anchorcolor=BurntOrange,
citecolor=MidnightBlue,
citecolor=OliveGreen,
filecolor=BlueViolet,
menucolor=Yellow,
urlcolor=OliveGreen
}

\newsavebox{\pullback}
\sbox\pullback{%
\begin{tikzpicture}%
\draw (0,0) -- (1ex,0ex);%
\draw (1ex,0ex) -- (1ex,1ex);%
\end{tikzpicture}}

\newsavebox{\pullbackdl}
\sbox\pullbackdl{%
\begin{tikzpicture}%
\draw (-1ex,0ex) -- (0ex,0ex);%
\draw (0ex,-1ex) -- (0ex,0ex);%
\end{tikzpicture}}

\newsavebox{\pushoutdr}
\sbox\pushoutdr{%
\begin{tikzpicture}%
\draw (-1ex,-1ex) -- (-1ex,0ex);%
\draw (-1ex,0ex) -- (0ex,0ex);%
\end{tikzpicture}}

\newcommand{\cred}{\color{red}}
\newcommand{\cora}{\color{black}}

\newcommand{\rup}[1]{\lceil #1 \rceil}
\newcommand{\rdown}[1]{\lfloor #1 \rfloor}
\newcommand{\mustata}{Musta{\c{t}}{\u{a}}}
\renewcommand{\mod}{\ \textrm{mod}\ }
 
\newcommand{\Z}{\mathbb{Z}}
\newcommand{\Q}{\mathbb{Q}}
\newcommand{\R}{\mathbb{R}}
\newcommand{\F}{\mathbb{F}}

\newcommand{\cHom}{\mathcal{H}om}
\newcommand{\cExt}{\mathcal{E}xt}

\newcommand{\Tr}{\mathrm{Tr}}

\newcommand{\bP}{\mathbb{P}}
\newcommand{\bQ}{\mathbb{Q}}

\newcommand{\bZ}{\mathbb{Z}}

\newcommand{\cB}{B}

\newcommand{\cD}{\mathcal{D}}
\newcommand{\cE}{\mathcal{E}}
\newcommand{\cF}{\mathcal{F}}
\newcommand{\cG}{\mathcal{G}}
\newcommand{\cH}{\mathcal{H}}

\newcommand{\cI}{\mathcal{I}}

\newcommand{\cL}{\mathcal{L}}

\newcommand{\cO}{\mathcal{O}}

\newcommand{\cX}{\mathcal{X}}

\newcommand{\MO}{\mathcal{O}}

\newcommand{\m}{\mathfrak{m}}

\DeclareMathOperator{\image}{Im} 
\DeclareMathOperator{\kernel}{Ker}
\DeclareMathOperator{\Ker}{Ker}
\DeclareMathOperator{\cokernel}{Coker}
\DeclareMathOperator{\res}{res}

\DeclareMathOperator{\adj}{adj}
\DeclareMathOperator{\reg}{reg}

\DeclareMathOperator{\Supp}{Supp}
\DeclareMathOperator{\Spec}{Spec}

\DeclareMathOperator{\Hom}{Hom}
\DeclareMathOperator{\Ext}{Ext}

\DeclareMathOperator{\Exc}{Exc}

\DeclareMathOperator{\Tor}{Tor}
\DeclareMathOperator{\Diff}{Diff}

\newcommand{\mydot}{{{\,\begin{picture}(1,1)(-1,-2)\circle*{2}\end{picture}\ }}}

\theoremstyle{plain}
\newtheorem{theorem}{Theorem}[section]

\newtheorem{proposition}[theorem]{Proposition}
\newtheorem{lemma}[theorem]{Lemma}
\newtheorem{corollary}[theorem]{Corollary}
\newtheorem{conjecture}[theorem]{Conjecture}

\newtheorem{claim}[theorem]{Claim}
\newtheorem*{claim*}{Claim}

\newtheorem{theoremIntro}{Theorem}

\theoremstyle{definition}
\newtheorem{definition}[theorem]{Definition}

\newtheorem{assumption}[theorem]{Assumption}
\newtheorem{setting}[theorem]{Setting}
\newtheorem{example}[theorem]{Example}
\newtheorem*{setup*}{Setup}
\newtheorem{step}{Step}

\theoremstyle{remark}
\newtheorem{remark}[theorem]{Remark}
\newtheorem{notation}[theorem]{Notation}

\numberwithin{equation}{theorem}

\AtBeginDocument{%
 \def\MR#1{}
}



\newif\ifshowColoursAndTodoes
\showColoursAndTodoesfalse

\ifshowColoursAndTodoes
\def\todo#1{\textcolor{Mahogany}%
{\footnotesize\newline{\color{Mahogany}\fbox{\parbox{\textwidth-15pt}{\textbf{todo: } #1}}}\newline}}
\def\commentbox#1{\textcolor{Mahogany}%
{\footnotesize\newline{\color{Mahogany}\fbox{\parbox{\textwidth-15pt}{\textbf{comment: } #1}}}\newline}}

\else
\def\todo#1{}
\def\commentbox#1{}
\renewcommand{\st}[1]{}
\colorlet{red}{black!100} 
\colorlet{teal}{black!100} 
\colorlet{brown}{black!100} 
\colorlet{blue}{black!100} 
\colorlet{magenta}{black!100} 
\colorlet{purple}{black!100} 
\colorlet{cyan}{black!100} 
\fi

\title[Quasi-F-splittings in birational geometry]{Quasi-F-splittings in birational geometry 
(Quasi-F-scindages 
en géométrie birationnelle)}

\author{Tatsuro Kawakami}
\address{Department of Mathematics, Graduate School of Science, Kyoto University, Kyoto 606-8502, Japan} 
\email{tkawakami@math.kyoto-u.ac.jp}
\author{Teppei Takamatsu}
\address{Department of Mathematics, Graduate School of Science, Kyoto University, Kyoto 606-8502, Japan}
\email{teppeitakamatsu.math@gmail.com}
\author{Hiromu Tanaka} 
\address{Graduate School of Mathematical Sciences, 
The University of Tokyo, 
3-8-1 Komaba, Meguro-ku, Tokyo 153-8914, JAPAN} 
\email{tanaka@ms.u-tokyo.ac.jp}
\author{Jakub Witaszek} 
\address{Princeton University, Department of Mathematics, Fine Hall, Washington Road, Princeton NJ 08544-1000, USA}
\email{jwitaszek@princeton.edu}
\author{Fuetaro Yobuko}
\address{Graduate School of Mathematics, Nagoya University, Furocho, Chikusaku, Nagoya, Japan}
\email{soratobumusasabidesu@gmail.com}
\author{Shou Yoshikawa}
\address{Tokyo Institute of Technology, Tokyo 152-8551, Japan}
\email{yoshikawa.s.al@m.titech.ac.jp}

\begin{document}



\begin{abstract}
We develop the theory of quasi-$F$-splittings in the context of birational geometry. 
Amongst other things, we obtain results on liftability of sections and establish a criterion for whether a scheme is quasi-$F$-split employing the higher Cartier operator. As one of the applications of our theory, we prove that three-dimensional klt singularities in large characteristic are quasi-$F$-split, and so, in particular, they lift modulo $p^2$.

\bigskip 

\noindent R{\tiny \'ESUM\'E}. 
Nous d\'eveloppons la th\'eorie des quasi-$F$-scindages 
dans le 
contexte de la g\'eom\'etrie birationnelle. Entre autres, nous obtenons des r\'esultats sur le rel\`evement des sections et \'etablissons un crit\`ere pour d\'eterminer si un sch\'ema est quasi-$F$-scind\'e en utilisant l'op\'erateur de Cartier sup\'erieur. 
Comme application de notre th\'eorie, nous prouvons que les singularit\'es klt tridimensionnelles en grande caract\'eristique sont quasi-$F$-scind\'ees et, en particulier, qu'elles se rel\`event modulo $p^2$.
\end{abstract}

\maketitle

\setcounter{tocdepth}{1}
\tableofcontents

\input section1.tex

\input section2.tex

\input section3.tex

\input section4.tex

\input section5.tex

\input section6.tex

\input section7.tex

\bibliographystyle{skalpha}
\bibliography{bibliography.bib}
\end{document}

%% file: section1.tex
\section{Introduction}\label{s-intro}

What allowed for many recent developments in both commutative algebra and birational geometry (see, for example, \cite{hx13}) was the study of \emph{$F$-split} varieties, that is, varieties $X$ in positive characteristic for which the Frobenius homomorphism $F \colon \cO_X \to F_*\cO_X$ splits. The importance of this notion stems from {the fact that it interpolates} between arithmetic and geometric properties of algebraic varieties. Specifically, 
{\cora projective} 
$F$-split varieties are always \emph{{weakly} ordinary} ({\cora the} Frobenius acts bijectively on $H^{\dim X}(X, \cO_X)$) and they have non-positive Kodaira dimension. In the case of Calabi-Yau varieties, being $F$-split and ordinary is equivalent.

Unfortunately, $F$-splittings are irrelevant {to} the study of non-{\cora weakly }ordinary varieties, even though such varieties may still carry interesting arithmetic information. 
To wit, every Calabi-Yau variety $X$ can be assigned its \emph{Artin-Mazur height} $\mathrm{ht}(X) \in \bZ_{>0} \cup \{\infty\}$ which is a crystalline-type invariant measuring  supersingularity. Having height one is equivalent to ordinarity, {whilst} having finite height may be thought of as mild non-ordinarity. Elliptic curves can only have height $1$ (ordinary) or height $2$ (supersingular), but the situation is more complicated in higher dimensions. For example, if $X$ is a K3 surface, then $\mathrm{ht}(X) \in \{1, \ldots, 10\} \cup \{\infty\}$ and every such height occurs. In view of the above, it is natural to wonder if one could find a similar notion to an $F$-splitting in the context of finite Artin-Mazur heights, 
and indeed, such a notion was found and introduced by the {fifth} author in \cite{yobuko19}. {One calls} a variety \emph{$n$-quasi-$F$-split} if there exists a dashed arrow rendering the following diagram 
\begin{equation} \label{diagram:intro-definition}
\begin{tikzcd}
W_n\cO_X \arrow{r}{{F}} \arrow{d}{R^{n-1}} & F_* W_n \arrow[dashed]{ld}{\exists} \cO_X \\
\cO_X
\end{tikzcd}
\end{equation}
{\cora in the category of $W_n\mathcal{O}_X$-modules}
commutative, where $W_n\cO_X$ is a sheaf of Witt vectors. 
We will say that $X$ is \emph{quasi-$F$-split} if it is $n$-quasi-$F$-split for some $n \in {\Z_{>0}}$. 
In analogy to $F$-splittings, 
if $X$ is a Calabi-Yau variety, then {the Artin-Mazur height} $\mathrm{ht}(X)$ is equal to 
the {smallest integer} $n \in {\Z_{>0}}$ for which $X$ is $n$-quasi-$F$-split. What is striking about quasi-$F$-splittings is that, despite them being much less restrictive, they share many properties with $F$-splittings. For instance, quasi-$F$-split varieties lift modulo $p^2$ \cite{yobuko19,achinger-zdanowicz21} and satisfy Kodaira vanishing \cite{NY21}.

The main objective of our article is to generalise the theory of $F$-splittings to quasi-$F$-splittings in the context of birational geometry. In particular, as one of the applications, we show the following.

\begin{theoremIntro}[{Theorem \ref{thm:qFsplit-3dim-klt}}] \label{thm:intro-3dim-klt}
There exists a positive integer $p_0>0$ such that the following holds. 
Let $X$ be a three-dimensional {$\bQ$-factorial} 
{affine} klt 
{variety} 
over a perfect field of characteristic $p\geq p_0$. Then $X$ is quasi-$F$-split. In particular, it lifts modulo $p^2$.
\end{theoremIntro}

\noindent Moreover, assuming a natural generalisation of the work of Arvidsson-Bernasconi-Lacini on log liftability of log del Pezzo pairs with standard coefficients (Conjecture \ref{conj:log-liftability}), we show that $p_0=43$ works in the above theorem (see Remark \ref{rem:main-theorem-for-42}). 
Note that  {Theorem \ref{thm:intro-3dim-klt}} is false for $F$-splittings as shown in {\cite[Theorem 1.1]{CTW15a}}. 
We also point out that it is usually very hard to prove liftability of singularities as their deformation theory is governed by the cotangent complex which is rarely computable in practice.

We introduce two techniques for verifying whether a {variety} is quasi-$F$-split: inversion of adjunction and a criterion employing the higher Cartier operator. By using inversion of adjunction and the local-to-global correspondence between klt singularities and log Fano varieties, one can reduce Theorem \ref{thm:intro-3dim-klt} to the study of {whether} log del Pezzo surfaces are quasi-$F$-split. 
In this context, we can show the following.
\begin{theoremIntro}[{Corollary \ref{cor:del-Pezo-type-surfaces-are-quasi-F-split}, {cf.\ Theorem \ref{thm:qFsplit-del-Pezzo-in-large-characteristic}}}] \label{thm:intro-quasi-F-split-del-Pezzo}\label{Introthm B}
Let $X$ be a klt del Pezzo surface over a perfect field of characteristic $p>5$. Then $X$ is quasi-$F$-split.
\end{theoremIntro}

\noindent In fact, we prove that if $(X,\Delta)$ is a log del Pezzo pair in characteristic ${p \gg 0}$ 
with $\Delta$ having standard coefficients (that is of the form $1-\frac{1}{n}$ for $n \in {\Z_{>0}}$), then $(X,\Delta)$ is $2$-quasi-$F$-split {(Theorem \ref{thm:qFsplit-del-Pezzo-in-large-characteristic})}. {These results are based on \cite{CTW15b} and \cite{ABL}.}

One of the foundational results in the theory of $F$-splittings is the theorem of Hara \cite{hara98}, who proved that {affine} klt surfaces are $F$-split in characteristic $p>5$. We show a stronger result for quasi-$F$-splittings. 

\begin{theoremIntro}[Corollary \ref{c-QFS-2dim-rel}]\label{Introthm C}
Let $X$ be an affine klt surface
over a perfect field of characteristic $p>0$. Then $X$ is quasi-$F$-split.
\end{theoremIntro}

\noindent In fact, we prove that every two-dimensional {affine} klt pair $(X,\Delta)$ in every positive characteristic is quasi-$F$-split. 
Surprisingly, one need not even assume that $\Delta$ has standard coefficients. 

\begin{remark}
    {\cora In Theorems \ref{thm:intro-3dim-klt}, \ref{Introthm B}, and \ref{Introthm C}, the choice of splittings are not natural amongst  
 other varieties. Indeed, these theorems depend on Theorem \ref{thm:intro-higher-Cartier-criterion-for-quasi-F-split} (cf.~Theorem \ref{thm:higher-Cartier-criterion-for-quasi-F-split}), which ensures only the existence of splittings. }
\end{remark}

Before moving on to discussing the aforementioned two key concepts (Subsection \ref{ss:intro-lifting-sections} and Subsection \ref{ss:intro-higher-Cartier}), we point out that multiple {other}  foundational results on quasi-$F$-splittings are contained in this article, for example, invariance under birational and finite maps (Corollary \ref{cor:pulling-back-quasi-F-splittings} and Corollary \ref{prop:quasi-F-splittings-under-finite-maps}), non-positivity of Kodaira dimension (Proposition \ref{prop:non-negative-Kodaira-dimension-for-quasi-F-split-pairs}), 
and Kawamata-Viehweg vanishing for quasi-$F$-split pairs (Theorem \ref{thm:KVV-for-quasi-F-split-pairs}). Furthermore, one of the starting problems of our work is finding a correct definition of a quasi-$F$-splitting for log pairs $(X,\Delta)$. 
To this end, we build on the construction, due to the {third} author, of Witt vectors with $\bQ$-boundary (see \cite{tanaka22}).

\subsection{Lifting sections and inversion of adjunction} \label{ss:intro-lifting-sections}
Inversion of adjunction is a key tool in the study of $F$-splittings. Briefly speaking, in the local setting, it says that if $X$ is a normal variety and {\cora $S \subseteq X$} 
is an $F$-pure normal prime divisor 
such that {\cora $K_X+S$} is Cartier, then $X$ {is $F$-pure along {\cora $S$}} \cite[main Theorem]{Sch09}. In the global setting, it stipulates that if $X$ is projective 
over a field of positive characteristic, {\cora $-(K_X+S)$} is ample, and $S$ is $F$-split, then $X$ is $F$-split, too.

Unfortunately, local inversion of adjunction is false for quasi-$F$-splittings (see \cite[Example 7.7]{KTY22}). 
What is true however is that global quasi-$F$-split inversion of adjunction holds when $-S$ is semiample. {\cora Typical} situations when this is the case is when $X$ admits a projective morphism $\pi \colon X \to Z$ such that either
$\pi$ is birational and $-S$ is ample, or such that $Z$ is a curve and $S$ is a fibre of $\pi$.

\begin{theoremIntro}[Corollary \ref{cor:inversion-of-adjunction-anti-ample}] \label{thm:intro-inversion-of-adjunction}
Let $X$ be a smooth variety admitting a projective  morphism $\pi \colon X \to Z$ to an affine normal variety $Z$ over a perfect field of characteristic $p>0$. Let $S$ be a smooth prime divisor on $X$. Assume that $S$ is quasi-$F$-split, 
$-S$ is {semiample}, $A := -(K_X+S)$ is ample, and 
\[
H^1(X, \cO_X(K_X+p^iA-kS))=0
\]
for every $i\geq 1$ and $k\geq 0$. Then $X$ is quasi-$F$-split {over an open neighbourhood of $\pi(S)$}. 
\end{theoremIntro}

\noindent Our result is more general; instead of smoothness, it is enough to assume that 
{\cora $X$ is  $F$-pure and $\MO_X(D)$ is Cohen-Macaulay for every Weil divisor $D$ on $X$} 
 (see  Corollary \ref{cor:inversion-of-adjunction-anti-ample}). We also point out that the assumption on the vanishing of $H^1$ is not very restrictive; in most cases, it can be dropped by methods of \cite{HW17}.\\ 

The key to showing the above theorem is lifting sections from quasi-$F$-split varieties. In the case of $F$-splittings, this works as follows \cite{schwede14}. Given a {Cartier divisor $L$} on a variety $X$, 
one defines {linear} subspaces $S^0(X, S; \cL) \subseteq H^0(X,\cL)$ and $S^0(S, \cL|_S) \subseteq H^0(S,\cL|_S)$ for $\cL := \cO_X(L)$ such that
\[
S^0(X, S; \cL) \to S^0(S, \cL|_S)
\]
is surjective, whenever $X$ is projective, $S$ is a smooth Cartier divisor, and $L-(K_X+S)$ is ample. When $S$ is $F$-split, the right hand side equals $H^0(S, \cL|_S)$, which implies that the restriction map $H^0(X,\cL) \to H^0(S, \cL|_S)$ is surjective. 

Analogously, we define linear subspaces $qS^0_n(X, \cL)$ and $qS^0_{n,\adj}(X,S;\cL)$  of  $H^0(X,\cL)$ such that the following holds.

\begin{theoremIntro}[{Theorem \ref{thm:lifting-sections-basic}}] \label{thm:intro-quasi-F-split-lifting}
Let $X$ be a smooth projective variety 
over a perfect field $k$ of characteristic $p>0$. Let $S$ be a smooth prime divisor on $X$ and let $L$ be a Cartier divisor {on $X$}. 
Set $A:= L-(K_X+S)$. Assume that
\begin{enumerate}
    \item $H^1(X, \cO_X(K_X + p^iA))=0$ and 
    \item $qS^0_{n-i}(X, \cO_X(K_X + p^iA)) = H^0(X, \cO_X(K_X+p^iA))$
\end{enumerate}
for all $i \geq 1$. Then
\[
qS^0_{n,\adj}(X,S; \cO_X(L)) \to qS^0_{n}(S, \cO_S(L|_S))
\]
is surjective.  
\end{theoremIntro}
\noindent In contrast to $F$-splittings and $S^0$, it is necessary to distinguish between the adjoint and non-adjoint variants of $qS^0_n$. 
{\cora For example, the following are equivalent for a normal variety $X$  and an effective $\Q$-divisor $\Delta$ (cf.\ Subsection \ref{ss log F-split}): 
\begin{enumerate}
\item $(X, \Delta)$ is naively keenly $F$-split, i.e., $F: \MO_X \to F_*\MO_X(p\Delta)$ splits. 
\item $F: \MO_X \to F_*\MO_X(p\Delta)$ splits and $\rdown{\Delta}=0$. 
\item $(X, \Delta)$ is $1$-quasi-$F$-split. 
\end{enumerate}
}

\subsection{Criterion using higher Cartier operator}
\label{ss:intro-higher-Cartier}

For any $i \geq 0$, {we} define 
\[ 
B_1\Omega^i_X := \mathrm{Im}(F_*d \colon F_*\Omega_X^{i-1} \to F_*\Omega_X^i).
\]
Such an object features in the Cartier isomorphism $C^{-1} \colon \Omega^i_X \xrightarrow{\cong} \cH^i(F_* \Omega^{\mydot}_X)$, which is one of the most powerful tools in positive characteristic algebraic geometry. Having {introduced}  this definition, we can now state our {aforementioned} criterion.

\begin{theoremIntro}[{Theorem \ref{thm:higher-Cartier-criterion-for-quasi-F-split}}] \label{thm:intro-higher-Cartier-criterion-for-quasi-F-split} Let $X$ be a $d$-dimensional smooth projective Fano variety over a perfect field of characteristic $p>0$. Suppose that
\begin{enumerate}
    \item  $H^{d-2}(X,\Omega^1_X \otimes \omega_X) = 0$ and
    \item $H^{d-2}(X,B_1\Omega^2_X \otimes \omega^{p^k}_X) = 0$ for every $k \geq 0$.
\end{enumerate}
Then $X$ is quasi-$F$-split.
\end{theoremIntro}
\begin{remark}
    Since $H^2(X,T_X)\cong H^{d-2}(X,\Omega^1_X \otimes \omega_X)$, the assumption (1) implies $X$ lifts to $W(k)$ \cite[Theorem 8.5.19]{fga2005}.
\end{remark}
\noindent In fact, by an application of the usual Cartier operator, 
one can check that Assumption (2) is valid if $H^{d-i-1}(X,\Omega^i_X \otimes  \omega^{p^k}_X)=H^{d-i}(X,\Omega^i_X \otimes  \omega^{p^{k+1}}_X) =0$ for every $k \geq 0$ and $i \geq 2$ (see Remark \ref{remark:log-criterion-no-B}).

As before, our result is much more general and applies to arbitrary relative log Fano pairs over affine bases. We point out that Assumption (1) is a special case of Akizuki-Nakano vanishing which holds in positive characteristic when $X$ lifts modulo $p^2$. On the other hand, Assumption (2) is much more subtle. Fortunately, when $d=2$, it is enough to verify that $H^0(X, \Omega^2_X \otimes  \omega^{p^{k+1}}_X) = 0$ which often follows for geometric reasons. 

The proof of the above theorem is based on the calculation of $H^{d-1}(X, B_n\Omega^1_X \otimes \omega_X)$ via the higher Cartier operator.

\subsection{Further directions}
\label{ss:intro-further-directions}
There are multiple new directions of research that  follow from, or are related to, the content of our paper. 

We summarise a few of such topics below. The list of natural questions at the intersection of $F$-splittings, quasi-$F$-splittings, and birational geometry is much longer, and we hope to gradually address them in the coming years.

\subsubsection*{Explicit bound and applications to the logarithmic extension of differential forms}
In {\cora the follow-up} article \cite{KTTWYY2} building on all the results of this paper and careful geometric considerations, we show that Theorem \ref{thm:intro-3dim-klt} holds for $p_0=43$ unconditionally (without assuming Conjecture \ref{conj:log-liftability}). Moreover, we prove that this bound on  $p$ is sharp, that is, there exists a three-dimensional klt singularity in characteristic $p=41$ which is not quasi-$F$-split.

On top of that, we extend some of the results of \cite{Kawakami22} on logarithmic extensions of one-forms to quasi-$F$-split singularities. As a consequence, we obtain the following theorem.
\begin{theoremIntro}[{\cora \textup{\cite[Theorem E]{KTTWYY2}}}]
Let $X$ be a three-dimensional variety over a perfect field of characteristic {\cora $p\geq 43$}. 
Assume that one of the following holds. 
\begin{enumerate}
\item $X$ is terminal. 
\item 
The singular locus of $X$ is zero-dimensional, and $X$ is $\Q$-factorial and klt. 
\end{enumerate}
Then $X$ satisfies the logarithmic extension theorem for one-forms.
\end{theoremIntro}

\subsubsection*{Minimal model program}
The fact that klt surface singularities are $F$-split in general only for $p>5$ is the key obstruction for generalising the three-dimensional Minimal model program of Hacon and Xu to low characteristics. We hope that such a generalisation can be now obtained by replacing the use of $F$-splittings with quasi-$F$-splittings in \cite{hx13}. However, given that the  quasi-$F$-split lifting result (Theorem \ref{thm:intro-quasi-F-split-lifting}) is not as strong as its analogue for $F$-splittings, this calls for a far more careful handling of the construction of flips, and so is beyond the scope of our article.

Similarly, Theorem \ref{thm:intro-3dim-klt} and Theorem \ref{thm:intro-quasi-F-split-del-Pezzo} should play an important role in establishing the four-dimensional Minimal model program in full generality for {\cora $p \geq 43$}.

\subsubsection*{Fedder's criterion} 
One of the advantages of the theory of $F$-split singularities (as opposed to, say, klt singularities) is that there is an explicit criterion (called \emph{Fedder's criterion}) for verifying whether a system of polynomial equations is $F$-split. An analogue of this result for quasi-$F$-splittings is obtained in an independent work of the first, {second}, and sixth authors \cite{KTY22}. Most notably, amongst many applications, the work contains the proof of the fact that there exist Calabi-Yau varieties of arbitrarily high Artin-Mazur height.

\subsubsection*{Quasi-$F$-regularity}  
When dealing with $F$-splittings, one often considers not a single Frobenius, but a multiple ${F^e} \colon \cO_X \to F^e_*\cO_X$ thereof as this allows for a direct use of Serre's vanishing in proofs. Unfortunately, naively replacing $F$ by $F^e$ in Diagram \ref{diagram:intro-definition} yields a meaningless definition. Nevertheless, the third, fourth, and fifth authors {\cora \cite{TWY24}}
find a way to rectify that, which, in turn, allows for defining quasi-analogues of $F$-regularity and $+$-regularity. Therewith, one can weaken some of the assumptions of the above theorems (unfortunately, the condition on $qS^0_{n-i}$ in Theorem \ref{thm:intro-quasi-F-split-lifting}, which is the most restrictive in terms of applications, cannot be dropped).  

\subsubsection*{Fano varieties}
{\cora As mentioned above, quasi-$F$-split varieties satisfy Kodaira vanishing. This fact enables us to show Kodaira vanishing for a much broader class of varieties in positive characteristic.

For example, the first and third authors \cite{Kawakami-Tanaka(Fano3-fold)} proved that smooth Fano threefolds are all quasi-$F$-split if the Picard rank or the Fano index is greater than one. As an application, we clarified that Kodaira vanishing holds for all smooth Fano threefolds (without assumptions on Picard number or Fano index).

Since there exist non-$F$-split smooth Fano threefolds whose Picard rank or Fano index is greater than one in each characteristic $p=2$, $3$, and $5$, 
we essentially need the quasi-$F$-splitting of Fano threefolds in the above argument. For this reason, we believe that the theory of quasi-$F$-splitting can be more useful than the theory of $F$-splitting.

We also show that smooth del Pezzo varieties in any dimension are all quasi-$F$-split, though they are not necessarily  $F$-split in each characteristic $p=2$, $3$, and $5$. We refer to \cite{Kawakami-Tanaka(delPezzoVar)} for details.}

\medskip
\noindent {\bf Overview.}
As our work builds on many technical results, we provided additional foundational material to make it easier for the reader to understand the content of this article. In particular, Subsections \ref{ss:witt-vectors}, \ref{ss:quasi-F-splittings}, \ref{ss:iterated-Cartier-operator}, \ref{ss:intro-to-lifting-quasi-stable-sections}, \ref{ss:sketch-of-proof-of-higher-Cartier-criterion}, \ref{ss:calculation-local-cohomology}, \ref{ss:canonical-lift} are mostly self-contained and cover topics such as Witt vectors, basics of quasi-$F$-splittings, higher Cartier operator,  quasi-$F$-stable sections, explicit calculations, and canonical liftings of quasi-$F$-split varieties modulo $p^2$. These subsections do not utilise log pairs; in particular, they do not require any understanding of log $F$-splittings or Witt divisorial sheaves.  The reader might consider exploring these subsections first before delving into the technical aspects of our article.

\medskip
\noindent {\bf Acknowledgements.}
The authors thank Mircea {\mustata} for valuable conversations related to the content of the paper.
They are also grateful to Shunsuke Takagi for helpful comments.
{\cora They thank the anonymous referees for many valuable comments that have improved the article.}
\begin{itemize}
    \item Kawakami was supported by JSPS KAKENHI Grant number JP22KJ1771 and JP24K16897.
    \item Takamatsu was supported by JSPS KAKENHI Grant number JP22J00962.
    \item Tanaka was supported by JSPS KAKENHI Grant numbers JP18K13386, 
     JP22H01112, and JP23K03028.
    \item Witaszek was supported by NSF research grant DMS-2101897.
    \item Yobuko was supported by JSPS KAKENHI Grant number JP19K14501.
    \item Yoshikawa was supported by JSPS KAKENHI Grant number JP20J11886. 
\end{itemize}

%% file: section2.tex
\section{Preliminaries}\label{s-Prelim}

\subsection{Notation}
\label{ss:notation}

In this subsection, we summarise notation and basic definitions used in this {article}. 
\begin{enumerate}
\item Throughout the paper, $p$ denotes a prime number and $\F_p \coloneqq \Z/p\Z$. 
\item 
We say that a scheme $X$ is {\em of characteristic $p$} if $X$ is an $\F_p$-scheme, 
that is, the induced morphism $X \to \Spec \Z$ factors through the closed immersion $\Spec \F_p \hookrightarrow \Spec \Z$. We denote by $F \colon X \to X$ the absolute Frobenius morphism on an $\F_p$-scheme $X$.
\item For an integral scheme $X$, 
we define the {\em function field} $K(X)$ of $X$ 
as the stalk $\MO_{X, \xi}$ at the generic point $\xi$ of $X$. 

\item An effective Cartier divisor $D \subseteq X$ on a Noetherian scheme $X$ is called \emph{simple normal crossing} if for every $x \in D$, the local ring $\cO_{X,x}$ is regular and there exists a regular system of parameters $x_1,\ldots, x_d$ in the maximal ideal $\m$ of $\cO_{X,x}$ and $1 \leq r \leq d$ such that $D$ is defined by $x_1\cdots x_r$ in $\cO_{X,x}$ (cf.\ \cite[\href{https://stacks.math.columbia.edu/tag/0BI9}{Tag 0BI9}]{stacks-project} and \cite[\href{https://stacks.math.columbia.edu/tag/0BIA}{Tag 0BIA}]{stacks-project}). 
For the definition of effective Cartier divisors, we refer to \cite[Chapter II, Section 6]{hartshorne77}. 
\item {Given an integral normal Noetherian scheme $X$, a projective birational morphism $\pi \colon Y \to X$ is called \emph{a log resolution (of singularities) of $X$} if $Y$ is regular and $\mathrm{Exc}(f)$ is a simple normal crossing divisor.}
\item Given a coherent sheaf $\cF$ on a scheme $X$, we define $\cF^{\vee} \coloneqq \cHom_{\cO_X}(\cF, \cO_X)$. We warn the reader that $(-)^{\vee}$ also denotes Matlis duality (see Subsection \ref{ss:Matlis-duality}) but it {should} always {be} clear from the context which one we refer to.

\item A coherent sheaf $\cF$ on an integral normal Noetherian scheme $X$ is \emph{divisorial} if $\cF \cong \cO_X(D)$ for some Weil divisor $D$.

    \item We say that $X$ is a {\em variety} (over a field $k$) if $X$ is an integral scheme 
    that is separated and of finite type over $k$. We say that $X$ is a {\em curve} if $X$ is a variety of dimension one, and that $X$ is a {\em surface} if $X$ is a variety of dimension two.
\item 
Let $T$ be a Noetherian scheme. 
We say that $(X, D)$ is {\em log smooth over} $T$ if $X$ is a smooth scheme over $T$, 
$D$ is an effective Cartier divisor on $X$, 
and for any $x \in X$, there exist an open neighbourhood $X'$ of $x \in X$ and an \'etale morphism 
$g \colon X' \to \mathbb A^d_T = T \times_{\Spec \Z} \Spec \Z[x_1, ..., x_d]$ such that 
the equality $D|_{X'} = g^{-1}E$ of {closed subschemes} holds, where $E = V(x_1 \cdots x_r)$ for some $1 \leq r \leq d$. 
{\cora Note that if $(X, D)$ is log smooth over $T$ and 
$T' \to T$ is a morphism of Noetherian schemes, 
then the base change $(X \times_T T', D \times_T T')$ is log smooth over $T'$ \cite[the paragraph immediately after Definition 9.3.6]{fga2005}.}
\item 
{Let $k$ be a perfect field $k$ of characteristic $p>0$ and let $(R,m)$ be a Noetherian local ring with residue field $k$. 
We say that a normal variety $X$ over $k$ \emph{lifts to $R$} if there exists a flat separated morphism $\cX \to \Spec R$ of finite type such that 
$ \cX \times_{\Spec R} \Spec k$ is isomorphic to $X$ over $k$. 

Let $(X, D)$ be a log smooth projective pair over $k$. 
Write $D = \sum_{i=1}^r D_i$, where $D_i$ are prime divisors. 
We say that $(X,D)$ \emph{lifts to $R$} 
if there exist a log smooth projective pair $(\cX, \cD)$ over $R$, 
effective Cartier divisors $\cD_1, ..., \cD_r$ on $\cX$ flat over $R$, and an isomorphism $\phi \colon \cX \times_{\Spec R} \Spec k \xrightarrow{\cong} X$ over $k$ 
such that the equality $\cD = \cD_1 + \cdots + \cD_r$ of Cartier divisors holds and 
$\phi( \cD_i \times_{\Spec R} \Spec k) =D_i$ for every $1 \leq i \leq r$. 
}
\item We say that a scheme $X$ is {\em excellent} if it is Noetherian and all the stalks 
$\MO_{X, x}$ are excellent. We note that the regular locus of an integral normal excellent scheme $X$ 
is an open dense subset of $X$ which contains all the points of codimension one
(see \cite[\href{https://stacks.math.columbia.edu/tag/07P7}{Tag 07P7}]{stacks-project} and \cite[\href{https://stacks.math.columbia.edu/tag/0BX2}{Tag 0BX2}]{stacks-project}).

    \item \label{sss-Ffinite} {We say that an $\F_p$-scheme $X$ is {\em $F$-finite} if 
$F \colon X \to X$ is a finite morphism,  
 and  we say that an $\F_p$-algebra $R$ is {\em $F$-finite} if $\Spec R$ is $F$-finite. Such schemes admit many good properties.

\begin{enumerate}
\renewcommand{\labelenumi}{(\roman{enumi})}
    \item If $R$ is an $F$-finite Noetherian $\F_p$-algebra, 
    then it is a homomorphic image of a regular ring of finite Krull dimension \cite[Remark 13.6]{Gab04}; 
    in particular, $R$ is excellent, it admits a dualising complex, and $\dim R < \infty$.  
    \item $F$-finite rings are stable under localisation and ideal-adic completions \cite[Example 9]{Has15}.
    \item If a scheme $X$ is of finite type over an $F$-finite Noetherian $\F_p$-scheme $Y$, then it is also $F$-finite. 
\end{enumerate}}

\item \label{it:essential}
{We recall the basics of morphisms essentially of finite type and refer to \cite{Nay09} for details. 
Let $\varphi \colon A \to B$ be a ring homomorphism and let $f \colon X \to Y$ be a morphism of Noetherian schemes.
\begin{enumerate} 
    \item We say that $\varphi$ is {\em essentially of finite type}  if $B$ is isomorphic to $S^{-1}C$ as an $A$-algebra for some finitely generated $A$-algebra $C$ and a multiplicatively closed subset $S \subseteq C$. 
    \item We say that $\varphi$ is {\em localising} 
    if $B$ is isomorphic to $S^{-1}A$ as an $A$-algebra for some multiplicatively closed subset $S \subseteq A$. 
    \item We say that $f$ is {\em essentially of finite type} or {\em localising} if it is such affine-locally 
    (explicitly, for any $x \in X$,  there exists an affine open neighbourhood $U \subseteq X$ of $x$ which maps into an affine open $V \subseteq Y$ such that the induced {ring homomorphism} $\cO_Y(V) \to \cO_X(U)$ satisfies the required property, cf.\ \cite[\href{https://stacks.math.columbia.edu/tag/01SS}{Tag 01SS}]{stacks-project}). We say that $f$ is a {\em localising immersion} if it is an injective localising morphism. 
    \item A composition of morphisms essentially of finite type between Noetherian schemes is essentially of finite type.
    \item If $f$ is separated and essentially of finite type, then there exists a factorisation $f \colon X \xrightarrow{j} X' \xrightarrow{f'} Y$
    such that $j$ is a localising immersion and $f'$ is separated and of finite type \cite[Theorem 3.6]{Nay09}.
\end{enumerate}}

\item 
Given an integral normal Noetherian scheme $X$ and a $\bQ$-divisor $D$, 
we define the subsheaf $\MO_X(D)$ of the constant sheaf $K(X)$ on $X$ 
{by the following formula}
\[
\Gamma(U, \MO_X(D)) = 
\{ \varphi \in K(X) \mid 
\left({\rm div}(\varphi)+D\right)|_U \geq 0\}
\]
{for every} open subset $U$ of $X$. 
In particular, 
$\cO_X(\rdown{{D}}) = \cO_X({D})$.
\end{enumerate}

\begin{remark} \label{remark:hom-reflexive}
In this article, we will repeatedly use the fact that $\cHom_{\cO_X}(\cF, \cG)$ is reflexive, when $\cF$ and $\cG$ are coherent sheaves on an integral Noetherian scheme $X$, and $\cG$ is reflexive (note that when $\Supp \cG=X$ {and} $X$ is normal, $\cG$ is $S_2$ if and only if it is reflexive {\cora \cite[Proposition 1.1]{Langer24b}}).
In particular, if $j \colon U \to X$ is the inclusion of the regular locus $U$ of a normal excellent scheme $X$, then $j_*\cHom_{\cO_U}(\cF|_U, \cG|_U) = \cHom_{\cO_X}(\cF, \cG)$. 
\end{remark}

Last, we emphasise the following two non-standard definitions used in  the article.
\begin{definition}
We say that an integral normal Noetherian scheme $X$ is \emph{divisorially Cohen-Macaulay} if every divisorial sheaf on $X$ is Cohen-Macaulay. 
\end{definition}
\noindent In particular, divisorially Cohen-Macaulay schemes are Cohen-Macaulay. 
\begin{remark} \label{remark:divisorially-Cohen-Macaulay}
{Note} that $X$ is divisorially Cohen-Macaulay if {$\dim X=2$} or 
$X$ is a strongly $F$-regular Noetherian  $F$-finite $\Q$-factorial  scheme  (\cite[Corollary 3.3]{patakfalvi-schwede14}). 
\end{remark}



\begin{definition} \label{def:p-compatible-divisors} 
For a prime number $p$ and a normal Noetherian scheme $X$, we say
\begin{enumerate}
    \item that a $\bQ$-divisor $\Delta$ on $X$ is \emph{$p$-compatible} if $\{\Delta\} \geq \{p^i\Delta\}$ for every $i \in \Z_{>0}$,
    \item {that a log pair $(X,\Delta)$ is \emph{$p$-compatible} if $\Delta$ is $p$-compatible and $K_X + \{p^i\Delta\}$ is $\bQ$-Cartier for every $i \in \Z_{\geq 0}$.}
\end{enumerate} 
\end{definition}
\noindent The two key cases in which a $\bQ$-divisor $\Delta$ is $p$-compatible is when it has standard coefficients or its  coefficients are of the form $\frac{a}{p}$ for $a \in \bZ$.

\subsubsection{Dualising complexes and canonical divisors}

For basic properties {of} dualising complexes, we refer to \cite[\href{https://stacks.math.columbia.edu/tag/0A85}{Tag 0A85}]{stacks-project} and \cite{Har66}. 

Throughout the article, whenever we consider a dualising complex, we implicitly work over a fixed excellent base scheme, which is denoted by $B$ in this subsection.
We assume that $B$ admits a dualising complex and we implicitly make a choice of {one}, say $\omega^{\mydot}_B$. It {automatically} holds that $\dim B < \infty$ \cite[Corollary V.7.2]{Har66}.

Given an integral normal excellent scheme $X$ and a separated morphism $\pi \colon X \to B$ of finite type, we set $\omega^{\mydot}_X := \pi^! \omega^{\mydot}_B$. Then $\omega^{\mydot}_X$ is a dualising complex of $X$ (cf.\ \cite[\href{https://stacks.math.columbia.edu/tag/0AU5}{Tag 0AU5}]{stacks-project}). In what follows, if a base scheme is not specified, then we implicitly take $X=B$ and $\pi = \mathrm{id}$. 

There exists a unique $e \in \Z$ such that 
$\mathcal H^{-e}(\omega^{\mydot}_X) \neq 0$ and 
$\mathcal H^{i}(\omega^{\mydot}_X) = 0$ for $i< -e$. We say that $\omega^{\mydot}_X$ is \emph{normalised} if $e = \dim X$. We set $\omega_{X} := \mathcal H^{-e}(\omega^{\mydot}_X)$. The sheaf $\omega_X$ is called a \emph{dualising sheaf} of $X$, 
and it is known that  
it is reflexive \cite[\href{https://stacks.math.columbia.edu/tag/0AWE}{Tag 0AWE}]{stacks-project} and invertible on the regular locus of $X$ \cite[\href{https://stacks.math.columbia.edu/tag/0AWX}{Tag 0AWX}]{stacks-project}. 
We fix a Weil divisor $K_X$ such that $\cO_X(K_X) \cong \omega_X$. Any such Weil divisor is called a {\em canonical divisor}\footnote{In our situation, we always work with a fixed $\omega^{\mydot}_B$, but if one allowed $\omega^{\mydot}_B$ to vary, then a canonical divisor $K_X$ on $X$ would only be defined up to linear equivalence over $B$, that is, 
if $K_X$ and $K'_X$ are canonical divisors on $X$ as above, 
then there exists an invertible sheaf $L$ on $B$ such that $\MO_X(K_X) \cong \MO_X(K'_X) \otimes \pi^*L$ \cite[Theorem V.3.1]{Har66}.} on $X$.

Throughout the article, we always assume that the dualising complex $\omega^{\mydot}_B$ on the base scheme $B$ is normalised.

\begin{proposition}[{\cora cf.\ \cite[\href{https://stacks.math.columbia.edu/tag/0AWE}{Tag 0AWE}]{stacks-project} and \cite[\href{https://stacks.math.columbia.edu/tag/0AWI}{Tag 0AWI}]{stacks-project}}]
{\cora
Let $\pi \colon X \rightarrow B$ be a proper morphism of integral excellent schemes such that $\pi$ is surjective or $\dim B = \dim \MO_{B, b}$ for every closed point $b \in B$.
Assume that the dualising complex $\omega^{\mydot}_B$ is normalized as above.
Then the dualizing complex $\omega^{\mydot}_X = \pi^!\omega_B^{\mydot}$ is also  normalised.
}
\end{proposition}

\begin{proof}
{\cora
Before proceeding with the proof, we note that if $x \in X$ is a closed point, then $b := \pi(x) \in B$ is also a closed point as $\pi$ is proper. In particular, $x$ is also a closed point of the fibre $X_b$, and so the field extension $\kappa(b) \subseteq \kappa(x)$ is finite by \cite[\href{https://stacks.math.columbia.edu/tag/01TF}{Tag 01TF}]{stacks-project}. 

First, we show that $\cH^{-\dim X}(\omega^{\mydot}_X) \neq 0$. To this end, it is enough to prove that $\cH^{-\dim X}(\omega^{\mydot}_{X,x}) \neq 0$ for any closed point $x \in X$ such that $\dim X = \dim \cO_{X,x}$. Set $b := \pi(x)$. 
If $\pi$ is surjective, by the dimension formula \cite[\href{https://stacks.math.columbia.edu/tag/02JU}{Tag 02JU}]{stacks-project} and \cite[\href{https://stacks.math.columbia.edu/tag/02JX}{Tag 02JX}]{stacks-project}, we have $\dim \cO_{B,b} = \dim B$.
Therefore, we may assume that $\dim \cO_{B,b} = \dim B$.
In this case, $\omega^{\mydot}_{B,b}$ is normalised as  a dualising complex on $\Spec \cO_{B,b}$. Therefore, by \cite[\href{https://stacks.math.columbia.edu/tag/0AWL}{Tag 0AWL}]{stacks-project}, $\omega^{\mydot}_{X,x}$ is normalised as dualising complex on $\Spec \cO_{X,x}$, and we have $\cH^{-\dim X}(\omega^{\mydot}_{X,x}) \neq 0$.

Second, we show that $\cH^{i}(\omega^{\mydot}_X)=0$ for $i < -\dim X$. To this end, since $\Supp \cH^{i}(\omega^{\mydot}_X)$ is closed and $X$ is Noetherian, it is enough to prove that $\cH^{i}(\omega^{\mydot}_{X,x})=0$ for every closed point $x \in X$ and $i < -\dim X$ 
(cf.\ \cite[\href{https://stacks.math.columbia.edu/tag/02IL}{Tag 02IL}]{stacks-project}). Set $b := \pi(x)$. By definition, $\omega_{B,b}^{\mydot}[\dim \cO_{B,b} - \dim B]$ is normalised, 
and so $\omega^{\mydot}_{X,x}[\dim \cO_{B,b} - \dim B]$ is normalised as well by \cite[\href{https://stacks.math.columbia.edu/tag/0AWL}{Tag 0AWL}]{stacks-project}. 
Thus $\cH^{i}(\omega^{\mydot}_{X,x})=0$ for all 
\[
i < -\dim \cO_{X,x} + \dim \cO_{B,b} - \dim B.
\]
If $\pi$ is surjective, then by the dimension formula \cite[\href{https://stacks.math.columbia.edu/tag/02JU}{Tag 02JU}]{stacks-project} and \cite[\href{https://stacks.math.columbia.edu/tag/02JX}{Tag 02JX}]{stacks-project} again, the right-hand side is equal to $- \dim X$.
On the other hand, if $\dim B = \dim \cO_{B,b}$, then the right-hand side is greater than or equal to $- \dim X$.
Therefore, we have $\cH^{i}(\omega^{\mydot}_X)=0$ for $i < -\dim X$ in both cases.
}
\end{proof}

\subsubsection{Singularities of minimal model program}\label{ss-mmp-sing}
We will freely use the standard notation in birational geometry, for which we refer to \cite{kollar13} and \cite{KM98}. In particular, we will use the abbreviated names for singularities such as klt or lc.  {Moreover, we repeatedly use the results on the minimal model program for surfaces as developed in  \cite{tanaka12,tanaka16_excellent}.}

We summarise some key definitions in birational geometry below.
\begin{enumerate}
    \item We say that $(X,\Delta)$ is a {\em log pair} if $X$ is an integral  normal excellent scheme admitting a dualising complex and $\Delta$ is an effective $\bQ$-divisor such that $K_X+\Delta$ is $\bQ$-Cartier\footnote{Note that some authors do not assume the last condition; moreover, in contrast to \cite[Definition 2.8]{kollar13}, we always require that $\Delta$ is effective, and so a klt pair $(X,\Delta)$ is always a log pair.}.
    \item A $\bQ$-divisor $\Delta$ has {\em standard coefficients} if 
    all the coefficients of $\Delta$ are contained in the set $\{1 - \frac{1}{n} \mid n \in\Z_{>0} \} \cup \{1\}$.
    \item A pair $(X, \Delta)$ is {\em log del Pezzo} if $(X, \Delta)$ is a two-dimensional klt pair projective over a field $k$ and such that $-(K_X+\Delta)$ is ample. 
    \item {A scheme 
$X$ is a surface {\em of del Pezzo type} over $k$ 
if $X$ is a normal projective surface over a field $k$ such that $(X, \Delta)$ is log del Pezzo for some effective $\Q$-divisor $\Delta$ on $X$.}
\end{enumerate}

\begin{remark} 
 In \cite{kollar13}, it is assumed that all schemes $X$ are of finite type over a regular base scheme $B$. This assumption was only introduced so that $B$, and so $X$,  admit a dualising complex, but the results of \cite{kollar13} go through for general excellent schemes admitting dualising complexes. 
We {point out} that 
the discrepancy $a(E, X, \Delta)$ for a log pair $(X,\Delta)$ and a prime divisor $E$ over $X$ does not depend on the choice of the dualising complex $\omega^{\mydot}_X$.


Analogous assumptions to that in \cite{kollar13} are stated in \cite{tanaka16_excellent} which covers the minimal model program for two-dimensional excellent schemes. As above, the results of \cite{tanaka16_excellent} work in our more general setting. 
\end{remark}

\subsection{Witt vectors} \label{ss:witt-vectors}
For the convenience of the reader, 
we give a very quick introduction to Witt vectors and 
refer to \cite[Ch. 0, Section 1]{illusie_de_rham_witt} and \cite[Ch. II, \S 6]{Ser79} for more details.

Given an $\F_p$-algebra $A$ 
and $n\in \Z_{>0}$,  
one can define the ring of Witt vectors $W_nA$ of length $n$. 
We have the following {equality} 
\[
W_nA = \{(a_0, \ldots, a_{n-1}) \mid a_i \in R\}
\]
of sets, 
but the addition and multiplication are not componentwise. For example, when $n=2$, 
one defines
\begin{align*}
(a_0,a_1) + (b_0,b_1) &:= \Big(a_0 + b_0, a_1+b_1 - \frac{(a_0+b_0)^p - a_0^p - b_0^p}{p}\Big), \text{ and }\\
(a_0,a_1) \cdot (b_0,b_1) &:= (a_0b_0, a_0^pb_1 + b_0^pa_1), 
\end{align*}
where $\frac{(a_0+b_0)^p - a_0^p - b_0^p}{p}$ is a short way of writing $\sum^{p-1}_{i=1}\frac{(p-1)!}{i!(p-i)!}a_0^ib_0^{p-i}$.
{\cora Note that, on a formal level, dividing by $p$ does not make sense in a ring of characteristic $p$.}

Witt vectors come with the following three maps 
\begin{alignat*}{4}
&(\textrm{Frobenius}) &&F \colon W_nA \to W_nA, \qquad &&F(a_0, a_1, \ldots, a_{n-1}) &&:= (a_0^p,a_1^p, \ldots, a_{n-1}^p),\\ 
&(\textrm{Verschiebung}) \quad  && V \colon W_nA \to W_{n+1}A, \qquad  &&V(a_0, a_1, \ldots, a_{n-1}) &&:= (0, a_0, a_1,\ldots,  a_{n-1}),\\
&(\textrm{Restriction}) && R \colon W_{n+1}A \to W_{n}A, \qquad  &&R(a_0, a_1, \ldots, a_{n}) &&:= (a_0, a_1, a_2,\ldots,  a_{n-1}),
\end{alignat*}
where $F$ and $R$ are ring homomorphisms and $V$ is an additive homomorphism. 
By abuse of notation, we will denote the composition
{$FR \colon W_{n} A \xrightarrow{R} W_{n-1}A \xrightarrow{F} W_{n-1}A$ 
\[
(a_0, a_1, \ldots, a_{n-1}) \xmapsto{FR} (a_0^p,a_1^p, \ldots, a_{n-2}^p)
\]
also by $F$. {\cora Note that the order of composition does not matter as $FR$ is equal to $RF \colon W_{n} A \xrightarrow{F} W_{n}A \xrightarrow{R} W_{n-1}A$.}}


{Using the above convention, } Frobenius and Verschiebung are connected by the key identities $VF = FV = p$ {on $W_{n}A$}; more precisely
\[
(0,a^p_0, \ldots, a^p_{n-2}) = VFx = FVx = p \cdot x \in W_{n}A,
\]
for every $x = (a_0,a_1,\ldots,a_{n-1}) \in W_{n}A$. Here $p=(0,1,0,\ldots,0) \in W_{n}A$. Since it is often a point of confusion, we emphasise that $F(p) = p$. Last, {we point out that up to}  identifying $W_nA$ with $F_*W_nA$ as sets, the Verschiebung is {not just an additive homomorphism, but} a $W_{n+1}A$-module homomorphism $V \colon F_*W_nA \to W_{n+1}A$.

Witt vectors lie in the following short exact sequence of $W_nA$-module homomorphisms 
\[
0 \to F_*W_{n-1}A \xrightarrow{V} W_{n}A \xrightarrow{{\cora R^{n-1}}} A \to 0.
\]
This short exact sequence may be visualised through the identity of $A$-modules:
\[
\mathrm{Gr}_V(W_nA) = A \oplus F_*A \oplus \cdots \oplus F^{n-1}_*A,
\]
where $\mathrm{Gr}_V(W_nA)$ denotes the graded sum {with respect to} the filtration on $W_nA$ given by the inclusions $W_n A \xhookleftarrow{V} F_*W_{n-1}A \xhookleftarrow{V} F^2_*W_{n-2}A \xhookleftarrow{V} \cdots$.

Last, we define the Teichm\"uller lift 
\[
A \to W_nA, \qquad a \mapsto [a] := (a,0,\ldots, 0). 
\]
This is a multiplicative, but \emph{not} an additive, map. In this language, the definition of addition in $W_2(R)$ implies the following important identity:
\begin{equation} \label{eq:addition-in-W2-using-V}
[a] + [b] = [a+b] - V\Big(\frac{(a+b)^p - a^p - b^p}{p}\Big) \in {\cora W_2A}
\end{equation}
for $n=2$ and $a,b \in A$.

The construction of Witt vectors explained above localises. 
In particular, given an $\F_p$-scheme $X$,  
we define the Witt vector sheaf $W_n\cO_X$ on the topological space $|X|$ by $\Gamma(U, W_n\MO_X) := W_n(\Gamma(U, \MO_X))$. 
It is easy to see that $W_nX :=(|X|, W_n\MO_X)$ is actually a scheme and 
the morphism $X \to W_nX$ induced by the ring homomorphism $R^{n-1} \colon W_n\MO_X \to \MO_X$ is a surjective closed immersion {(a `thickening')}. All the properties of Witt vectors stated above extend to $W_n\cO_X$. 
In particular, the Frobenius operator induces a morphism of schemes $F \colon W_nX \to W_nX$, and we have the following short exact sequence of $W_n\cO_X$-module homomorphisms 
\begin{equation} \label{eq:key-sequence-for-witt-vectors}
0 \to F_*W_{n-1}\cO_X \xrightarrow{V} W_n \cO_X \xrightarrow{R^{n-1}} \cO_X \to 0.    
\end{equation}
{Throughout the article, we will implicitly use the fact that $F \colon W_nX \to W_nX$ is a finite morphism when $X$ is $F$-finite.}

\begin{remark}
We warn the reader that 
if $X$ is a variety over a perfect field $k$ of  characteristic $p>0$ and $\dim X >0$, 
then $W_nX$ is \emph{not} a lift of $X$ to $W_n(k)$ {as explained by either of the following observations.}
\begin{enumerate}
    \item $W_n\cO_X$ is not necessarily flat over $W_n(k)$.  
    \item The induced ring homomorphism $W_n\MO_X/pW_n\MO_X \to \MO_X$ is not an isomorphism, 
    that is, the fibre of the induced morphism $W_nX \to \Spec W_n(k)$ over the closed point {of $\Spec W_n(k)$} does not coincide with $X$. 
\end{enumerate}
\end{remark}

\subsection{Quasi-$F$-splittings} \label{ss:quasi-F-splittings}

In this subsection, we state basic properties of quasi-$F$-splittings following \cite{yobuko19}. In what follows, $X$ is a Noetherian $\F_p$-scheme. 

First, we restate the definition of  a quasi-$F$-split variety from the introduction (see (\ref{diagram:intro-definition})) to make it easier to work with. Consider the following pushout diagram of $W_n\cO_X$-modules 
{\cora (in the category of $W_n\MO_X$-modules not the one of rings):}  
\begin{center}
\begin{tikzcd}
W_n\cO_X \arrow{r}{F} \arrow{d}{R^{n-1}} & F_* W_n\cO_X \arrow{d}  \\
\cO_X \arrow{r}{\Phi_{X, n}} & Q_{X,n}. \arrow[lu, phantom, "\usebox\pushoutdr" , very near start, yshift=0em, xshift=0.6em, color=black] 
\end{tikzcd}
\end{center}
By the universal property of a pushout square, $X$ is quasi-$F$-split if and only if there exists a $W_n\MO_X$-module homomorphism $\alpha \colon Q_{X,n} \to \cO_X$ such that the composition $\cO_X \xrightarrow{\Phi_{X, n}}  Q_{X,n} \xrightarrow{\alpha} \cO_X$ is the identity. Thus, we get the following.

\begin{proposition}[{\cite[Remark 4.3]{yobuko19}}] \label{prop:intro-definition-via-splitting}
A Noetherian $\F_p$-scheme $X$ is quasi-$F$-split if and only if 
$\Phi_{X, n} \colon \cO_X \to Q_{X,n}$ splits as a $W_n\MO_X$-module homomorphism.
\end{proposition}

\noindent Note that $Q_{X,1} = F_*\cO_X$, and so being $F$-split is equivalent to being $1$-quasi-$F$-split. Further, the restriction map $F_*W_m \cO_X \to F_*W_n \cO_X$ for $m \geq n$ induces a map $Q_{X,m} \to Q_{X,n}$. Thus if $X$ is $n$-quasi-$F$-split, then it is $m$-quasi-$F$-split for every $m \geq n$.

We can calculate the $W_n\cO_X$-module $Q_{X,n}$ explicitly. 
\begin{proposition} \label{prop:intro-properties-of-C}
Let $X$ be a Noetherian $\F_p$-scheme. Then the following hold. 
\begin{enumerate}
    \item {The projection $F_*W_n\cO_X \to Q_{X,n}$ induces an isomorphism of $W_n\MO_X$-modules:}
    \[
    Q_{X,n} \cong \frac{F_*W_n\cO_X}{p F_*W_n\MO_X}.
    \]
    \item $V(W_n\cO_X)$ acts trivially on $Q_{X,n}$, and so $Q_{X,n}$ is naturally a 
quasi-coherent $\cO_X$-module. 
{Hence}, if $X$ is $F$-finite, then $Q_{X,n}$ is a coherent $\MO_X$-module.
\end{enumerate}
\end{proposition}

\begin{proof}
Assertion (1) follows from 
\[
Q_{X,n} = \frac{F_*W_n\MO_X}{F(\Ker {R^{n-1}})} = 
\frac{F_*W_n\MO_X}{  F( V(F_*W_{n-1}\MO_X))} = 
\frac{F_*W_n\MO_X}{  p F_*W_{n-1}\MO_X},
\]
where the first equality holds by definition of pushouts 
and the third one {by} $FV=p$. 

Let us show (2). 
Pick an open subset $U \subseteq X$ and $s \in W_n\cO_X(U)$. Then 
\[
Vs \cdot Q_{X,n}(U) = F_*(FVs \cdot W_n \cO_X(U)/p) =F_*(ps \cdot W_n\cO_X(U)/p) = 0,
\]
and so $Q_{X,n}$ is an $\cO_X$-module. Since all the sheaves in question are quasi-coherent {$W_n\cO_X$-modules}, we get that $Q_{X,n}$ is a quasi-coherent {$\cO_X$-module}. Moreover, if $X$ is $F$-finite, then $F_*W_n\cO_X$ is a coherent {$W_n\cO_{X}$-module}, and so $Q_{X,n}$ is a coherent {$\cO_X$-module}.
\end{proof}
{From now on, assume that $X$ is reduced.  This is needed to ensure that Frobenius $F \colon \cO_X \to F_*\cO_X$ is injective.} We state a few key short exact sequences that will be used throughout this article {and refer to Lemma \ref{lem:key-sequences-in-the-log-case} for a proof in a more general setting}. First, by definition of a pushout again, the following sequence of $\cO_X$-modules
\begin{equation} \label{eq:intro-C-quotient-sequence}
0 \to \cO_X \to Q_{X,n} \to F_*W_n \cO_X / W_n\cO_X \to 0
\end{equation}
is exact. Here, we implicitly treat  $W_n\cO_X$ as a submodule of $F_*W_n\cO_X$ embedded by 
Frobenius $F \colon W_n\cO_X \hookrightarrow F_*W_n\cO_X$.  
The term on the right lies naturally in the following sequence
\begin{equation} \label{eq:intro-Bn-sequence}
0 \to F_*(F_*W_{n-1}\cO_X/W_{n-1}\cO_X) \to F_*W_n \cO_X /W_n \cO_X \to F_* \cO_X / \cO_X \to 0
\end{equation}
induced from (\ref{eq:key-sequence-for-witt-vectors}). Finally, the restriction  $F_*W_n\cO_X \to F_*\cO_X$ induces  {an $\cO_X$-module homomorphism  $Q_{X,n} \to F_*\cO_X$ sitting inside the following sequence}
\begin{equation} \label{eq:intro-C-restriction-sequence}
0 \to F_*(Q_{X,n-1}/\cO_X) \to Q_{X,n} \to F_* \cO_X \to 0,
\end{equation}
where we note that $Q_{X,n-1}/\cO_X \cong F_*W_{n-1}\cO_X/W_{n-1}\cO_X$ by \eqref{eq:intro-C-quotient-sequence}.

\begin{proposition} \label{prop:C-locally-free} Let $X$ be a smooth scheme over a perfect field of characteristic $p>0$. Then $Q_{X,n}$ is a locally free $\MO_X$-module.  
\end{proposition}
\begin{proof}
The scheme $X$ is smooth, and so both $F^e_*\cO_X$ and $F^e_*\cO_X/F^{e-1}_*\cO_X$ are locally free for every $e \geq 1$. Since an extension of a locally free sheaf by a locally free sheaf is locally free (as local $\cExt$ vanish), a repeated use of (\ref{eq:intro-Bn-sequence}) shows that $F_*W_n\cO_X/W_n\cO_X$ is locally free. Therefore, $Q_{X,n}$ is locally free by (\ref{eq:intro-C-quotient-sequence}).
\end{proof}

The $\MO_X$-module 
homomorphism ${\Phi_{X, n}} \colon \cO_X \to Q_{X,n}$ splits if and only if the evaluation map
\begin{equation}
\label{eq:intro-dual-defin-quasi-F-split}
\Psi_{X, n} 
\colon Q^\vee_{X,n} = \cHom(Q_{X,n}, \cO_X) \to \cHom(\cO_X, \cO_X) = \cO_X
\end{equation}
is surjective on global sections. This yields cohomological criterions for checking whether a scheme is  quasi-$F$-split. 
\begin{lemma} \label{lem:intro-cohomological-quasi-F-spliteness}
Let $X$ be an integral Noetherian $\F_p$-scheme of dimension $d$. Then $X$ is quasi-$F$-split if and only if
\begin{enumerate}
    \item $H^d(X, \omega_X) \to H^d(X, Q_{X,n} \otimes \omega_X)$ is injective, when $X$ is smooth and projective over a perfect field $k$ of characteristic $p>0$,
    \item $H^d_{\m}(X, \omega_X) \to H^d_{\m}(X, Q_{X,n} \otimes \omega_X)$ is injective, when $X = \Spec R$ for an excellent normal local domain $(R,\m)$  admitting a dualizing complex.
\end{enumerate}
\end{lemma}
\noindent The above maps are given by tensoring $\cO_X \!\to Q_{X,n}$ with $\omega_X$ and applying cohomology. 
\begin{proof}
 (1) follows by applying Serre duality to (\ref{eq:intro-dual-defin-quasi-F-split}). As for (2), the statement follows by Matlis duality (cf.\ Proposition \ref{prop:Matlis-duality})  as it interchanges injectivity with surjectivity. {For details we refer to the proof of a more general Lemma \ref{lem:cohomological-criterion-for-log-quasi-F-splitting}.}
\end{proof}

\begin{remark} \label{remark:elliptic-curves-quasi-F-split}
{Let $E$ be a smooth projective curve over an $F$-finite field of characteristic $p>0$ such that $K_E$ is trivial. Then $E$ is $2$-quasi-$F$-split. Indeed, this can be checked after a base change of the field $k$ (see Corollary \ref{c-descent2} or \cite[Corollary 2.18]{KTY22}), 
so that we may assume that $E$ is an elliptic curve. In this case, $E$ is $2$-quasi-$F$-split by \cite[Corollary IV.7.5]{silverman-elliptic-curves-book} in view of the definition of the quasi-$F$-split height using formal groups (see  \cite[Theorem 4.5]{yobuko19}).}
\end{remark}

\subsection{Higher Cartier operator} \label{ss:iterated-Cartier-operator}
A standard reference for the results of this subsection is \cite{illusie_de_rham_witt}. In what follows, we assume that $X$ is a smooth variety over a perfect field $k$ of characteristic $p>0$. The starting point of our discussion is the Cartier isomorphism $C^{-1} \colon \Omega_X^i \xrightarrow{\cong}  \cH^{i}(F_*\Omega_X^{\mydot})$. Note that $\Omega_X^{\mydot}$ is \emph{not} an $\cO_X$-linear complex, but $F_*\Omega_X^{\mydot}$ is, as $d(f^p\omega)= f^pd\omega$ for local sections $f$ and $\omega$ of $\cO_X$ and $\Omega_X^i$, respectively. 

We define $B_1\Omega_X^i$ and $Z_1\Omega_X^i$ to be the boundaries and the cycles of $F_*\Omega_X^{\mydot}$ at $i$: 
\begin{align} 
\label{eq:definition-of-B1} B_1 \Omega_X^i &:= {\rm Im}(F_*\Omega_X^{i-1} \xrightarrow{F_*d}  F_*\Omega_X^{i}),\\ 
Z_1 \Omega_X^i &:= \Ker(F_*\Omega_X^{i} \xrightarrow{F_*d}  F_*\Omega_X^{i+1}) \nonumber.
\end{align}
Then the Cartier isomorphism gives the following short exact sequence
\begin{equation} \label{eq:Cartier-for-Z1}
0 \to B_1\Omega^i_X \to Z_1\Omega^i_X \xrightarrow{C}  \Omega^i_X \to 0. 
\end{equation}

In order to capture more refined information on the cohomology of $F^n_*\Omega_X^{\mydot}$ one can iterate the above construction. Specifically, we will explain how to detect $\Omega_X^i$ as a subquotient of $\cH^i(F^n_*\Omega_X^{\mydot})$.
\begin{remark} \label{remark:iterated-partial-map}
One can think of $C$ as a \emph{partial map} $F_*\Omega^i_X \rightharpoonup \Omega^i_X$ with the domain of definition equal to $Z_1\Omega^i_X$. Then the higher Cartier operator
\[
C^n \colon F^n_*\Omega_X^i \stackrel{F^{n-1}_*C}{\rightharpoonup} F^{n-1}_*\Omega_X^i \stackrel{F^{n-2}_*C}{\rightharpoonup} \ldots \stackrel{C}{\rightharpoonup} \Omega_X^i  
\]
will be the composition of the partial maps $F^{n-1}_*C, \ldots, C$. The sheaf $Z_n\Omega^i_X$ will be the domain of definition of $C^n$ (equivalently, $(C^n)^{-1}(\Omega_X^i)$), and the sheaf $B_n\Omega^i_X$ will be the kernel of $C^n$ (equivalently, $(C^n)^{-1}(0)$). In what follows, we formalise this explanation without using the notion of a partial map.
\end{remark}

Set $C_n := F^{n}_{*}C \colon F^n_* Z_1\Omega_X^i \to F^n_*\Omega_X^i$.
\begin{definition}
We define { a coherent $\cO_X$-module $Z_n\Omega^i_X$} inductively by the formula
\[
Z_n\Omega_X^i := (C_{n-1})^{-1}(Z_{n-1}\Omega_X^i) \subseteq F^{n}_* \Omega^i_X,
\]
where $Z_0\Omega^i_X=\Omega^i_X$. 
\end{definition}
This construction may be visualised through the following pullback diagrams:
\[
\begin{tikzcd}
F_*Z_1\Omega_X^i \arrow{r}{C_1} & F_* \Omega_X^i\\
Z_2\Omega_X^i \arrow[ru, phantom, "\usebox\pullbackdl" , very near start, yshift=-0.3em, xshift=-0.6em, color=black] \arrow[hook]{u} \arrow{r} & Z_1\Omega^i_X \arrow[hook]{u} 
\end{tikzcd}
\hspace{1em}
\begin{tikzcd}
F^2_*Z_1\Omega_X^i \arrow{r}{C_2} & F^2_* \Omega_X^i\\
Z_3\Omega_X^i \arrow[ru, phantom, "\usebox\pullbackdl" , very near start, yshift=-0.3em, xshift=-0.6em, color=black] \arrow[hook]{u} \arrow{r} & Z_2\Omega^i_X \arrow[hook]{u} 
\end{tikzcd}
\hspace{1em}
\cdots
\hspace{1em}
\begin{tikzcd}
F^{n-1}_*Z_1\Omega_X^i \arrow{r}{C_{n-1}} & F^{n-1}_* \Omega_X^i\\
Z_n\Omega_X^i \arrow[ru, phantom, "\usebox\pullbackdl" , very near start, yshift=-0.3em, xshift=-0.6em, color=black] \arrow[hook]{u} \arrow{r} & Z_{n-1}\Omega^i_X. \arrow[hook]{u} 
\end{tikzcd}
\]
Since $C_{n-1}$ is surjective, so is $C_{n-1}|_{Z_n\Omega_X^i} \colon Z_{n}\Omega_X^i \to Z_{n-1}\Omega_X^i$.

\begin{definition}
We define a {coherent $\cO_X$-module $B_n\Omega_X^i$}  inductively by the formula
\[
B_n\Omega_X^i := (C_{n-1})^{-1}(B_{n-1}\Omega_X^i) \subseteq F^{n}_* \Omega^i_X,
\]
where $B_0\Omega^i_X=0$. 
\end{definition}

\begin{definition} We define the \emph{higher  Cartier operator} as the composition 
\[
C^{n} :=(C_0|_{Z_1\Omega_X^i}) \circ \cdots \circ (C_{n-1}|_{Z_n\Omega_X^i})  \colon Z_n \Omega^i_X \to \Omega^i_X
\]
of the above surjective maps $Z_n\Omega^i_X \to Z_{n-1}\Omega^i_X \to \ldots \to \Omega^i_X$.
\end{definition}
 This name is justified by the following short exact sequence
\begin{equation} \label{eq:higher-Cartier-operator}
0 \to B_n \Omega^i_X \to Z_n \Omega^i_X \xrightarrow{C^{n}} \Omega^i_X 
\to 0,
\end{equation}
which exists, and is exact, as $B_n\Omega^i_X \subseteq Z_n \Omega^i_X$, the map $C^{n}$ is surjective, and 
\[
B_n\Omega^i_X =  (C_{n-1})^{-1}(B_{n-1}\Omega^i_X) = \cdots = (C_{0}\circ \cdots \circ C_{n-1})^{-1}(B_{0}\Omega^i_X) = (C^{n})^{-1}(0).
\]
Here, $C_{0}\circ \cdots \circ C_{n-1}$ does not make sense as a map, but $(C_{0}\circ \cdots \circ C_{n-1})^{-1}(M)$ for any subsheaf $M \subseteq \Omega^i_X$ is well-defined.   
\begin{lemma} \label{lemma:key-ses-for-Zn} With notation as above, the following sequence is exact.
\[
0 \to Z_n\Omega^i_X \to F_*Z_{n-1}\Omega^i_X \xrightarrow{F_*d \circ F_*C^{n-1}} B_1\Omega^{i+1}_X \to 0.
\]
\end{lemma}
\noindent This lemma will be crucial for us later on.
\begin{proof}Since taking a pull-back of a map {by a surjection} preserves its cokernel, the above sequence is constructed through the following diagram
\begin{center}
\begin{tikzcd}
0 \arrow{r} & Z_1\Omega^i_X \arrow{r} & F_*\Omega^i_X \arrow{r}{F_*d} & B_1\Omega^{i+1}_X \arrow{r} & 0\\
0 \arrow{r} & Z_n\Omega^i_X \arrow[ru, phantom, "\usebox\pullbackdl" , very near start, yshift=-0.3em, xshift=-0.6em, color=black] \arrow{u} \arrow{r} & F_*Z_{n-1}\Omega^i_X \arrow{u}[swap]{F_*C^{n-1}} \arrow{r} & B_1\Omega^{i+1}_X \arrow{u}[swap]{=} \arrow{r} & 0,
\end{tikzcd}    
\end{center}
in which the upper row is exact by definition of cycles and boundaries and the left square is Cartesian because by definition of $Z_n\Omega^i_X$ we have the equality
\[
Z_n\Omega^i_X = (C_{1}\circ \cdots \circ C_{n-1})^{-1}(Z_1\Omega^i_X) = (F_*C^{n-1})^{-1} (Z_1\Omega^i_X). \qedhere
\]
\end{proof}
Last, we recall the famous observation of Serre (cf.\ Lemma \ref{lem:log-Serre's map}) that the map
\begin{equation}\label{eq:Serre-map}
F_*W_n\cO_X/W_n\cO_X \xrightarrow{\cong} B_n\Omega^1_X.
\end{equation}
given by $(a_0,a_1,\ldots, a_{n-1}) \mapsto \frac{1}{p^{n-1}}d(a_0^{p^{n-1}} + pa_1^{p^{n-2}} + \cdots + p^{n-1}a_{n-1})$ is an isomorphism. Since $F_*W_n\cO_X/W_n\cO_X$ is equal to the cokernel of $\cO_X \to {\cora Q_{X,n}}$, this observation connects the study of quasi-$F$-splittings with that of the higher Cartier operator.

\subsection{Log $F$-splittings}\label{ss log F-split}

We say that an $\F_p$-scheme $X$ is \emph{$F$-split} if $F \colon \cO_X \to F_*\cO_X$ splits as an $\cO_X$-module homomorphism.
There are three natural ways to define $F$-splittings for a pair $(X,\Delta)$.  
\begin{definition}\label{d-Fsplit}
Let $X$ be an integral normal Noetherian scheme and 
let $\Delta$ be an effective $\Q$-divisor. 
We say that $(X,\Delta)$ is 
\begin{align*}
 \textrm{\emph{F-split}} &\iff \cO_X \xrightarrow{F^e} F^e_*\cO_X(\rdown{(p^e{-}1)\Delta})  \textrm{ splits for all } e>0.  \\
  \textrm{\emph{sharply $F$-split}} &\iff  \cO_X \xrightarrow{F^e} F^e_*\cO_X(\rup{(p^e{-}1)\Delta}) \textrm{ splits for some } e>0. \\
 \textrm{\emph{keenly $F$-split}} &\iff  \cO_X \xrightarrow{F^e} F^e_*\cO_X(\rdown{p^e\Delta}) \textrm{ splits for all } e>0.
\end{align*}
\end{definition}
We have the following implications in the case when $\rdown{\Delta} = 0$:
\begin{equation}\label{e1-Fsplit}
\textrm{sharply }F\textrm{-split} \implies \textrm{keenly }F\textrm{-split} \implies F\textrm{-split}.
\end{equation}

\begin{lemma}
(\ref{e1-Fsplit}) holds. 
\end{lemma}

\begin{proof}
{The latter implication is clear. 
Let us show the former one. 
Assume that $(X, \Delta)$ is sharply $F$-split, 
 that is, $F^e \colon \MO_X \to F_*^e \MO_X(\rup{(p^e{-}1)\Delta})$ splits for some $e>0$. 
Recall that 
also $F^{en} \colon \MO_X \to F_*^{en} \MO_X(\rup{(p^{en}{-}1)\Delta})$ 
splits for any $n \in \Z_{>0}$ \cite[Proposition 3.8(b)]{schwedesmith10}. 
Fix $n \in \Z_{>0}$ and $\delta \in \Q$ {satisfying} $0<\delta <1$. 
Since the splitting of 
$F^d \colon \MO_X \to F_*^d \MO_X(\rdown{p^d \Delta})$ implies 
the splitting of 
$F^{d-1} \colon \MO_X \to F_*^{d-1} \MO_X(\rdown{p^{d-1} \Delta})$, 
it suffices to show that 
\begin{equation}\label{e2-Fsplit}
\rdown{p^{en}\delta} \leq \rup{(p^{en}-1)\delta}. 
\end{equation}
If $(p^{en}-1)\delta \in \Z$, then (\ref{e2-Fsplit}) holds by 
\[
\rdown{p^{en}\delta} = 
\rdown{(p^{en}-1)\delta + \delta} = (p^{en}-1)\delta = \rup{(p^{en}-1)\delta}. 
\]
If $(p^{en}-1)\delta \not\in \Z$, then (\ref{e2-Fsplit}) follows from  
\[
\rdown{p^{en}\delta} \leq \rdown{(p^{en}-1)\delta} +1 
=\rup{(p^{en}-1)\delta}. \qedhere
\]}
\end{proof}

The notion of a keen $F$-splitting is less common; the name has been suggested to us by Schwede and Smith who are preparing a book on $F$-splitting.
The notation is introduced by the paper \cite{Hara-Watanabe} and called strong $F$-purity in the paper. 
One big advantage of keen $F$-splitting is that it has good properties without the necessity of  perturbing the coefficients to make $(p^e-1)\Delta$ integral. From our point of view, keen $F$-splitting is the notion that generalises well to log quasi-$F$-splittings. Since we do not work with powers of Frobenius, we also need to make the following definition.

\begin{definition}
Let $X$ be an integral normal Noetherian $\F_p$-scheme and 
let $\Delta$ be an effective $\Q$-divisor. We say that $(X,\Delta)$ is \emph{naively keenly $F$-split} if $F \colon \cO_X \to F_*\cO_X(\rdown{p\Delta})$ splits as an $\MO_X$-module homomorphism. We say that $(X,\Delta)$ is \emph{naively keenly $F$-pure} {if it is naively keenly $F$-split affine-locally (explicitly, for every $x \in X$, there exists an affine open neighbourhood $U \subseteq X$ of $x$ such that $\cO_X(U)$ is naively keenly $F$-split, cf.\ \cite[\href{https://stacks.math.columbia.edu/tag/01OQ}{Tag 01OQ}]{stacks-project}).}
\end{definition}

Unfortunately, keen $F$-splittings do not behave well when some of the coefficients of $\Delta$ are equal to $1$. Thus, in order to set up adjunction properly, we need to make the following definition (motivated by $F$-regularity and pure $F$-regularity).

\begin{definition}
Let $X$ be an integral normal Noetherian scheme{, let $S$ be a prime divisor, and let $B$ be an effective $\bQ$-divisor such that $S \not \subseteq \Supp B$}. We say that $(X,{S+B})$ is \emph{purely naively keenly $F$-split {(along $S$)}} if $F \colon \cO_X \to F_*\cO_X((p-1)S + \lfloor pB \rfloor)$ splits. We say that $(X,\Delta)$ is \emph{purely naively keenly $F$-pure {(along $S$)}} if it is purely naively keenly $F$-split (along $S$) affine-locally.
 

\end{definition}

\begin{lemma} \label{lem:p-compatible-splitting}
Let $X$ be an integral normal excellent $\F_p$-scheme and 
let $\Delta$ be a $\Q$-divisor. Then for every $i\geq 0$:
\begin{enumerate}
\setlength\itemsep{0.3em}
    \item $\rdown{p^{i+1}\Delta} - p\rdown{p^{i}\Delta} = \rdown{p\{p^{i}\Delta\}}$, and 
\item $F \colon \cO_X(\rdown{p^i\Delta}) \to F_*\cO_X(\rdown{p^{i+1}\Delta})$  splits as an $\cO_X$-module homomorphism provided that $(X,\{p^{i}\Delta\})$ is naively keenly $F$-split.
\end{enumerate}
\end{lemma}
\begin{proof}
{\cora We first prove (1).
Since the coefficients of $p\rdown{p^{i}\Delta}$ are integers, it follows that
\[\rdown{p^{i+1}\Delta}-p\rdown{p^{i}\Delta}=\rdown{p^{i+1}\Delta-p\rdown{p^{i}\Delta}}.\]
Then (1) 
holds as
\[
\rdown{p^{i+1}\Delta} - p\rdown{p^{i}\Delta} = \rdown{p^{i+1}\Delta-p\rdown{p^{i}\Delta}}= \rdown{p(p^{i}\Delta-\rdown{p^{i}\Delta})} = \rdown{p\{p^{i}\Delta\}}. 
\]}
(2) follows immediately from (1) {(cf.\ Remark \ref{remark:hom-reflexive})}. 
\end{proof}

\subsection{Witt divisorial sheaves}


For the convenience of the reader, we briefly review the construction and properties of Witt divisorial sheaves as constructed by the third author, and refer to \cite[Section 3]{tanaka22} for details. Throughout this subsection, $X$ is an integral normal $F$-finite Noetherian $\F_p$-scheme. 
Recall that $K(X)$ denotes the function field of $X$. 

For a $\bQ$-divisor $D$, we define a subsheaf $W_n\cO_X(D)$ of the constant sheaf $W_n(K(X))$ on $X$ by the formula
\[
W_n\cO_X(D)(U) := \{(\varphi_0, \ldots, \varphi_{n-1}) \in W_n(K(X)) \mid \varphi_i \in {\cora\Gamma(U, F^i_*\cO_X(p^iD))}\text{ for every }i\},
\]
for every open subset $U \subseteq X$ (we remind the reader that $\cO_X(p^iD) = \cO_X(\rdown{p^iD})$). One can verify that this is in fact a sheaf (\cite[Lemma 3.5(1)]{tanaka22}), that 
it is a coherent $W_n\MO_X$-module  (\cite[Proposition 3.8]{tanaka22}), {and that it comes equipped with the following three $W_n\cO_X$-modules homomorphisms
\begin{alignat*}{2}
&(\textrm{Frobenius}) &&F \colon 
W_n\cO_X(D) \longrightarrow F_*W_n\cO_X(pD),  \\ 
&(\textrm{Verschiebung}) \qquad &&V \colon F_*W_{n}\cO_X(pD) \longrightarrow W_{n+1}\cO_X(D), \\
&(\textrm{Restriction}) &&R \colon W_{n+1}\cO_X(D) \longrightarrow W_{n}\cO_X(D).
\end{alignat*}
These are induced by  $F \colon W_n(K(X)) \to  F_*W_n(K(X))$,  $V \colon F_*W_{n}(K(X))\to W_{n+1}(K(X))$, and  $R \colon W_{n+1}(K(X)) \to W_{n}(K(X))$, respectively. Finally,} we warn the reader that $W_n\cO_X(D)$ need not be equal to $W_n\cO_X(\rdown{D})$. 
{\cora For example, $W_2\MO_X((1/p)D) \neq 
W_2\MO_X(\rdown{(1/p)D}) = W_2\MO_X$ for a prime divisor $D$ on $X$.} 

Now (\ref{eq:key-sequence-for-witt-vectors}) induces the following exact sequence \cite[Proposition 3.7]{tanaka22}:
\begin{equation} \label{eq:key-sequence-for-log-witt-vectors}
0 \to F_*W_{n-1}\cO_X(pD) \xrightarrow{V} W_n \cO_X(D) \xrightarrow{R^{n-1}} \cO_X(D) \to 0,
\end{equation}
which can be visualised through the identity
\[
\mathrm{Gr}_V(W_n\cO_X(D)) = \cO_X(D) \oplus F_*\cO_X(pD) \oplus \cdots \oplus F^{n-1}_*\cO_X(p^{n-1}D).
\]
Here $\mathrm{Gr}_V(W_n\cO_X(D))$ denotes the graded sum obtained 
by the filtration on $W_n\MO_X(D)$ given by the inclusions  $W_n\cO_X(D) \xhookleftarrow{V} F_*W_{n-1}\cO_X(pD) \xhookleftarrow{V} F^2_*W_{n-2}\cO_X(p^2D) \xhookleftarrow{V} \cdots$.


Since $X$ is normal, the above short exact sequence 
(\ref{eq:key-sequence-for-log-witt-vectors}) 
shows that $W_n\cO_X(D)$ is $S_2$. 
In particular, if $D$ {and $E$} are Weil divisors on $X$ and
$j \colon X_{\reg} \hookrightarrow X$ is the open immersion 
from the regular locus $X_{\reg}$ of $X$, then 
\[
W_n\cO_X(D + E) \cong j_*j^*(W_n\cO_X(D) \otimes_{W_n\MO_X} W_n\cO_X(E)).
\]
Here we use the fact that when $E$ is Cartier, $W_n\cO_X(E)$ is an invertible $W_n\cO_X$-module (a line bundle on $W_nX$), which is in fact equal to the Teichm\"{u}ller lift\footnote{{Specifically, if $E|_{U_i}= \mathrm{div}(f_i)$ for an affine open cover $\{U_i\}$ of $X$ and sections $f_i \in \cO_X(U_i)$, then there are natural identifications $\cO_X(E)|_{U_i} = f_i^{-1}\cO_{U_i}$ and $W_n\cO_X(E)|_{U_i} =[f_i]^{-1}W_n\cO_{U_i}.$}} of $\cO_X(E)$ (see \cite[Proposition 3.12]{tanaka22}). 

\medskip


We now introduce a $W_n\MO_X$-submodule $W_n\cI_S(D) \subseteq W_n\MO_X(D)$, which is a generalisation of $W_n\MO_X(D)$ constructed above. 
For an integral normal $F$-finite Noetherian $\F_p$-scheme $X$, 
a $\Q$-divisor $D$, and an effecitve $\Q$-divisor $S$, 
 we define 
 a subpresheaf $W_n\cI_S({D}) \subseteq W_n\cO_X({D})$ by the formula
\[
W_n\cI_S({D})(U) := \{(\varphi_0,\ldots,\varphi_{n-1}) \mid 
\varphi_i \in F^i_*\cO_X(-S+p^i{D}) \text{ for every } i\} \subseteq W_n\cO_X({D})(U)
\]
for every open subset $U \subseteq X$. 
By the same argument as in \cite[Subsection 3.1]{tanaka22}, 
the following hold for $n,m \in \Z_{>0}$. 
\begin{enumerate}
\item 
$W_n\cI_S(D)$ is an $S_2$ coherent $W_n\MO_X$-submodule of $W_n\MO_X(D)$ 
(cf.\ {\cora a proof below and} \cite[Proposition 3.8]{tanaka22}). 
\item There exist $W_n\cO_X$-module homomorphisms
{
\begin{alignat*}{2}
&(\textrm{Frobenius}) &&F \colon 
W_n\cI_S(D) \longrightarrow F_*W_n\cI_S(pD),  \\ 
&(\textrm{Verschiebung}) \qquad &&V \colon F_*W_{n}\cI_S(pD) \longrightarrow W_{n+1}\cI_{S}(D), \\
&(\textrm{Restriction}) &&R \colon W_{n+1}\cI_S(D) \longrightarrow W_{n}\cI_{S}(D),
\end{alignat*}
induced by  $F \colon W_n(K(X)) \to  F_*W_n(K(X))$,  $V \colon F_*W_{n}(K(X))\to W_{n+1}(K(X))$, and  $R \colon W_{n+1}(K(X)) \to W_{n}(K(X))$, respectively.}
\item 
We have the following exact sequence (cf.\ \cite[Proposition 3.7]{tanaka22}): 
\begin{equation} \label{eq:key-sequence-for-WnI}
0\to F_*^nW_m\cI_S(p^nD) \xrightarrow{V^n}  
W_{n+m}\cI_S(D) \xrightarrow{R^m} W_n\cI_S(D) \to 0.
\end{equation}
\item 
{\cora It holds that $W_n\cO_X(-S) \subsetneq W_n\cI_S$ 
for the case when $n \geq 2$ and $S$ is a nonzero effective Cartier divisor.} 
\end{enumerate}

{\cora 
\begin{proof}[Sketch of {\rm (1)}]
For the reader's convenience, 
we include a proof of the fact that $W_n\cI_S(D)$ is closed under addition. 
The following argument is extracted from \cite[the proof of Lemma 3.5(1)]{tanaka22}. 

We may assume that $X$ is affine. 
Fix $\varphi, \psi \in \Gamma(X, W_n\cI_S(D))$. 
Let us prove $\varphi+\psi \in \Gamma(X, W_n\cI_S(D))$. 
We can write 
$$\varphi=(\varphi_0, \varphi_1, \ldots, \varphi_{n-1}), \quad \psi=(\psi_0, \psi_1, \ldots, \psi_{n-1})$$
for some $\varphi_m, \psi_m \in \Gamma(X, \MO_X(p^mD-S))$. 
By \cite[Lemma 2.3]{tanaka22}, 
we have that 
$$\varphi+\psi=(S_0(\varphi_0, \psi_0), S_1(\varphi_0, \psi_0, \varphi_1, \psi_1), 
\ldots, S_{n-1}(\varphi_0, \psi_0, \varphi_1, \psi_1, \ldots, \varphi_{n-1}, \psi_{n-1}))$$
for some polynomials 
$$S_m(x_0, y_0, \ldots, x_m, y_m) \in \Z[x_0, y_0, \ldots, x_m, y_m]$$
satisfying the properties listed in \cite[Lemma 2.3]{tanaka22}. 
We equip the polynomial ring $\Z[x_0, y_0, x_1, y_1, \ldots]$ with the structure of graded $\Z$-algebra defined in the statement of \cite[Lemma 2.4]{tanaka22}, 
i.e., we consider $\Z[x_0, y_0, x_1, y_1, \ldots]$ as a wighted polynomial ring 
such that $\deg x_i=\deg y_i=p^i$. 
Pick a monomial 
$x_0^{a_0}y_0^{b_0} \cdots x_m^{a_m}y_m^{b_m}$ appearing in 
the monomial decomposition of $S_m(x_0, y_0, \ldots, x_m, y_m)$. 
Since $S_m$ is homogeneous of degree $p^m$ \cite[Lemma 2.4]{tanaka22}, 
it holds that 
$$\sum_{i=0}^m p^i(a_i+b_i)=p^m.$$
By 
$${\rm div} (\varphi_i) \geq -p^iD +S
\qquad {\rm and}\quad \quad {\rm div} (\psi_i) \geq -p^iD+S,$$
we have that 
\begin{eqnarray*}
{\rm div} (\varphi_0^{a_0}\psi_0^{b_0} \cdots \varphi_m^{a_m}\psi_m^{b_m})
&=&\sum_{i=0}^m (a_i{\rm div} (\varphi_i)+b_i{\rm div} (\psi_i))\\
&\geq & \sum_{i=0}^m a_i(-p^iD+S) + b_i(-p^iD+S)\\
&=& - \sum_{i=0}^m p^i(a_i+b_i)D + \sum_{i=0}^m (a_i+b_i)S\\
&=& - p^mD + \sum_{i=0}^m (a_i+b_i)S\\
&\geq &- p^mD +S.\\
\end{eqnarray*}
In other words, we obtain 
$\varphi_0^{a_0}\psi_0^{b_0} \cdots \varphi_m^{a_m}\psi_m^{b_m} 
\in \Gamma(X, \MO_X(p^mD-S))$. 
Therefore, it holds that $\varphi+\psi \in \Gamma(X, W_n\cI_S(D))$. 
\end{proof}
}

\subsection{Matlis duality} \label{ss:Matlis-duality}

In this subsection, we recall the foundations of Matlis duality and refer to \cite[\href{https://stacks.math.columbia.edu/tag/08XG}{Tag 08XG}]{stacks-project} for details. 

Let $(R, \m)$ be a local Noetherian ring admitting a normalised dualising complex $\omega^{\mydot}_R$. 
Let $E$ be the injective hull of the residue field $R/\m$ over $R$ (the precise definition thereof will not be {used} in our paper; the reader should 
 be aware {though} that $E$ is an injective Artinian $R$-module). 
The operation $(-)^{\vee} := \Hom_R(-, E)$ is called \emph{Matlis duality}. 
If $(R, \m)$ is complete, then it is an anti-equivalence between the  category of Noetherian $R$-modules and the category of Artinian $R$-modules. Whether $R$ is complete or not, $(-)^\vee$ is exact, and so this operation extends canonically to the corresponding derived categories;
we will denote such an extension by the same symbol $(-)^{\vee}$.

A key property of Matlis duality is that it turns the local cohomology of a complex into its Grothendieck dual
as indicated by the following result. This is a generalisation of Serre duality to the relative setting. In what follows, $(-)^{\wedge}$ denotes the derived $\m$-adic completion (see \cite[\href{https://stacks.math.columbia.edu/tag/0922}{Tag 0922}]{stacks-project} for the definition in the Noetherian case, cf.\ \cite[\href{https://stacks.math.columbia.edu/tag/091V}{Tag 091V}]{stacks-project}), which can be dropped from the statement if $(R,\m)$ is complete.

\begin{proposition}[{cf.\ \cite[\href{https://stacks.math.columbia.edu/tag/0A84}{Tag 0A84}]{stacks-project} and \cite[\href{https://stacks.math.columbia.edu/tag/0AAK}{Tag 0AAK}]{stacks-project}}] \label{prop:Matlis-duality}
Let $(R, \m)$ be a local Noetherian ring and 
let $K \in D^{\mathrm{b}}_{\mathrm{coh}}(R)$. Then
\[
    R\Hom_R(K, \omega^{\mydot}_R)^{\wedge} \cong R\Gamma_{\m}(K)^{\vee}.
\]
In particular, if $K$ is {a finitely generated} $R$-module, then
\[
 \Ext^{-i}_R(K, \omega^{\mydot}_R)^{\wedge} \cong H^{i}_{\m}(K)^{\vee}.
\]
\end{proposition}
\noindent 
In the second equation, $(-)^{\wedge}$ denotes the usual $\m$-adic completion, 
but there is no ambiguity here as the derived {$\m$-adic} completion
agrees with the usual {$\m$-adic} completion for {finitely generated} modules over Noetherian rings by \cite[\href{https://stacks.math.columbia.edu/tag/0EEU}{Tag 0EEU}]{stacks-project}. 

In order to make proofs easier to read, we introduce the following notation used throughout the article.
\begin{notation} \label{notation:global-local-cohomology}
Given $i \in \Z$, 
a {Noetherian} local ring $(R,\m)$, a proper morphism of schemes $g \colon Y \to \Spec R$, and a coherent sheaf $\cF$ on $Y$, we define 
\[
H^i_{\m}(Y,\cF) := H^iR\Gamma_{\m}R\Gamma({Y},\cF), 
\]
which is an Artinian $R$-module. 
Here, $R\Gamma({Y},\cF)$ may be identified with  $Rg_*\cF$.
\end{notation}

Given a short exact sequence of coherent sheaves on $Y$, 
the functors $H^i_{\m}(Y, -)$  induce 
a long exact sequence of cohomologies. {Note that} if $Y=\Spec R$ {and $g = {\rm id}$}, then $H^i_{\m}(Y,\cF) $ agrees with the usual local cohomology, and if $R=k$ is a field and $\m = (0)$, then $H^i_{\m}(Y,\cF)  = H^i(Y,\cF)$ is the coherent cohomology.

\begin{lemma}  \label{lem:Matlis-duality-for-highest-cohomology}
Let $(R, \m)$ be a Noetherian local domain admitting a dualising complex. 
Let $X$ be {a $d$-dimensional} integral scheme  which is proper over $\Spec R$ and let $\cF$ be a coherent sheaf on $X$. Then
\[
{H^d_{\m}(X,\cF)^{\vee}} \cong \Hom_{\cO_X}(\cF, \omega_X)^{\wedge}.
\]
\end{lemma}
\begin{proof}
In view of {Proposition \ref{prop:Matlis-duality}}, the proof is exactly the same as that of \cite[Lemma 2.2]{BMPSTWW20} (wherein it is assumed that $(R,\m)$ is complete). For the convenience of the reader, we briefly summarise the argument. Consider the following isomorphisms
\begin{alignat*}{1}
\big(H^dR\Gamma_{\m}R\Gamma(X,\cF)\big)^{\vee} &\cong  \big(H^{-d}R\Hom_R(R\Gamma(X,\cF),\omega^{\mydot}_R)\big)^{\wedge} \\
&\cong \big(H^{-d}R\Gamma \circ  R\cHom_{\cO_X}(\cF,\omega^{\mydot}_X)\big)^{\wedge} \cong \big(H^{-d}R\Hom_{\cO_X}(\cF,     \omega_X^{\mydot})\big)^{\wedge},
\end{alignat*}
where the first one is Matlis duality ({Proposition \ref{prop:Matlis-duality}}) 
and \cite[\href{https://stacks.math.columbia.edu/tag/0A06}{Tag 0A06}]{stacks-project},  whilst the second one is Grothendieck duality. If $X$ is Cohen-Macaulay, then $\omega_X^{\mydot} = \omega_X[d]$ and the statement of the lemma follows immediately. However, we are taking the lowest cohomology $H^{-d}$, and so in general we can argue that the higher cohomologies of $\omega_X^{\mydot}$ do not interfere (we refer to the proof of \cite[Lemma 2.2]{BMPSTWW20} for details).
\end{proof}

\begin{lemma} \label{lem:dual-of-vanishing}
Let $X = \Spec R$ for an excellent domain $R$ admitting a dualising complex, let $f \colon Y \to X$ be a proper birational morphism from an integral normal divisorially Cohen-Macaulay scheme $Y$, and let $L$ be a Weil divisor on $Y$. Suppose that
\[
f_*\cO_Y(K_Y+L) \text{ is Cohen-Macaulay,} \, \text{ and } \ R^if_*\cO_Y(K_Y+L) = 0 \text{ for } i>0.
\]
Then $R^if_*\cO_Y(-L)=0$ for every $i>0$.
\end{lemma}
\begin{proof}
To prove this lemma, we may assume that $(R,\m)$ is a  local domain with maximal ideal $\m$. Set $d := \dim R$. 

It is enough to show that {the complex} 
\begin{align*}
Rf_*\cO_Y(-L)&= Rf_*\cHom (\cO_Y(K_Y+L),\omega_Y) \\
&= Rf_*R\cHom(\cO_Y(K_Y+L),\omega_Y)\\
&\cong R\cHom(Rf_*\cO_Y(K_Y+L), \omega^{\mydot}_X[-d])
\end{align*}
is supported in degree $0$. Here, the first equality $\cO_Y(-L)=\cHom(\cO_Y(K_Y+L),\omega_Y)$ can be checked on the regular locus ({see Remark \ref{remark:hom-reflexive}}) in which case it is clear, the second equality follows from 
{\cite[Theorem 3.3.10(d)]{BH93} and }
the fact that $Y$ is divisorially Cohen-Macaulay (and so {$Y$ and} $\cO_Y(K_Y+L)$ are  Cohen-Macaulay), and the 
{isomorphism} {at the end} is Grothendieck duality.

Now, Matlis duality ({Proposition \ref{prop:Matlis-duality}}) yields: 
\[
R\cHom(Rf_*\cO_Y(K_Y+L), \omega^{\mydot}_X)^{\wedge} \cong (R\Gamma_{\mathfrak{m}}Rf_*\cO_Y(K_Y+L))^{\vee},
\]
and so it is enough to show that 
$R\Gamma_{\mathfrak{m}}Rf_*\cO_Y(K_Y+L)$ is supported in degree $d$. Indeed, then
\[
\cH^i(R\cHom(Rf_*\cO_Y(K_Y+L), \omega^{\mydot}_X))^{\wedge} = 0,
\]
for $i \neq {-d}$ by \cite[\href{https://stacks.math.columbia.edu/tag/0A06}{Tag 0A06}]{stacks-project}, where $(-)^{\wedge}$ denotes {here} the usual $\m$-adic completion. Thus $\cH^i(R\cHom(Rf_*\cO_Y(K_Y+L), \omega^{\mydot}_X)) = 0$ as the operation of taking completion on a local ring is fully faithful (combine  \cite[\href{https://stacks.math.columbia.edu/tag/00MA}{Tag 00MA}]{stacks-project} and \cite[\href{https://stacks.math.columbia.edu/tag/00MC}{Tag 00MC}]{stacks-project}).

By assumptions, $Rf_*\cO_Y(K_Y+L)$ is a Cohen-Macaulay coherent sheaf (supported in degree zero), and so
\[
R\Gamma_{\mathfrak{m}}Rf_*\cO_Y(K_Y+L)
\]
is supported in degree $d$ (see \cite[\href{https://stacks.math.columbia.edu/tag/0AVZ}{Tag 0AVZ}]{stacks-project}), 
concluding the proof.
\end{proof}

{\cora Now, we can show the following dual variant of Kodaira/Kawamata-Viehweg vanishing.}

\begin{corollary}
Let $f \colon Y \to X$ be a proper birational morphism of normal surfaces and 
let $(Y, \Delta)$ be a klt pair. 
If $L$ is a Weil divisor such that $L-\Delta$ is $f$-nef and $f_*\MO_Y(K_Y+L)$ is  reflexive, 
then $R^if_*\MO_Y(-L)=0$ for every $i>0$. 
\end{corollary}
\noindent {Note that $Y$ is automatically $\bQ$-factorial by \cite[Theorem 3.3]{tanaka16_excellent}.}
\begin{proof}
It follows from 
\cite[Theorem 3.3]{tanaka16_excellent} that $R^if_*\cO_Y(K_Y+L) = 0$ for $i>0$. 
Hence, the assertion holds by  Lemma \ref{lem:dual-of-vanishing}. 
\end{proof}

\subsection{Generalities on restrictions}

One of the key strategies used in many proofs in birational geometry is to restrict a line bundle to a prime divisor and lift sections therefrom. However, this becomes more complicated when we replace the line bundle by a divisorial sheaf (cf.\ \cite[Proposition 3.1]{HW17} or \cite[Subsection 4.1]{HW19b}). 
The following result is independent of the other parts of the paper, and likely useful in other applications. 


We start by defining a restriction of a $\bQ$-divisor.
\begin{definition} \label{definition:restricting-divisor}
Let $X$ be an integral normal Noetherian scheme, let $S$ be a {\cora normal} prime divisor, and let $D$ be a $\bQ$-Cartier $\bQ$-divisor such that $S \not\subseteq \Supp D$. 
For a positive integer $m$ such that $mD$ is Cartier, we define 
\[
D|_S := \frac{1}{m} ( (mD)|_S), 
\]
which is a $\Q$-Cartier $\Q$-divisor on $S$. 
\end{definition}

\noindent Note that $D|_S$ is independent of the choice of $m$. 
Moreover, the reader should be aware that even when $D$ is Weil, $D|_S$ itself need not be Weil.

\begin{proposition} \label{prop:basic_restriction}
Let $X$ be an integral normal divisorially Cohen-Macaulay excellent scheme admitting a dualising complex and let $S$ be a normal prime divisor. Let $D$ be a $\bQ$-Cartier  $\bQ$-divisor 
such that $S \not\subseteq \Supp D$ and 
$(X,S+\{D\})$ is plt. Then 
there exists a unique $\MO_X$-module homomorphism 
\[
\res \colon \cO_X(\lfloor D \rfloor) 
\to \cO_S(\lfloor D|_S \rfloor) 
\]
such that $\res|_{X\, \setminus\, \Supp D}$ coincides with 
the restriction homomorphism $\rho \colon \MO_{X \,\setminus\, \Supp D} \to \MO_{S\, \setminus\, \Supp D}$. 
Moreover, $\res$ is surjective. 
\end{proposition}

\begin{proof}
The uniquenss of $\res$ is clear {as $X$ and $S$ are integral}. 
Then the problem is local, 
hence we may assume that $X$ is affine and $mD = {\rm div}(f)$ for some $m \in \Z_{>0}$ and $f \in K(X) \setminus \{0\}$.

\setcounter{step}{0} 
\begin{step}\label{s1:basic_restriction}
We prove that there exists an $\MO_X$-module homomorphism $\res$ such that 
$\res|_{X \setminus \Supp D} =\rho$. 
\end{step}


{To this end, take a section} $g \in \Gamma(X, \cO_X(\lfloor D \rfloor))$,
and define its restriction 
\[
g|_S := \rho(g|_{ X\, \setminus \,\Supp D}) \in \Gamma(S\,\! \setminus\,\! \Supp D, \MO_{S\, \setminus \,\Supp D}).
\]
It suffices to show that ${\rm div}(g|_S) \geq - \rdown{D|_S}$ {as this implies that} $g|_S \in \Gamma(S, \MO_S( \rdown{D|_S}))$. 
We have 
\[
\mathrm{div}(g) \geq -\lfloor D \rfloor \geq -D.
\]
Then $\mathrm{div}(g^m) \geq -\mathrm{div}(f)$, where $mD = \mathrm{div}(f)$ locally, and so 
\[
\mathrm{div}(g^m|_S) \geq -\mathrm{div}(f|_S) = -mD|_S.
\]
Hence $\mathrm{div}(g|_S) \geq \lceil -D|_S \rceil = - \lfloor D|_S \rfloor$, as required.  

\begin{step}\label{s2:basic_restriction}
We show that $\res$ is surjective at codimension $2$ points of $X$.
\end{step}

To this end, we can localise and assume that $\dim X = 2$. 
Let $f \colon Y \to X$ be a log resolution of $(X,S+D)$. 
Set $S_Y := f^{-1}_*S$ and $D_Y := f^*D$. 
We define a $\Q$-divisor $B_Y$ by $K_Y+S_Y+B_Y= f^*(K_X+S +\{D\})$. 
{Since} $\rdown{D_Y|_{S_Y}} = \rdown{D_Y}|_{S_Y}$, we have the following short exact sequence
\[
0 \to \cO_Y(\lfloor D_Y - S_Y \rfloor) \to \cO_Y(\lfloor D_Y \rfloor) \to \cO_{S_Y}(\lfloor D_Y|_{S_Y} \rfloor) \to 0.
\]
We claim that \[
R^1f_*\cO_Y(\lfloor D_Y - S_Y \rfloor) = 0 \quad 
\text{ and } \quad f_*\cO_Y(\lfloor D_Y \rfloor) = \cO_X(\lfloor D \rfloor).
\]
{Assuming this claim, Step 2 is concluded as the assertion} $(f|_{S_Y})_*\cO_{S_Y}(\lfloor D_Y|_{S_Y}\rfloor) = \cO_S(\lfloor D|_S \rfloor)$ {is automatic. Indeed,} 
$\dim S_Y=1$, which gives that $f|_{S_Y} \colon S_Y \to S$ is an isomorphism.

By $D_Y = f^*D$, we have 
$f_*\cO_Y(\lfloor D_Y \rfloor) =f_*\MO_Y(D_Y) = \MO_X(D) = \cO_X(\lfloor D \rfloor)$. 
To show 
\[
R^1f_*\cO_Y(\lfloor D_Y - S_Y \rfloor) = 0,
\]
we apply Lemma \ref{lem:dual-of-vanishing} for $L = \rup{S_Y -D_Y}$. 
To check its assumptions, we note that $(Y,S_Y + \{D_Y\})$ is plt as $f$ is a log resolution of $(X,S+D)$, and so
\begin{align*}
R^if_*\cO_Y(K_Y + \rup{S_Y-D_Y}) &= R^if_*\cO_Y(K_Y + S_Y + \{D_Y\} - D_Y) \\
&= R^if_*\cO_Y(K_Y + S_Y + \varepsilon G + \{D_Y\} - D_Y - \varepsilon G) \\
&=  0
\end{align*}
for $i>0$ by the Kawamata--Viehweg vanishing theorem for birational morphisms of surfaces (\cite[Theorem 10.4]{kollar13}), where $G$ is an effective $f$-exceptional $f$-anti-ample $\bQ$-divisor and $0 < \varepsilon \ll 1$.
Moreover,
\begin{equation}\label{e1:basic_restriction}
f_*\cO_Y(K_Y + \lceil S_Y-D_Y \rceil) = \cO_Y(K_X + \lceil S -D \rceil).
\end{equation}
Indeed, we have that
\begin{align*}
K_Y + \lceil S_Y -D_Y \rceil - f^*(K_X + \lceil S - D \rceil) &= K_Y+S_Y + \{D_Y\} - f^*(K_X+S +\{D\}) \\
&= -B_Y + \{D_Y\}.
\end{align*}
The round-up of this $\bQ$-divisor is effective, 
because $\lceil -B_Y \rceil \geq 0$ follows from the fact that $(X,S+\{D\})$ is plt. This yields  (\ref{e1:basic_restriction}) (cf.\ \cite[Lemma 2.31]{BMPSTWW20}).

In particular, $f_*\cO_Y(K_Y + \lceil S_Y-D_Y \rceil)$ is reflexive, and so Cohen-Macaulay as $X$ is a normal two-dimensional scheme. This concludes the verification of the assumptions of Lemma \ref{lem:dual-of-vanishing}, and so of the proof of Step \ref{s2:basic_restriction}.

\begin{step}\label{s3:basic_restriction}
We finish the proof of Proposition \ref{prop:basic_restriction}. 
\end{step}

By Step \ref{s2:basic_restriction}, 
we may assume that $\dim X \geq 3$. Consider the following exact sequence
\[
0 \to \cO_X(\lfloor D - S \rfloor) \to \cO_X(\lfloor D \rfloor) \to \cE \to 0,
\]
{defining} $\cE$ as the cokernel of $\cO_X(\lfloor D - S \rfloor) \hookrightarrow \cO_X(\lfloor D \rfloor)$. 
Since $X$ is divisorially Cohen-Macaulay, we have that $\cE$ is a reflexive $\MO_S$-module 
by \cite[Lemma 2.60]{kollar13} (cf.\  \cite[Lemma 2.4]{HW17}). Moreover, by Step 2, {the sheaf} $\cE$ agrees with $\cO_S(\lfloor D|_S \rfloor)$ on codimension one points of $S$. Since both are reflexive sheaves on $S$, this implies that $\cE \cong \cO_S(\lfloor D|_S \rfloor)$, concluding the proof. 
\end{proof}

\subsection{Properties of two-dimensional plt and lc pairs}
The following lemma will be used throughout the article.

\begin{lemma} \label{lem:everything-about-plt-surface-pairs}
Let $(S,C+B)$ be a two-dimensional plt pair, where $S= \Spec R$ {is a spectrum of a local ring}, 
$C$ is a prime divisor, and $B=\sum a_i E_i$ for prime divisors  $E_i$ and rational numbers $\frac{1}{2} \leq a_i < 1$. 
Let $P$ be the closed point of $S$.

Define $n := \det(\Gamma_P)$, where $\Gamma_P$ is the intersection matrix of the minimal resolution of $P \in S$. 
We set $n=1$ if $P \in S$ is regular. {Then the} following hold.
\begin{enumerate}
    \item $P \in S$ is a regular point or a cyclic singularity, that is, the dual graph of the exceptional divisors of the minimal resolution is a chain.
    \item $(K_S+C)|_C= K_C+ (1-\frac{1}{n})P$.
    \item For every Weil divisor $D$ on $S$, we have that $nD$ is Cartier.
    \item If $B \neq 0$, then there exists a prime divisor $E$ such that
        \begin{enumerate}
            \item $B=aE$ for some rational number $\frac{1}{2} \leq a <1$, and
            \item $(S,C+E)$ is log canonical. 
        \end{enumerate}
    \item $(K_S+C+B)|_C = K_C + (1-\frac{1-a}{n})P$, where we set $a=0$ if $B=0$.
\end{enumerate}
\end{lemma}
\noindent Note that $S$ is automatically $\bQ$-factorial by \cite[Corollary 4.11]{tanaka16_excellent} and $C$ is normal by \cite[3.35]{kollar13}. We emphasise that, by definition, prime divisors are closed and non-zero, and so they must pass through $P$.

\begin{proof}
Assertions (1) and  (2) follow from \cite[Theorem 3.36]{kollar13}. We now show (3). 
Let $\pi \colon W \to X$ be the minimal resolution of $X$. 
It follows from the same argument as in \cite[Lemma 2.2]{CTW15b}  
that $n \pi^*D$ is a Cartier divisor on $W$. 
Note that the Cartier divisor $n\pi^*D = K_W + \Delta_W -(K_W+\Delta_W) + n\pi^*D$ is semiample  for $K_W + \Delta_W = \pi^*K_S$. 
By the base point free theorem (\cite[Theorem 4.2]{tanaka16_excellent}), 
$bnD$ is a Cartier divisor for any sufficiently large integer $b \gg 0$, and so 
the divisor $nD = (b+1)nD -bnD$ is Cartier. Thus (3) holds.
As for {Assertions (4) and (5)}, since $(S,C+B)$ is plt, the cases (1) and (2) from \cite[Corollary 3.45]{kollar13} cannot happen, and so we are in the case (3) of op.cit.\ which states that our {A}ssertions (4a) and (5) are valid.

In what follows, we show (4b). To this end, observe that by (2) and (5)
\begin{align*}
aE|_C &= (K_S+C+aE)|_C - (K_S+C)|_C \\ &= K_C+  \Big(1-\frac{1-a}{n}\Big)P - \Big(K_C + \Big(1-\frac{1}{n}\Big)P\Big)\\
&= \frac{a}{n}P.
\end{align*}
Therefore, $E|_C = \frac{1}{n}P$, and so $(K_S+C+E)|_C = (K_S+C)|_C + E|_C =  K_C + P$ by (2). In particular, $(S,C+E)$ is log canonical by inversion of adjunction (\cite[Theorem 5.1(1)]{tanaka16_excellent}). 
\end{proof}

For future reference, we observe that in (5) of the above lemma, we have that
\begin{equation} \label{eq:coefficients-increase-in-adjunction}
1 - \frac{1-a}{n} \geq a,
\end{equation}
and so the coefficients of the different {$\Diff_C(B)$} of $(S,C+B)$ are greater {than} or equal to the coefficients of $B$. Moreover, Lemma \ref{lem:everything-about-plt-surface-pairs} shows that if $B$ has standard coefficients, then $B=(1-\frac{1}{m})E$ for {some} $m \in \bZ_{>0}$ (we set $m=1$ if $B=0$), and
\begin{equation} \label{eq:qdlt-discrepancy}
(K_S+C+B)|_C = K_C + \Big(1-\frac{1}{nm}\Big)P.
\end{equation}

{ 
\begin{lemma} \label{lem:surface-lc-pair-with-high-coeffs-is-plt}
{Let $(S,C+B)$ be a two-dimensional lc pair, where $S= \Spec R$ {is a spectrum of a local ring}, 
$C$ is a prime divisor, and $B=\sum_{i \in I} a_i E_i$ for prime divisors  $E_i$ and rational numbers $\frac{1}{2} \leq a_i < 1$.}
Suppose that $a_j>\frac{1}{2}$ for some $j \in I$. Then $(S,C+B)$ is plt.
\end{lemma}
\begin{proof}
By {applying} Lemma \ref{lem:everything-about-plt-surface-pairs}(4a) to 
{a plt pair} $(S,C+B-\varepsilon E_j)$ for $0 < \varepsilon \ll 1$, we get that $B = aE$ for a prime divisor $E$ and a rational number $\frac{1}{2} < a < 1$. The assumptions of \cite[Corollary 3.45]{kollar13} are satisfied, and since Case (1) and (2) in op.cit.\ cannot hold as $a>\frac{1}{2}$, we must satisfy Case (3). 
By inversion of adjunction \cite[Theorem 5.1(2)]{tanaka16_excellent} {and the fact that $a<1$, the pair} $(S,C+B)$ is plt. 
\end{proof}}

For the sake of future references, we state the following fact.
\begin{remark} \label{remark:ACC-P^1}
A sum of numbers in the set $\{1-\frac{1}{n}\mid {n \in \Z_{>0}}
\} \cup (\frac{5}{6},1 ]$ is never equal to $2$ if at least one of them {belongs} to $(\frac{5}{6},1)$. This follows by a simple calculation (note that $t + \frac{1}{2} + \frac{1}{2} < 2 < t + \frac{1}{2} + \frac{2}{3}$ when $\frac{5}{6} < t < 1$). 
\end{remark}

In Section \ref{section-log-del-Pezzo}, we often use \cite[Corollary 3.44(1)]{kollar13}, which {says} that log canonical thresholds on surfaces cannot lie in the set $(\frac{5}{6},1)$. {We provide a statement and a proof in a slightly more general case of log pairs with standard coefficients. {In what follows}, we denote the set of coefficients of a $\bQ$-divisor $\Delta$ by $\mathrm{coeff}(\Delta)$.}

\begin{lemma}  \label{lem:kollar-acc} Let $(X,\Delta)$ be a two-dimensional lc pair. {Suppose that $\Delta=\Delta_1+\Delta_2$, where $\Delta_1$ and $\Delta_2$ are effective $\bQ$-divisors {without common components} and such that $\mathrm{coeff}(\Delta_1) \subseteq \{1-\frac{1}{n} \mid 1 \leq n \leq 6\}$ and $\mathrm{coeff}(\Delta_2) \subseteq (\frac{5}{6},1]$. Then $(X, \Delta_1 + \Supp \Delta_2)$ is log canonical.}
\end{lemma}

\begin{proof}
{Write $\Delta = tC+B$ for $\frac{5}{6} < t < 1$, a prime divisor $C$, and an effective $\Q$-divisor $B$ such that $C \not\subseteq \Supp B$ and $\mathrm{coeff}(B) \subseteq \{1-\frac{1}{n} \mid 1 \leq n \leq 6\} \cup (\frac{5}{6},1]$. We claim that 
\begin{equation} \label{eq:lct-claim}
(X,C+B) \quad \textrm{is log canonical. }
\end{equation}
Assuming the claim and subsequently repeating this procedure for $(X,C+B)$ and another prime divisor in $\Supp \Delta_2$, the proof of the lemma is concluded.}


{Therefore, it is enough to show Claim (\ref{eq:lct-claim})}. The question is local, so we may assume that $X= \Spec R$ for a local ring $R$ with maximal ideal $\m$. 
Suppose that $(X,C+B)$ is not log canonical. {Then, after} replacing $t$ by the log canonical threshold of $(X,B)$ with respect to $C$, we may assume that  $\m$ 
is a log canonical centre of $(X, \Delta := tC+B)$ {and still $\frac{5}{6}<t<1$}. 
By a standard argument, 
we can construct a projective birational morphism $f \colon Y \to X$ from a two-dimensional integral normal $\Q$-factorial excellent scheme $Y$ such that $E:= \Exc(f)$ is a prime divisor and 
\[
K_Y + E+f_*^{-1}\Delta = f^*(K_X+\Delta)
\footnote{ 
Such a birational morphism $f\colon Y \to X$ can be constructed as follows: 
we first take a dlt blowup $g \colon Z \to X$ of $(X, \Delta)$ \cite[Theorem 4.7]{tanaka16_excellent}, so that $(Z, \Exc(g)+g_*^{-1}\Delta)$ is dlt and $K_Z +\Exc(g) +g_*^{-1}\Delta =g^*(K_X+\Delta)$. 
For $0 < \varepsilon \ll 1$, we pick $0 < \delta \ll 1$ such that 
$(Z, (1+\delta)\Exc(g)  + g_*^{-1}\Delta -\varepsilon g^*\Delta)$ is klt. 
{Since} $K_Z + (1+\delta)\Exc(g)  + g_*^{-1}\Delta -\varepsilon g^*\Delta \equiv_{g} 
K_Z+(1+\delta)\Exc(g)  + g_*^{-1}\Delta$, 
we may run a  $(K_Z+(1+\delta)\Exc(g)  + g_*^{-1}\Delta)$-MMP over $X$. 
By the negativity lemma, $X$ itself is the end result of this MMP, and 
so $Y$ can be chosen as the last step of this MMP.}. 
\]
Then $(Y,E+f^{-1}_*\Delta)$ is lc, because so is $(X, {\Delta})$. 

Write 
\[
(K_Y + E + f^{-1}_*\Delta)|_E = K_E + \Delta_E \quad {\rm and} \quad  \Delta_E = \sum_{i \in I} a_i P_i
\] 
for rational numbers $0 < a_i  \leq 1$ and distinct closed  points $P_i \in E$.  
For every $i \in I$, we have that $a_i \in \{1 - \frac{1}{n} \mid  1 \leq n \leq 6\} \cup (\frac{5}{6}, 1]$; indeed, if $P_i$ is a non-plt point of $(Y,E+f^{-1}_*\Delta)$, then $a_i=1$ by inversion of adjunction, otherwise apply Lemma \ref{lem:everything-about-plt-surface-pairs}(5), (\ref{eq:coefficients-increase-in-adjunction}), and  (\ref{eq:qdlt-discrepancy}). 

Since $f^{-1}_*C$ intersects $E$, we may pick a closed point $P_j \in f^{-1}_*C \cap E$ for some $j \in I$. 
As $(Y,E+f^{-1}_*\Delta)$ is {lc} at $P_j$ and $\Delta \geq tC$ for $t > \frac{1}{2}$, Lemma \ref{lem:surface-lc-pair-with-high-coeffs-is-plt} shows that $(Y,E+f^{-1}_*\Delta)$ is plt at $P_j$. 
Therefore, $\frac{5}{6} < a_j < 1$ by Lemma \ref{lem:everything-about-plt-surface-pairs}(5) and (\ref{eq:coefficients-increase-in-adjunction}). Furthermore, $K_E +\Delta_E \sim_{\bQ} 0$, so $\deg \Delta_E = 2$ and $\sum_{{i \in I}} a_i = 2$. This contradicts Remark \ref{remark:ACC-P^1}.    
\end{proof}

%% file: section3.tex
\section{Foundations of log quasi-$F$-splittings} \label{s:foundations}

This section is devoted to defining log quasi-$F$-splittings and the study of their properties. However,  before delving  into the technicalities of this section, we recommend the reader to consult Subsection \ref{ss:calculation-local-cohomology} in the Appendix, where we provide an example of how to check whether an explicitly defined singularity is quasi-$F$-split.

\subsection{Definition of log quasi-$F$-splitting}\label{ss-log-QFS}

Until the end of this section, we work under the assumptions of Setting \ref{setting:most-general-foundations-of-log-quasi-F-splitting}.

\begin{setting} \label{setting:most-general-foundations-of-log-quasi-F-splitting}
Let $X$ be a $d$-dimensional 
integral normal {Noetherian} 
$F$-finite 
$\F_p$-scheme admitting a dualising complex $\omega^{\mydot}_X$. 
\end{setting}
{In some cases, we might need to assume that our scheme is projective over a fixed affine base as described by Setting \ref{setting:foundations-of-log-quasi-F-splitting} to which we will refer whenever necessary.} 

\begin{setting} \label{setting:foundations-of-log-quasi-F-splitting}
Let $R$ be {an $F$-finite Noetherian domain of characteristic $p>0$}. 
Let $X$ be a $d$-dimensional integral normal scheme which is projective over $\Spec R$. 
\end{setting}

\begin{remark}
{In the situation of} Setting \ref{setting:most-general-foundations-of-log-quasi-F-splitting}, 
it {automatically} holds that $X$ is excellent (Subsection \ref{ss:notation}(\ref{sss-Ffinite})) and 
$d = \dim X <\infty$ \cite[Ch.\ V, Corollary 7.2]{Har66}. 
{In the situation of} Setting \ref{setting:foundations-of-log-quasi-F-splitting}, 
we have that $\dim R<\infty$, $R$ is excellent, and $R$ admits a dualising complex (Subsection \ref{ss:notation}(\ref{sss-Ffinite})). 
\end{remark}

Let $\Delta$ be a (non-necessarily effective) $\bQ$-divisor on $X$. We define $Q_{X,\Delta,n}$ and $\Phi_{X, \Delta, n}$ by the following pushout diagram of $W_n\MO_X$-module homomorphisms: 
\begin{equation} \label{diagram:quasi-F-split-definition}
\begin{tikzcd}
W_n\cO_X(\Delta) \arrow{r}{F} \arrow{d}{R^{n-1}} & F_* W_n \cO_X(p\Delta)  \arrow{d}\\
\cO_X(\Delta) \arrow{r}{\Phi_{X, \Delta, n}}& \arrow[lu, phantom, "\usebox\pushoutdr" , very near start, yshift=0em, xshift=0.6em, color=black] Q_{X, \Delta, n}.
\end{tikzcd}
\end{equation}
We remind the reader that 
$\MO_X(\Delta) = \MO_X(\rdown{\Delta})$, 
but it is \emph{not} true in general that $W_n\cO_X(\Delta) = W_n\cO_X(\rdown{\Delta})$.  

We define a $W_n\cO_X$-module 
\begin{align}
\label{eq:definition-of-log-B} B_{X, \Delta, n} &:= {\rm Coker}( W_n\MO_X(\Delta) \xrightarrow{F} F_*W_n\MO_X(p\Delta))\\
&\hphantom{:}= F_*W_n\MO_X(p\Delta)/F(W_n\MO_X(\Delta)) \nonumber,
\end{align}
{and set $B_{X,\Delta} := B_{X,\Delta,1}$ to simplify the notation.}
\begin{remark}
{The key properties of the construction of $Q_{X,\Delta,n}$,  $B_{X,\Delta,n}$, and $\Phi_{X,\Delta,n}$ may be encapsulated by the following diagram}
\begin{equation}\label{e-big-BC-diagram}
\begin{tikzcd}
& 0 \arrow{d} & 0 \arrow{d} & & \\
&F_*W_{n-1}\MO_X(p\Delta) \arrow{d}{V} \arrow[r,dash,shift left=.1em] \arrow[r,dash,shift right=.1em] & F_*W_{n-1}\MO_X(p\Delta) \arrow{d}{FV} & & \\
0 \arrow{r} & W_n\MO_X(\Delta) \arrow{r}{F} \arrow{d}{R^{n-1}} & F_*W_n\MO_X(p\Delta) \arrow{r} \arrow{d} & B_{X, \Delta, n} \arrow[d,dash,shift left=.1em] \arrow[d,dash,shift right=.1em] \arrow{r}  & 0 \\
0 \arrow{r} & \MO_X(\Delta) \arrow{r}{\Phi_{X, \Delta, n}} \arrow{d} &  \arrow[lu, phantom, "\usebox\pushoutdr" , very near start, yshift=0em, xshift=0.6em, color=black] Q_{X, \Delta, n} \arrow{r} \arrow{d} & B_{X, \Delta, n} \arrow{r} & 0.\\
& 0 & 0 & &
\end{tikzcd}
\end{equation}
{All  the horizontal and vertical sequences are exact, as} $F \colon W_n\MO_X(\Delta) \to F_*W_n\MO_X(p\Delta)$ is injective, $R^{n-1} \colon W_n\MO_X(\Delta) \to \MO_X(\Delta)$ is surjective, and 
(\ref{diagram:quasi-F-split-definition}) is a pushout diagram. 
\end{remark}

In what follows, we shall often assume that $\rdown{\Delta}=0$, which is equivalent to 
the condition that the coefficients of $\Delta$ are contained in $[0, 1)$. 
In this case, $\cO_X(\Delta) =\MO_X$. 
\begin{definition}
\label{definition:log-quasi-F-split}
In the situation of Setting  \ref{setting:most-general-foundations-of-log-quasi-F-splitting}, let $\Delta$ be a $\bQ$-divisor on $X$ such that $\rdown{\Delta}=0$. 
We say that $(X,\Delta)$ is \emph{$n$-quasi-$F$-split} if there exists 
a $W_n\MO_X$-module homomorphism $\alpha \colon F_*W_n\cO_X(p\Delta) \to \cO_X$ 
such that  $R^{n-1} = \alpha \circ F$:

\begin{equation*} 
\begin{tikzcd}
W_n\cO_X(\Delta) \arrow{r}{F} \arrow{d}{{R^{n-1}}} & F_* W_n\cO_X(p\Delta) \arrow[dashed]{ld}{\alpha} \\
\cO_X.
\end{tikzcd}
\end{equation*}
We call $(X,\Delta)$ \emph{quasi-$F$-split} if it is $n$-quasi-$F$-split for some $n \in \Z_{>0}$.
{We say that $(X,\Delta)$ is \emph{$n$-quasi-$F$-pure (quasi-$F$-pure, resp.)} if it is $n$-quasi-$F$-split (quasi-$F$-split, resp.) affine-locally.} 
\end{definition}

\begin{proposition}[{cf.\ Propostion \ref{prop:intro-properties-of-C}}] \label{prop:properties-of-C}
In the situation of Setting \ref{setting:most-general-foundations-of-log-quasi-F-splitting}, let $\Delta$ be a $\bQ$-divisor on $X$. 
Then the following hold. 
\begin{enumerate}
    \item The map $F_*W_n\cO_X(p\Delta) \to Q_{X,\Delta,n}$ induces an isomorphism of $W_n\MO_X$-modules 
    \[
    Q_{X,\Delta,n} \cong \frac{F_*W_n\MO_X({ p}\Delta)}{pF_*W_n\MO_X({ p}\Delta)}.
    \] 
    \item 
    $Q_{X,\Delta,n}$ is naturally an $\cO_X$-module. In particular, {it} is a coherent $\MO_X$-module and $\Phi_{X, \Delta, n} \colon \MO_X \to Q_{X, \Delta, n}$ is an $\MO_X$-module homomorphism. 
    \item $B_{X, \Delta, n}$ is naturally an $\cO_X$-module. In particular, {it} is a coherent $\MO_X$-module. 
\end{enumerate}
\end{proposition}

\begin{proof}
Assertion (1) holds by  
\[
Q_{X,\Delta,n} \cong \frac{F_*W_n\MO_X({ p}\Delta)}{FV(F_*W_{n-1}\MO_X(p\Delta)) } = \frac{F_*W_n\MO_X({ p}\Delta)}{pF_*W_{n}\MO_X({ p}\Delta) }, 
\]
where the isomorphism follows from (\ref{e-big-BC-diagram}) and 
we get the equality {by the identity $FV = p$ in $W_nK(X)$; specifically, we have following factorisation with $R$ surjective}
\[
p \colon F_*W_n\MO_X({p}\Delta ) \xrightarrow{R} F_*W_{n-1}\MO_X({ p}\Delta ) \xrightarrow{V} 
W_{n}\MO_X(\Delta ) \xrightarrow{F} F_*W_n\MO_X({ p}\Delta ). 
\]

Let us show (2). 
By (1), it suffices to show 
the inclusion 
\begin{equation}\label{e:properties-of-C}
    V(W_{n-1}\cO_X) \cdot (F_*W_n\MO_X({p}\Delta)) \subseteq p F_*W_n\MO_X({p}\Delta).
\end{equation}
Pick an open subset $U \subseteq X$ and  $s \in W_{n-1}\cO_X(U)$. 
Then the required inclusion (\ref{e:properties-of-C}) follows from 
\[
(Vs) \cdot F_*W_n\MO_X({p}\Delta)(U) = F_*( (FVs) \cdot W_n\MO_X({ p}\Delta)(U)) 
\]
\[
= p F_*( s \cdot W_n\MO_X({p}\Delta)(U)) 
\subseteq 
p F_*W_n\MO_X({p}\Delta)(U). 
\]
Thus, (2) holds. Assertion (3) follows from (2) and the surjection $Q_{X, \Delta, n} \to B_{X, \Delta, n}$ in (\ref{e-big-BC-diagram}). 
\end{proof}

\begin{proposition} \label{prop:definition-via-splitting}
In the situation of Setting \ref{setting:most-general-foundations-of-log-quasi-F-splitting}, let $\Delta$ be a $\bQ$-divisor on $X$ such that $\lfloor \Delta \rfloor =0$. Then 
the following are equivalent. 
\begin{enumerate}
    \item $(X,\Delta)$ is $n$-quasi-$F$-split. 
    \item $\Phi_{X, \Delta, n}\colon\cO_X \to Q_{X,\Delta,n}$ splits as a $W_n\MO_X$-module homomorphism. 
    \item $\Phi_{X, \Delta, n} \colon \cO_X \to Q_{X,\Delta,n}$ splits as an $\MO_X$-module homomorphism. 
\end{enumerate}
\end{proposition}

\begin{proof}
Since (\ref{diagram:quasi-F-split-definition}) is a pushout diagram of $W_n\MO_X$-modules, 
(1) and (2) are equivalent. 
By Proposition \ref{prop:properties-of-C}, (2) and (3) are equivalent. 
\end{proof}

One can now check that 
\begin{equation}\label{e-C-tensor}
Q_{X,\Delta+D,n} \cong Q_{X,\Delta,n} \otimes_{\MO_X} \cO_X(D)
\end{equation} 
for any Cartier divisor $D$. 
Indeed, we have $Q_{X,\Delta,n} \otimes_{\MO_X} \cO_X(D) \cong 
Q_{X,\Delta,n} \otimes_{W_n\MO_X} W_n\cO_X(D)$ (Proposition \ref{prop:properties-of-C}{\cora (2)}), and hence 
both sides of (\ref{e-C-tensor}) are isomorphic to the pushout of 
$F_*W_n\MO_X(p\Delta+pD) \xleftarrow{F} W_n\MO_X(\Delta+D) \xrightarrow{R^{n-1}} \MO_X(\Delta+D) $. 

We now state analogues of short exact sequences (\ref{eq:intro-C-quotient-sequence}), (\ref{eq:intro-Bn-sequence}), and (\ref{eq:intro-C-restriction-sequence}) for when $\Delta$ is an arbitrary $\bQ$-divisor on $X$.

\begin{lemma} \label{lem:key-sequences-in-the-log-case}
In the situation of Setting \ref{setting:most-general-foundations-of-log-quasi-F-splitting}, let $\Delta$ be a $\bQ$-divisor on $X$. 
Take $n \in \Z_{>0}$. 
Then we have the following exact sequences of coherent $\MO_X$-modules: 
\begin{align} \label{eq:C-quotient-sequence}
&0 \to \cO_X(\Delta) \xrightarrow{\Phi_{X, \Delta, n}}\, Q_{X,\Delta,n} \to \cB_{X,\Delta,n} \to 0 \\
\label{eq:Bn-sequence}
&0 \to F_* \cB_{X,p\Delta,n} \to \cB_{X,\Delta,n+1} \to {\cB_{X,\Delta}} \to 0 \\
\label{eq:C-restriction-sequence}
&0 \to 
F_*B_{X, p\Delta, n} \to Q_{X,\Delta,n+1} \to F_* \cO_X(p\Delta) \to 0. 
\end{align}
\end{lemma}

\begin{proof}
By (\ref{e-big-BC-diagram}), we have the exact sequence (\ref{eq:C-quotient-sequence}). 
The exact sequences (\ref{eq:Bn-sequence}) and (\ref{eq:C-restriction-sequence})  
are obtained by applying the snake lemma 
to the following commutative diagrams in which each horizontal sequence is exact
\[
\begin{tikzcd}
0 \arrow{r}& F_*W_{n}\MO_X(p\Delta) \arrow{d}{F}  \arrow{r}{V}& W_{n+1}\MO_X(\Delta) \arrow{d}{F} \arrow{r}{R^n}& \MO_X(\Delta)\arrow{d}{F} \arrow{r}& 0\\
0 \arrow{r}& F^2_*W_{n}\MO_X(p^2\Delta) \arrow{r}{V}& F_*W_{n+1}\MO_X(p\Delta) \arrow{r}{R^n}& F_*\MO_X(p\Delta) \arrow{r}& 0\mathrlap{\textrm{ and }}
\end{tikzcd}
\]
\vspace{0.3em}
\[
\begin{tikzcd}
& F_*W_{n}\MO_X(p\Delta) \arrow{d}{F} \arrow[r,dash,shift left=.1em] \arrow[r,dash,shift right=.1em] & F_*W_{n}\MO_X(p\Delta) \arrow{d}{VF} & &\\
0 \arrow{r}&  F_*^2W_{n}\MO_X(p^2\Delta) \arrow{r}{V} & F_*W_{n+1}\MO_X(p\Delta) \arrow{r}{R^n} & F_*\MO_X(p\Delta) \arrow{r} & 0.
\end{tikzcd}
\]
\end{proof}

\begin{proposition} \label{prop:C-cohen-macaulay} In the situation of Setting \ref{setting:most-general-foundations-of-log-quasi-F-splitting},  suppose that $X$ is divisorially Cohen-Macaulay and let $\Delta$ be a $\bQ$-divisor on $X$ such that $(X,\{p^j\Delta\})$ is naively keenly $F$-pure for every $j \geq 0$. Then $Q_{X,\Delta,n}$ and $\cB_{X,\Delta,n}$ are 
coherent Cohen-Macaulay $\MO_X$-modules for every $n \in \Z_{>0}$.  
\end{proposition}
\begin{proof}
We may assume that $X = \Spec R$ for a local ring $(R,\m)$. It is enough to show that $H^i_\m(Q_{X,\Delta,n})=0$ and $H^i_\m(\cB_{X,\Delta,n})=0$ for {every} $i<\dim {R}$ by  \cite[\href{https://stacks.math.columbia.edu/tag/0AVZ}{Tag 0AVZ}]{stacks-project}. 

Fix an integer  $i < \dim {R}$. 
Since $X$ is divisorially Cohen-Macaulay, we have that $H^i_\m(\cO_X(D))=0$ for any Weil divisor $D$, and so {also} $H^i_\m(F^{j+1}_*\cO_X(p^{j+1}\Delta))=0$ for every $j \geq 0$. Moreover, by Lemma \ref{lem:p-compatible-splitting} and 
the exact sequence 
\[
0 \to \cO_X(p^{j}\Delta) \xrightarrow{F} F_*\cO_X(p^{j+1}\Delta) \to { B_{X,p^j\Delta}} 
\to 0, 
\]
{the sheaf}  $F_*^j{B_{X,p^j\Delta}}$ is locally a direct summand of $F_*^{j+1}\cO_X(p^{j+1}\Delta)$, and thus
\[
H^i_{\m}(F^j_*{B_{X,p^j\Delta}})=0.
\]
Now, by  repeatedly applying (\ref{eq:Bn-sequence}), we get that $H^i_{\m}(\cB_{X,\Delta,n})=0$. Finally (\ref{eq:C-quotient-sequence}) implies $H^i_{\m}(Q_{X,\Delta,n})=0$.
\end{proof}

\noindent In particular, under the same assumptions as in the above proposition, $Q_{X,\Delta,n}$ is a coherent reflexive $\MO_X$-module, 
which in turn shows that 
$Q_{X,\Delta+D,n} \cong (Q_{X,\Delta,n} \otimes \cO_X(D))^{\vee\vee}$ for every Weil divisor $D$.\\

For any $\Q$-divisor $\Delta$ on $X$, 
we define $\Psi_{X, \Delta, n}$  by applying $\cHom_{\cO_X}(-, \cO_X(\Delta))$ to the $\cO_X$-module homomorphism $\Phi_{X, \Delta, n} \colon \cO_X({\Delta}) \to Q_{X,\Delta,n}$:
\begin{equation}\label{eq:dual-defin-quasi-F-split}
\Psi_{X, \Delta, n} \colon 
\cHom_{\MO_X}(Q_{X,\Delta,n}, \cO_X(\Delta)) \to \cHom_{\MO_X}(\cO_X(\Delta), \cO_X(\Delta)) = \cO_X,
\end{equation}
{where the last equality follows from Remark \ref{remark:hom-reflexive}. When $\rdown{\Delta}=0$, we get an $\MO_X$-module homomorphism} $\Psi_{X, \Delta, n} \colon Q^{\vee}_{X,\Delta,n} \to \cO_X$.


\begin{lemma} \label{lem:dual-criterion-for-log-quasi-F-splitting}
In the situation of  
Setting \ref{setting:most-general-foundations-of-log-quasi-F-splitting}, 
let $\Delta$ be an {arbitrary} $\Q$-divisor {on $X$}. 
{Then $\Phi_{X,\Delta,n} \colon \cO_X(\Delta) \to Q_{X,\Delta,n}$ splits if and only if the following map is surjective}
\[
H^0(X, \Psi_{X, \Delta, n}) \colon  \Hom_{\cO_X}(Q_{X,\Delta,n}, \cO_X(\Delta)) \to  H^0(X,\cO_X).
\]
\end{lemma}

\begin{proof}
If $\Phi_{X,\Delta,n} \colon \cO_X(\Delta) \to Q_{X,\Delta,n}$ splits, then 
$H^0(X, \Psi_{X, \Delta, n})$ is a split surjection,  because it is a dual of a split {injection}.

Conversely, if 
$
H^0(X, \Psi_{X, \Delta, n}) \colon  \Hom_{\cO_X}(Q_{X,\Delta,n}, \cO_X(\Delta)) \to  H^0(X,\cO_X)
$
is surjective, then there exists an $\MO_X$-module homomorphism $\alpha \colon Q_{X,\Delta,n} \to \cO_X(\Delta)$ such that 
$H^0(X, \Psi_{X, \Delta, n})(\alpha) = {\rm id}_{\MO_X(\Delta)}$, that is, {the composite homomorphism}
{\[
\cO_X(\Delta) \xrightarrow{\Phi_{X,\Delta,n}} Q_{X,\Delta,n} \xrightarrow{\alpha} \cO_X(\Delta)
\]}
is the identity. {This} implies that
$\Phi_{X,\Delta,n}$ splits. 
\end{proof}

 {In particular, if $\rdown{\Delta}=0$, then $(X,\Delta)$ is $n$-quasi-$F$-split if and only if \[
 H^0(X, \Psi_{X, \Delta, n}) \colon H^0(X,Q^{\vee}_{X,\Delta,n}) \to H^0(X,\cO_X)
 \]
 is surjective. 
 One can easily deduce from this that $X$ is $n$-quasi-$F$-pure if and only if all the stalks $\MO_{X, x}$ are $n$-quasi-$F$-split.}

{Whether $\Phi_{X,\Delta,n}$ splits depends only on the fractional part $\{\Delta\}$ of $\Delta$ as indicated by the following lemma.}
\begin{lemma}\label{l-inj-QFS1}
In the situation of {Setting  \ref{setting:most-general-foundations-of-log-quasi-F-splitting}}, 
let $\Delta$ be an {arbitrary} $\Q$-divisor on $X$ and {let $D$ be a Weil divisor on $X$. Then the following are equivalent. 
\begin{enumerate}
    \item $\Phi_{X, \Delta, n} \colon \cO_X(\Delta) \to Q_{X, \Delta, n}$ splits as an $\cO_X$-module homomorphism. 
    \item $\Phi_{X, D+\Delta, n} \colon \cO_X(D+\Delta) \to Q_{X, D+\Delta, n}$ splits  as an $\cO_X$-module homomorphism. 
\end{enumerate}}
\end{lemma}

\begin{proof}
{Consider the following commutative diagram
 \[
 \begin{tikzcd}[column sep = 6.5em]
 \Hom_{\cO_X}(Q_{X,\Delta,n}, \cO_X(\Delta)) \arrow[d,dash,shift left=.1em] \arrow[d,dash,shift right=.1em] \arrow{r}{H^0(X,\Psi_{X,\Delta,n})} &\Hom_{\cO_X}(\cO_X(\Delta),\cO_X(\Delta)) \arrow[d,dash,shift left=.1em] \arrow[d,dash,shift right=.1em] \\
  \Hom_{\cO_X}(Q_{X,D+\Delta,n}, \cO_X(D+\Delta)) \arrow{r}{H^0(X,\Psi_{X,D+\Delta,n})} & \Hom_{\cO_X}(\cO_X(D+\Delta),\cO_X(D+\Delta)).
 \end{tikzcd}
 \]
 Here, the vertical equalities can be checked on the regular locus by Remark \ref{remark:hom-reflexive}, in which case the left one follows from (\ref{e-C-tensor}) and the right one is obvious.
 
 Thus $H^0(X,\Psi_{X,\Delta,n})$ is surjective if and only if $H^0(X,\Psi_{X,D+\Delta,n})$ is surjective. In view of Lemma \ref{lem:dual-criterion-for-log-quasi-F-splitting}, this concludes the proof.}
\end{proof}

\begin{remark}\label{r-big-open-QFS}
In the situation of Setting  \ref{setting:most-general-foundations-of-log-quasi-F-splitting}, 
let $\Delta$ be an effective $\Q$-divisor on $X$ with ${\cora\rdown{\Delta} }=0$, and let $X'$ be an open subset of $X$ such that ${\rm codim}_X(X \setminus X') \geq 2$. 
By a similar argument to the proof of Lemma \ref{l-inj-QFS1}, 
we see that $(X, \Delta)$ is $n$-quasi-$F$-split {if and only if} $(X', \Delta|_{X'})$ is $n$-quasi-$F$-split.

\end{remark}
\vspace{0.2em}

{Finally, we obtain a  cohomological criterion for whether a pair is quasi-$F$-split  analogous to Lemma \ref{lem:intro-cohomological-quasi-F-spliteness}
but  which we state in a more general framework.}

\begin{lemma} \label{lem:cohomological-criterion-for-log-quasi-F-splitting}
In the situation of Setting  \ref{setting:foundations-of-log-quasi-F-splitting}, {assume that $R$ is a local ring with maximal ideal $\m$ and}
let $\Delta$ be a $\Q$-divisor with $\rdown{\Delta}=0$. 
{Then $(X, \Delta)$ is $n$-quasi-$F$-split if and only if the following map is injective}
\[
H^d_{\m}(\Phi_{X, K_X+\Delta, n}) \colon H^d_{\m}(X, \cO_X(K_X+\Delta)) \to H^d_{\m}(X, Q_{X,K_X+\Delta,n}).
\]
\end{lemma}
\noindent Here $\cO_X(K_X+\Delta)=\cO_X(K_X)$ as $\rdown{\Delta}=0$. 
For the definition of $H^d_{\m}(X, -)$, see Notation \ref{notation:global-local-cohomology}.

\begin{proof}
Consider the following exact sequence
\begin{equation}\label{e1:cohomological-criterion-for-log-quasi-F-splitting}
0 \to K \to H^d_{\m}(X, \cO_X(K_X+\Delta)) \xrightarrow{
H^d_{\m}(\Phi_{X, K_X+\Delta, n})} H^d_{\m}(X, Q_{X,K_X+\Delta,n}), 
\end{equation}
where $K$ is the kernel of $H^d_{\m}(\Phi_{X, K_X+\Delta, n})$. 
By applying Matlis duality $\Hom_R(-, E)$ (which is exact), 
Lemma \ref{lem:Matlis-duality-for-highest-cohomology} yields an exact sequence 
\begin{equation}\label{e2:cohomological-criterion-for-log-quasi-F-splitting}
0 \leftarrow K^{\vee} \leftarrow H^0(X,\cO_X)^{\wedge} \xleftarrow{H^0(X, \Psi_{X, \Delta, n})^{\wedge}} H^0(X, Q^{\vee}_{X,\Delta,n})^{\wedge}.
\end{equation}
{Then the following is true:}
\begin{align*}
(X, \Delta)\text{ is }n\text{-quasi-}F\text{-split} 
&\Longleftrightarrow  H^0(X, \Psi_{X, \Delta, n})\text{ is surjective}\\
&\Longleftrightarrow  H^0(X, \Psi_{X, \Delta, n})^{\wedge}\text{ is surjective}\\
&\Longleftrightarrow  H^d_{\m}(\Phi_{X, K_X+\Delta, n})\text{ is injective}.
\end{align*}
Here the first equivalence follows from Lemma \ref{lem:dual-criterion-for-log-quasi-F-splitting}, 
the second one holds by \cite[\href{https://stacks.math.columbia.edu/tag/00MA}{Tag 00MA}]{stacks-project} and the fact that the completion $\widehat{R}$ 
of $R$ is a faithfully flat $R$-module \cite[\href{https://stacks.math.columbia.edu/tag/00MC}{Tag 00MC}]{stacks-project}, and the last one follows from 
\eqref{e1:cohomological-criterion-for-log-quasi-F-splitting} and 
\eqref{e2:cohomological-criterion-for-log-quasi-F-splitting}. {In the last assertion, we used that $K=0$ if and only if $K^{\vee}=0$, which is true because $K^{\vee\vee} \cong K$}\footnote{As local cohomology is independent of taking completion, this can be checked after replacing $R$ by its $\m$-adic completion $\hat{R}$, in which case we can apply \cite[\href{https://stacks.math.columbia.edu/tag/08Z9}{Tag 08Z9}]{stacks-project}.}. 
\end{proof}

\subsection{Kodaira dimension}
In this subsection, we show that quasi-$F$-split pairs $(X, \Delta)$ have non-{positive} Kodaira dimension 
(Proposition \ref{prop:non-negative-Kodaira-dimension-for-quasi-F-split-pairs}). 
This result was previously unknown even when $\Delta =0$.

\begin{proposition}
\label{prop:non-negative-Kodaira-dimension-for-quasi-F-split-pairs}
Let $X$ be an integral normal projective scheme over a perfect field $k$ of  characteristic $p>0$. Let $\Delta$ be a $\bQ$-divisor such that $\lfloor \Delta \rfloor = 0$. Suppose that $(X,\Delta)$ is quasi-$F$-split and let $n$ be the minimal number for which $(X,\Delta)$ is $n$-quasi-$F$-split. Then  $\kappa(X,-(K_X+\frac{1}{p^n-1}\rdown{ p^n\Delta})) \geq 0$.
\end{proposition}
\noindent In particular, if $X$ is quasi-$F$-split, then $\kappa(X) \leq 0$.

\begin{proof}
Since $(X, \Delta)$ is $n$-quasi-$F$-split, there exists an $\MO_X$-module homomorphism $\alpha \colon F_*W_n\cO_X(p\Delta) \to \cO_X$ 
such that $\alpha(F_*1)=1$, that is,  $\alpha \circ F=R^{n-1}$:  
\[
R^{n-1} \colon {W_n\MO_X}  \xrightarrow{F} F_*W_n\cO_X(p\Delta) \xrightarrow{\alpha} \cO_X. 
\]
Now consider the following short exact sequence
\[
0 \to F^n_* \cO_X(p^n\Delta) \xrightarrow{V^{n-1}} 
F_*W_n\cO_X(p\Delta) \xrightarrow{R} F_*W_{n-1}\cO_X(p\Delta) \to 0.
\]
Since $n$ was chosen as minimal, the composite $\MO_X$-module homomorphism 
\[
F^n_*\cO_X(p^n\Delta) \xrightarrow{V^{n-1}} F_*W_n\cO_X(p\Delta) \xrightarrow{\alpha} \cO_X
\]
is nonzero. Indeed, otherwise, 
there would exist an $\MO_X$-module homomorphism $F_*W_{n-1}\cO_X(p\Delta) \to \cO_X$ sending $F_*1$ to $1$, {thus} contradicting the 
minimality of $n$.

Observe that
\[
\Hom_{\MO_X} (F^n_*\cO_X(p^n\Delta), \cO_X) \cong H^0(X, \cO_X(-(p^n-1)K_X - \rdown{p^n\Delta})) 
\]
by Grothendieck duality. Hence, there exists an effective Weil divisor $D$ such that 
\[
D 
\sim - (p^n-1) \left(K_X + \frac{1}{p^n-1}\rdown{p^n\Delta}\right), 
\]
which implies {$\kappa(X, -(K_X+\frac{1}{p^n-1}\rdown{p^n\Delta})) \geq 0$.}
\end{proof}

It is an open problem whether quasi-$F$-split varieties $X$ have singularities of log canonical type, or, more generally, whether there exists $\Delta \geq 0$ such that $(X,\Delta)$ is log canonical and $K_X+\Delta \sim_{\bQ} 0$.

\subsection{Kawamata--Viehweg vanishing}
In this subsection, we generalise \cite[Theorem 1.7]{NY21} 
to show the Kawamata--Viehweg vanishing for quasi-$F$-split log pairs.

\begin{theorem} \label{thm:KVV-for-quasi-F-split-pairs}
 {Let $X$ be an integral normal scheme projective over an $F$-finite field $k$ of characteristic $p>0$.} Let $\Delta$ be a $p$-compatible $\bQ$-divisor such that $K_X+\Delta$ is  $\bQ$-Cartier and $\lfloor \Delta \rfloor =0$. We assume that $X$ is divisorially Cohen-Macaulay 
 and  $(X,\Delta)$ is quasi-$F$-split and naively keenly $F$-pure. 
Let $L$ be a {$\Q$-Cartier Weil} divisor such that $L-(K_X+\Delta)$ is ample. Then
\[
H^j(X,\cO_X(L)) = 0
\]
for all $j>0$.
\end{theorem}

\begin{proof}
Set $d := \dim X$. 
By assumptions, $A := L-(K_X+\Delta)$ is an ample $\bQ$-Cartier $\Q$-divisor, $L = K_X + \lceil A \rceil$, and $\{-A\} = \Delta$ (in particular, $\rup{A} =  A + \Delta$).

By Serre duality, 
it is enough to show that
\[
H^j(X,\cO_X(-\rup{A}))=0
\]
for $j < d$. 
Fix $j \in \Z$ with $j <d$. 

Let $n\in \Z_{>0}$ be {\cora a positive integer} such that $(X,\Delta)$ is $n$-quasi-$F$-split. Pick a positive integer $m \gg n$. Since $\cO_X(-\rup{A}) \to Q_{X,\Delta-\rup{A},m}$ 
is a split injection 
(Lemma \ref{l-inj-QFS1}), it is enough to show that $H^j(X, Q_{X,-A,m}) = 0$. Here $-A = \Delta - \rup{A}$.

Set $r:=m-n$. 
By descending induction on $i$, we will prove that
\begin{equation} \label{eq:KVV-induction}
H^j(X, Q_{X, -p^iA, m-i}) = 0
\end{equation}
for any $0 \leq i \leq r$. 

For the base case of the induction, we first check that (\ref{eq:KVV-induction}) holds 
when $i=r =m-n$. 
Since $-A$ is $p$-compatible and $(X, \{ - A \})$ is naively keenly $F$-pure, 
it follows from Proposition \ref{prop:C-cohen-macaulay} that 
$Q_{X,-p^rA, n}$ is Cohen-Macaulay. 
We have 
\[
H^j(X, Q_{X, -p^rA, n})^{\vee} \cong 
{\rm Ext}^{d-j}(Q_{X, -p^rA, n}, \omega_X) 
\cong 
H^{d-j}(X, \cHom(Q_{X, -p^iA, n}, \omega_X))
\]
where $H^{\vee}$ denotes the dual vector space of $H$, 
the first isomorphism holds by Serre duality \cite[Ch. III, Theorem 7.6]{hartshorne77}, and the second one follows from the fact that 
$Q_{X, -p^rA, n}$ is Cohen-Macaulay \cite[Theorem 3.3.10]{BH93}. 
Fix a positive integer $l$ such that $lA$ is Cartier, and 
{define $b := p^r \mod l$; this is the unique integer satisfying} $0 \leq b \leq l-1$ and {$b \equiv p^r \! \pmod l$}. 
Then
\[
Q_{X,-p^{r}A, n} =
Q_{X,{\{p^r\Delta\}} - \rup{p^rA},n} = Q_{X,{ \{p^{r}\Delta\}}-\rup{bA},n} \otimes_{\MO_X} \cO_X(-(p^r-b)A),
\]
where the first equality holds by ${\{ p^r\Delta\}} -\rup{p^rA} = \{ - p^r A\} + \rdown{-p^rA} =-p^rA$ 
and 
the second one follows from (\ref{e-C-tensor}). 
Since $r = m-n\gg 0$ and there are only finitely many possibilities for $b$, $\{p^r\Delta\}$, {\cora and therefore for} $Q_{X,\{p^{r}\Delta\} -\rup{bA},n}$, we get that (\ref{eq:KVV-induction})  holds  for $i=r$ by 
the Serre vanishing theorem.




From now on, we assume that $0  \leq i \leq r-1$. 
To carry out the induction, we can suppose that
\begin{equation}\label{e2:KVV-for-quasi-F-split-pairs}
H^j(X, Q_{X, -p^{i+1}A, m-i-1}) = 0 
\end{equation}
and aim for showing (\ref{eq:KVV-induction}).

Our key exact sequence (\ref{eq:C-restriction-sequence}) specialises to:
\[
0 \to F_*B_{X, -p^{i+1}A, m-i-1} 
\to Q_{X,-p^iA,m-i} \to F_* \cO_X(-p^{i+1}A) \to 0.
\]
The left and the right terms are direct summands of $F_*Q_{X,-p^{i+1}A,m-i-1}$ 
by Lemma \ref{lem:p-compatible-quasi-splitting} and the exact sequence 
\[
0 \to \cO_X(-p^{i+1}A) \xrightarrow{\Phi_{X, -p^{i+1}A, {m-i-1}}} 
Q_{X,-p^{i+1}A,{m-i-1}} \to B_{X,-p^{i+1}A,{m-i-1}} 
\to 0. 
\] 
Therefore, their $j$-th cohomologies $H^j(X, -)$ are zero by the induction hypothesis  (\ref{e2:KVV-for-quasi-F-split-pairs}), concluding the proof of (\ref{eq:KVV-induction}). 
\end{proof}

In the above proof, we used the following result which is a variant of Lemma \ref{lem:p-compatible-splitting}.
\begin{lemma} \label{lem:p-compatible-quasi-splitting}
In the situation of  
Setting \ref{setting:most-general-foundations-of-log-quasi-F-splitting}, 
let $\Delta$ be a $p$-compatible $\bQ$-divisor. Suppose that $(X,\{\Delta\})$ is $n$-quasi-$F$-split.  Then
\[
\Phi_{X, p^m\Delta, n} \colon 
\cO_X(p^m\Delta) \to Q_{X,p^m\Delta,n}
\]
splits for every $m \geq 0$.
\end{lemma}

\begin{proof}
Fix  $m \geq 0$. 
Since $\Delta$ is $p$-compatible, we have $\{ p^m \Delta\} \leq \{ \Delta \} $. 
Hence, $(X,\{ p^m\Delta\})$ is $n$-quasi-$F$-split, which implies that 
\[
\Phi_{X, \{ p^m\Delta\}, n} \colon \MO_X \to Q_{X, \{p^m\Delta\},n}
\]
splits. Then $\Phi_{X, p^m\Delta, n} =\Phi_{X, \rdown{p^m \Delta} +\{ p^m\Delta\}, n}$ splits 
by Lemma \ref{l-inj-QFS1}. 
\end{proof}

\subsection{Finite covers}

In this subsection, we study how quasi-$F$-{splittings} behave under finite covers. As an application, we show that log Calabi-Yau pairs $(\bP^1, \Delta)$ with standard coefficients are quasi-$F$-split for $p>3$ (Proposition \ref{prop:p1-log-calabi-yau}). This result will be important in the proof of Theorem \ref{thm:intro-3dim-klt}.

\begin{lemma} \label{lemma:splitting-finite-covers} Let $f \colon Y \to X$ be a finite surjective morphism of 
integral normal Noetherian $\F_p$-schemes. 
Let $D$ be a $\bQ$-divisor on $X$. Suppose that the degree of $f$ is not divisible by $p$. Then the $\MO_X$-module homomorphism 
\[
{f^*} \colon \cO_X(D) \to f_*\cO_Y(f^*D)
\]
induced by pulling-back sections splits.
\end{lemma}

\begin{proof}
Let $K$ and $L$ be the function fields of $X$ and $Y$,  respectively. By assumptions the extension $K \subseteq L$ is separable. There exists a field trace $\Tr_{L/K} \colon L \to K$ sending $\alpha \in L$ to the trace of the multiplication-by-$\alpha$ map on the vector space $L$ over $K$. Explicitly, if we take the minimal polynomial of $\alpha$ over $K$ and choose its roots $\sigma_1(\alpha), \ldots, \sigma_r(\alpha)$ in the algebraic closure $\overline{L}$ of $L$, then
\[
\Tr_{L/K}(\alpha) = \deg_{L/K(\alpha)}\sum_{i=1}^r \sigma_i(\alpha),
\]
where $K(\alpha)$ is the extension of $K$ generated by $\alpha$.

We claim that $\Tr_{L/K}$ extends to a morphism $\Tr_{Y/X} \colon f_*\cO_Y(f^*D) \to \cO_X(D)$. To this end we need to verify that given a local section $s\in f_*\cO_Y(f^*D)$, we have that $\sum_{i=1}^r \sigma_i(s)$ is a local section of $\cO_X(D)$, where $\sigma_i(s) \in \overline{L}$ are the roots of the minimal polynomial of $s \in L$. This can be checked after replacing $Y$ by a finite cover, and so we may assume that $Y \to X$ is Galois. In this case $\sigma_i(s)$ are given by an action of the Galois group, and the statement is straightforward (see \cite[Step 2 of Lemma 4.4]{tanaka22} for details).

Now the lemma follows as $\frac{1}{\deg_{L/K}}\Tr_{Y/X}$ splits $\cO_X(D) \to f_*\cO_Y(f^*D)$.
\end{proof}

{\cora 
\begin{remark}
The above proof shows that, under the same assumption as in 
Lemma \ref{lemma:splitting-finite-covers}, 
the trace map 
$\Tr_{L/K} \colon L \to K$ induces an $\MO_X$-module homomorphism 
\[
\Tr_{Y/X} \colon f_*\cO_Y(f^*D) \to \cO_X(D).
\]
\end{remark}
}

\noindent We warn the reader that we could not  replace $Y$ by a finite cover so that $Y \to X$ is Galois at the \emph{beginning} of the proof because then $\deg_{L/K}$ could have become divisible by $p$. 

{\cora We now prove finite-cover descents of quasi-$F$-splitting, a variant of which for $F$-splitting can be found in \cite[Example 1.1.10 (2)]{BrionKumar}.}

\begin{proposition}
\label{prop:descent}
{
Let $f \colon Y \to X$ be a morphism of integral normal  Noetherian $F$-finite schemes over a field $k$ of  characteristic $p>0$.
Let $D$ and $D_Y$ be $\mathbb{Q}$-divisors on $X$ and $Y$, respectively.
Assume that one of the following holds.
\begin{enumerate}
    \item $f$ is finite and surjective, and the degree of $f$ is not divisible by $p$.     
    \item $f \colon Y = X \times_k k' \to X$ is the morphism induced by a base change {to} an algebraic separable field extension {$k'$}  of $k$.  
\end{enumerate}
Further suppose that $f^*D \sim D_Y$ and $(Y,\{D_Y\})$ is $n$-quasi-$F$-split. Then $(X,\{D\})$ {is $n$-quasi-$F$-split}.}
\end{proposition}

\begin{proof}
As $ \{ f^*D\} =\{ D_Y \}$, we may assume that $f^*D=D_Y$.
In both cases, by construction (cf.\  (\ref{diagram:quasi-F-split-definition})), we have the following commutative diagram
\begin{center}
\begin{tikzcd}
f_*\cO_Y(D_Y) \arrow{r}{f_*\Phi_{Y, D_Y, n}} &[9mm] f_*Q_{Y,D_Y,n} \\
\cO_X(D) \arrow{r}{\Phi_{X, D, n}} \arrow{u}{f^*} & Q_{X,D,n} \arrow{u}.
\end{tikzcd}
\end{center}
Since $(Y,\{D_Y\})$ is $n$-quasi-$F$-split, {by Lemma \ref{l-inj-QFS1}} there exists an $\cO_Y$-module homomorphism
\[
\alpha \colon Q_{Y,D_Y,n} \to \cO_Y(D_Y)
\]
such that ${\alpha}({F_*1})=1$, 
{where ${F_*1} \in Q_{Y, 0, n}(V) = Q_{Y,D_Y,n}(V)$ and 
$1 \in \MO_Y(V) = \MO_Y(D_Y)(V)$ for $V := Y \setminus \Supp D_Y$}.

Next, we show that the $\MO_X$-module homomorphism  
$f^* \colon \cO_X(D) \to f_*\cO_Y(D_Y)$ splits. 
Indeed, in  Case (1), {this} follows from Lemma \ref{lemma:splitting-finite-covers}.
In  Case (2), as $X$ is excellent, {the existence of a splitting may be checked on the regular locus, and so we may assume that $X$ is regular (see Remark \ref{remark:hom-reflexive}).}
{Then $\rdown{D}$ is Cartier, where} $D = \rdown{D} +\{D\}$. 
Since $f$ is a limit of \'etale covers, 
$f^*\rdown{D} = \rdown{f^*D} = \rdown{D_Y}$, 
{and} since field extensions split, $\MO_X \to f_*\MO_Y$ splits. 
By taking the tensor product with the invertible sheaf $\MO_X(\rdown{D})$, {the homomorphism}
$f^* \colon \MO_X(D) =  \MO_X(\rdown{D}) \to f_*(\MO_Y(\rdown{D_Y})) =f_*\MO_X(D_Y)$ splits as well.


{By the above paragraph,} there exists an $\cO_X$-module homomorphism
\[
\beta \colon f_*\cO_Y(D_Y) \to \cO_X(D)
\]
such that $\beta(1)=1$. Therefore, the composite homomorphism
\[
Q_{X,D,n} \xrightarrow{} f_*Q_{Y,D_Y,n} \xrightarrow{f_*\alpha} f_*\cO_Y(D_Y) \xrightarrow{\beta} \cO_X(D)
\]
gives a splitting of $\Phi_{X, D, n} \colon \cO_X(D) \to Q_{X,D,n}$. Indeed, any $\MO_X$-module homomorphism $\cO_X(D) \to \cO_X(D)$ sending $1 \in \Gamma(X  \setminus D, \cO_X(D))$ to $1$ is the identity. 
By Lemma \ref{l-inj-QFS1}, $(X,\{D\})$ is $n$-quasi-$F$-split.
\end{proof}

\begin{corollary} \label{prop:quasi-F-splittings-under-finite-maps}  Let $f \colon Y \to X$ be a finite surjective morphism of normal varieties 
over an $F$-finite field $k$ of characteristic $p>0$. Let $\Delta$ and $\Delta_Y$ be effective $\bQ$-divisors on $X$ and $Y$, respectively, such that $\rdown{\Delta}=\rdown{\Delta_Y}=0$ and the following linear equivalence $f^*(K_X+\Delta) \sim K_Y + \Delta_Y$ holds. Suppose that the degree of $f$ is not divisible by $p$.
In this case, if $(Y,\Delta_Y)$ is $n$-quasi-$F$-split, 
then so is $(X,\Delta)$. 
\end{corollary}

\begin{proof}
Apply Proposition \ref{prop:descent}(1) for 
$D := K_X+\Delta$ and $D_Y := K_Y + \Delta_Y$. 
\end{proof}

\begin{corollary}\label{c-descent2}
Let $X$ be a  normal variety over an $F$-finite field $k$ of characteristic $p>0$. 
Let $\Delta$ be an effective $\Q$-divisor {satisfying} $\rdown{\Delta} =0$, {and} let
 $k'$ be an algebraic separable extension of $k$ such that $X \times_k k'$ is integral. 
 
 {If} the base change $(X \times_k k', \Delta \times_k k')$ is $n$-quasi-$F$-split, 
then so is $(X, \Delta)$. 
\end{corollary}

\begin{proof}
Apply Proposition \ref{prop:descent}(2) for $Y := X \times_k k'$, $D := \Delta,$ and $D_Y := \Delta \times_k k'$. 
\end{proof}

\begin{proposition} \label{prop:p1-log-calabi-yau} Let $k$ be a perfect field of characteristic $p>3$ and 
let $\Delta$ be an effective $\bQ$-divisor on $\bP^1_k$ with standard coefficients such that $K_{\bP^1} + \Delta \sim_{\bQ} 0$ and $\rdown{\Delta}=0$. Then $(\bP^1,\Delta)$ is $2$-quasi-$F$-split.
\end{proposition}
\begin{proof}
By Corollary \ref{c-descent2}, we can replace $k$ by its algebraic closure, and so we may assume that $k$ is algebraically closed.
A simple calculation (cf.\ Remark \ref{remark:ACC-P^1}) shows that one of the following holds
\begin{enumerate}
    \item $\Delta = \frac{1}{2}P_1 + \frac{1}{2}P_2 + \frac{1}{2}P_3 + \frac{1}{2}P_4$,
    \item $\Delta = \frac{2}{3}P_1 + \frac{2}{3}P_2 + \frac{2}{3}P_3$,
    \item $\Delta = \frac{1}{2}P_1 + \frac{3}{4}P_2 + \frac{3}{4}P_3$, 
    \item $\Delta = \frac{1}{2}P_1 + \frac{2}{3}P_2 + \frac{5}{6}P_3$
\end{enumerate}
for four distinct points $P_1, P_2, P_3, P_4$. By acting via an automorphism of $\bP^1$ we may assume that $P_1 = (0:1)$, $P_2 = (1:1)$, $P_3 = (1:0)$, and $P_4 = (a:1)$ for $a \neq 0,1$. For the ease of notation, we shall denote $P_1$ by $0$ and $P_3$ by $\infty$.

We claim that we can find a finite map $f \colon E \to \bP^1$ from a smooth projective curve $E$ such that $f^*(K_X+\Delta) = K_E$. To this end, we consider the totalisation $\mathrm{Tot}(\cO_{\bP^1}(m))$ of the line bundle $\cO_{\bP^1}(m)$ for $m\in \bZ_{>0}$. Let 
\begin{itemize}
    \item $(x:1)$ and $(1:y)$ be the coordinates of the base curve $\bP^1$ {on} $\bP^1\,\backslash\, \{\infty\}$ and $\bP^1\,\backslash\, \{0\}$, respectively, so that $xy=1$,
    \item $u$ and $v$ be the natural coordinates of $\cO_{\bP^1}(m)$ over $\bP^1\, \backslash\, \{\infty\}$ and $\bP^1\,\backslash\, \{0\}$, respectively, so that $u = x^mv$.
\end{itemize}
{We will show that} the sought-after finite cover is given by normalising the curves described by the following equations.
\begin{alignat*}{6}
  &(1) \quad  &&u^2 &&= x(x-1)(x-a) \quad &&\text{ and } \quad &&v^2 = y(y-1)(ay-1), \qquad &&\text{ for } m=2,\\
  &(2)   &&u^3 &&= x(x-1) &&\text{ and } &&v^3 = -y(y-1), &&\text{ for } m=1,\\
  &(3)  &&u^4 &&= x^2(x-1) &&\text{ and } &&v^4 = -y(y-1), &&\text{ for } m=1,\\
  &(4)    &&u^6 &&= x^3(x-1)^2 &&\text{ and } &&v^6 = {y(y-1)^2}, &&\text{ for } m=1.
\end{alignat*}
Here (1), (2), (3), (4) correspond to the cases of $\Delta$ laid out in the first paragraph of the proof. In each line, the first equation describes the curve in $\mathrm{Tot}(\cO_{\bP^1}(m))$ over $\bP^1 \, \backslash \, \{\infty\}$ and the second equation over $\bP^1 \, \backslash \, \{0\}$. Since $u = x^mv$ and $y = \frac{1}{x}$, one immediately sees that  in each line the two equations agree over $\bP^1 \, \backslash \, \{0,\infty\}$.

The curve described in (1) is smooth by the Jacobian criterion and, as shown by the equation, ramifies of index $2$ at each of the points $P_1, P_2, P_3, P_4$. The curve described in (2) is smooth by the Jacobian criterion as well, and ramifies of index $3$ over $P_1, P_2, P_3$.  The normalisation of the curve described in (3) ramifies over $P_1, P_2, P_3$ with indices $2, 4, 4$, respectively.
\footnote{The curve described in (3) is smooth except over $P_1$, whereat up to localisation and finite extension, we may assume that $x-1$ is a unit admitting a fourth root, and so the equation of the curve is equivalent to $u^4=x^2$; its normalisation consists of two branches: $x=u^2$ and $x=-u^2$ each ramifying with index $2$ over $P_1$.}. The normalisation of the curve described in (4) ramifies over $P_1, P_2, P_3$ with indices $2, 3, 6$, respectively
\footnote{The curve described in (3) is smooth except over $P_1$ and $P_2$, whereat up to localisation and finite extension, we may assume that the equations are equivalent to $u^6=x^3$ and $u^6 = (x-1)^2$, respectively; over $P_1$ the normalisation consists of three branches ramifying with index two: $x=u^2$, $x=\zeta u^2$, and $x= \zeta^2u^2$ where $\zeta$ is the third root of unity, while over $P_2$ the normalisation consists of two branches ramifying with index three: $x-1 = u^3$ and $x-1 = -u^3$.}.  

 By the above calculation of ramification indices and the Riemann-Hurwitz formula, we immediately get that $f^*(K_X+\Delta) = K_E$. In particular, $K_E$ is trivial, and {so}  $E$ is $2$-quasi-$F$-split {by Remark \ref{remark:elliptic-curves-quasi-F-split}}.
Therefore, $(X,\Delta)$ is $2$-quasi-$F$-split by Corollary \ref{prop:quasi-F-splittings-under-finite-maps}. 
\end{proof}

\subsection{Definition of pure quasi-$F$-splitting}
\label{ss:pure-quasi-F-split}

\noindent 
{Throughout this subsection, we work under the assumptions of Setting \ref{setting:most-general-foundations-of-log-pure-quasi-F-splitting}. 

\begin{setting} \label{setting:most-general-foundations-of-log-pure-quasi-F-splitting}
Let $X$ be a $d$-dimensional 
integral normal $F$-finite
$\F_p$-scheme admitting a dualising complex $\omega^{\mydot}_X$. We fix a 
Weil divisor $S$ on $X$ which {either} is a prime divisor or is equal to zero. 
\end{setting}
\begin{setting} \label{setting:foundations-of-log-pure-quasi-F-splitting}
Let $R$ be 
an $F$-finite Noetherian domain of characteristic $p>0$. Let $X$ be a $d$-dimensional integral normal scheme which is projective over $\Spec R$. We fix a 
Weil divisor $S$ on $X$ which is either a prime divisor or is equal to zero.
\end{setting}}
We denote the ideal sheaf of $S$ by $\cI_S$. We remind the reader that $W_n\cI_S$ and $W_n\cO_X(-S)$ are \emph{not} equal in general.

Let $\Delta$ be a (non-necessarily effective) $\bQ$-divisor on $X$. 
We define $Q^S_{X,\Delta,n}$ and $\Phi^S_{X, \Delta, n}$ by the following pushout diagram of $W_n\MO_X$-modules: 
\begin{equation}\label{diagram:purely-quasi-F-split-definition}
\begin{tikzcd}
W_n\cI_S(\Delta) \arrow{r}{F} \arrow{d}{R^{n-1}} & F_* W_n \cI_S(p\Delta)  \arrow{d}\\
\cI_S(\Delta) \arrow{r}{\Phi^S_{X, \Delta, n}}& \arrow[lu, phantom, "\usebox\pushoutdr" , very near start, yshift=0em, xshift=0.6em, color=black] Q^S_{X, \Delta, n}.
\end{tikzcd}
\end{equation}
We remind the reader that
$\cI_S(\Delta) =  \MO_X(\Delta-S)$.

\begin{remark}
In our applications, we are mainly interested in the case {when $\Delta=S+B$ for a prime divisor $S$ and a $\bQ$-divisor $B$ such that $S \not \subseteq \Supp B$ and  $\rdown{B}=0$.}
Under these assumptions, $\cI_S(\Delta) = \MO_X$.
\end{remark}

We define a coherent $W_n\cO_X$-module
\begin{align*}
B^S_{X, \Delta, n} &:= {\rm Coker}(W_n\cI_S(\Delta) \xrightarrow{F} F_*W_n\cI_S(p\Delta)) \\
&\hphantom{:}= F_*W_n\cI_S({ p}\Delta)/F(W_n\cI_S({\Delta})), 
\end{align*} 
{and set $B^S_{X,\Delta} := B^S_{X,\Delta,1}$}.

\begin{remark}
{The key properties of the construction of $Q^S_{X,\Delta,n}$,  $B^S_{X,\Delta,n}$, and $\Phi^S_{X,\Delta,n}$ may be encapsulated by the following diagram}
\begin{equation}\label{e-big-BC-diagram2}
\begin{tikzcd}
& 0 \arrow{d} & 0 \arrow{d} & & \\
& F_*W_{n-1}\cI_S(p\Delta) \arrow{d}{V} \arrow[r,dash,shift left=.1em] \arrow[r,dash,shift right=.1em] &  F_*W_{n-1}\cI_S(p\Delta) \arrow{d}{ FV} &&\\
0 \arrow{r} &  W_n\cI_S(\Delta) \arrow{d}{R^{n-1}} \arrow{r}{F} & F_*W_n\cI_S(p\Delta) \arrow{r} \arrow{d} & B^S_{X, \Delta, n} \arrow[d,dash,shift left=.1em] \arrow[d,dash,shift right=.1em] \arrow{r} & 0 \\
0 \arrow{r} &  \cI_S(\Delta) \arrow{d} \arrow{r}{\Phi^S_{X, \Delta, n}} & Q^S_{X, \Delta, n} \arrow[lu, phantom, "\usebox\pushoutdr" , very near start, yshift=0em, xshift=1em, color=black] \arrow{d} \arrow{r} & B^S_{X, \Delta, n} \arrow{r} & 0.\\
& 0 & 0 &&
\end{tikzcd}
\end{equation}
{All the horizontal and vertical sequences are exact, as} $F \colon W_n\cI_S(\Delta) \to F_*W_n\cI_S(p\Delta)$ is injective, $R^{n-1} \colon W_n\cI_S(\Delta) \to \cI_S(\Delta)$ is surjective, and 
(\ref{diagram:purely-quasi-F-split-definition}) is a pushout diagram. 

\end{remark}

\begin{definition}
\label{definition:pure-quasi-F-split}
In the situation of Setting \ref{setting:most-general-foundations-of-log-pure-quasi-F-splitting}, 
let $B$ be a $\Q$-divisor on $X$ such that $S \not\subseteq \Supp B$ and $\rdown{B}=0$. We say that $(X, S+B)$ is \emph{purely $n$-quasi-$F$-split (along $S$)} if  there exists 
a $W_n\MO_X$-module homomorphism $\alpha \colon F_*W_n\cI_S(p(S+B)) \to \cI_S(S+B) = \MO_X$ such that  ${R^{n-1}} = \alpha \circ F$: 
\begin{equation*} 
\begin{tikzcd}
W_n\cI_S(S+B) \arrow{r}{F} \arrow{d}{{ R^{n-1}}} & F_* W_n\cI_S(p(S+B)) \arrow[dashed]{ld}{\alpha} \\
\cO_X.
\end{tikzcd}
\end{equation*}

We call $(X, S+B)$ \emph{purely quasi-$F$-split (along $S$)} if it is purely $n$-quasi-$F$-split for some $n \in \Z_{>0}$. 
{ We say that $(X,S+B)$ \emph{purely $n$-quasi-$F$-pure (purely quasi-$F$-pure, resp.)}  if it is purely $n$-quasi-$F$-split (purely quasi-$F$-split, resp.) affine-locally.}
\end{definition}

By exactly the same arguments as in Subsection \ref{ss-log-QFS}, we get the following results. 

\begin{proposition} \label{prop:adjoint-properties-of-C}
In the situation of Setting \ref{setting:most-general-foundations-of-log-pure-quasi-F-splitting},  let $\Delta$ be a $\bQ$-divisor on $X$. 
Then the following hold. 
\begin{enumerate}
    \item {The map $F_*W_n\cI_S(p\Delta) \to Q^S_{X,\Delta,n}$ induces an isomorphism of $W_n\MO_X$-modules }
    \[
    Q^S_{X,\Delta,n} \cong \frac{F_*W_n\cI_S({ p}\Delta)}{pF_*W_n\cI_S({p}\Delta)}.
    \] 
    \item 
    $Q^S_{X,\Delta,n}$ is naturally an $\cO_X$-module. In particular, 
    {it} is a coherent $\MO_X$-module and $\Phi^S_{X, \Delta, n} \colon \cI_S({\Delta}) \to Q^S_{X, \Delta, n}$ is an $\MO_X$-module homomorphism. 
    \item $B^S_{X, \Delta, n}$ is naturally an $\cO_X$-module. In particular, 
    {it} is a coherent $\MO_X$-module. 
\end{enumerate}
\end{proposition}

\begin{proposition} \label{prop:pure-definition-via-splitting}
In the situation of Setting \ref{setting:most-general-foundations-of-log-pure-quasi-F-splitting}, 
let $B$ be a $\Q$-divisor on $X$ such that $S \not\subseteq \Supp B$ and $\rdown{B}=0$. 
Then the following are equivalent. 
\begin{enumerate}
    \item $(X, S+B)$ is purely $n$-quasi-$F$-split. 
    \item $\Phi^S_{X, S+B, n} \colon \cO_X \to Q^S_{X, S+B,n}$ splits as a $W_n\MO_X$-module homomorphism. 
    \item $\Phi^S_{X, S+B, n} \colon \cO_X \to Q^S_{X,S+B,n}$ splits as an $\MO_X$-module homomorphism. 
\end{enumerate}
\end{proposition}

One can check that 
\begin{equation}\label{e-C-tensor2}
Q^S_{X,\Delta+D,n} \cong Q^S_{X,\Delta,n} \otimes_{\MO_X} \cO_X(D)
\end{equation} 
for any Cartier divisor $D$.

\begin{lemma}[{cf.\ (\ref{e-big-BC-diagram2}) and (\ref{eq:key-sequence-for-WnI})}]
In the situation of Setting \ref{setting:most-general-foundations-of-log-pure-quasi-F-splitting}, 
 let $\Delta$ be a $\bQ$-divisor on $X$. 
Take $n \in \Z_{>0}$. 
Then we have the following exact sequences of coherent $\MO_X$-modules: 
\begin{align}\label{eq:adjoint-C-quotient-sequence}
&0 \to \cI_S(\Delta) \xrightarrow{\Phi^S_{X, \Delta, n}} Q^S_{X,\Delta,n} \to B^S_{X,\Delta,n} \to 0 \\
\label{eq:adjoint-Bn-sequence}
&0 \to F_* \cB^S_{X,p\Delta,n} \to \cB^S_{X,\Delta,n+1} \to { \cB^S_{X,\Delta}} \to 0 \\
\label{eq:adjoint-C-restriction-sequence}
&0 \to 
F_*B^S_{X, p\Delta, n} \to Q^S_{X,\Delta,n+1} \to F_* \cI_S(p\Delta) \to 0. 
\end{align}
\end{lemma}

\begin{proposition}\label{prop:adjoint-C-cohen-macaulay} 
In the situation of Setting \ref{setting:most-general-foundations-of-log-pure-quasi-F-splitting},
let $B$ be a $\bQ$-divisor on $X$ 
such that $S \not\subseteq \Supp B$. 
Suppose that $X$ is divisorially Cohen-Macaulay and that $(X, S+ \{{p^j}B \})$ is purely naively keenly $F$-pure {for every $j\geq 0$}. 

Then $Q^S_{X, S+B,n}$ and $\cB^S_{X, S+B,n}$ are 
coherent Cohen-Macaulay $\MO_X$-modules for every $n \in \Z_{>0}$.  
\end{proposition}

\noindent In particular, under the same assumptions as in the above proposition, $Q^S_{X,\Delta,n}$ is a coherent reflexive $\MO_X$-module, and we get that for every Weil divisor $D$
\[
Q^S_{X,\Delta+D,n} \cong (Q^S_{X,\Delta,n} \otimes_{\MO_X} \cO_X(D))^{\vee\vee}.
\]

{For any $\Q$-divisor $\Delta$ on $X$,
we define $\Psi^S_{X, \Delta, n}$ by applying $\cHom_{\cO_X}(-, \cI_S(\Delta))$ to the $\cO_X$-module homomorphism $\Phi^S_{X, \Delta, n} \colon \cI_S(\Delta) \to Q^S_{X,\Delta,n}$: 
\begin{equation}\label{eq:dual-defin-pure-quasi-F-split}
\Psi^S_{X, \Delta, n} \colon  
\cHom_{\MO_X}(Q^S_{X,\Delta,n}, \cI_S(\Delta)) \to \cHom_{\MO_X}(\cI_S(\Delta), \cI_S(\Delta)) = \cO_X,
\end{equation}
where the last equality follows from Remark \ref{remark:hom-reflexive}. }

\begin{lemma} \label{lem:dual-criterion-for-pure-quasi-F-splitting}
In the situation of Setting \ref{setting:most-general-foundations-of-log-pure-quasi-F-splitting}, 
let $\Delta$ be an {arbitrary} $\Q$-divisor on $X$. 
{Then $\Phi^S_{X,\Delta,n} \colon \cI_S(\Delta) \to Q^S_{X,\Delta,n}$ splits if and only if the following map is surjective}
\[
H^0(X, \Psi^S_{X, \Delta, n}) \colon  \Hom_{\cO_X}(Q^S_{X,\Delta,n}, \cI_S(\Delta)) \to  H^0(X,\cO_X).
\]
\end{lemma}
 When $\Delta=S+B$ for a $\bQ$-divisor $B$ satisfying $\rdown{B}=0$ and $S \not \subseteq \Supp B$, we get a map $\Psi^S_{X, S+B, n} \colon (Q^S_{X,S+B,n})^{\vee} \to \cO_X$. Then the above lemma shows that $(X,S+B)$ is purely $n$-quasi-$F$-split if and only if $H^0(X, \Psi^S_{X, S+B, n}) \colon H^0(X,(Q^{S}_{X,S+B,n})^{\vee}) \to H^0(X,\cO_X)$ is surjective.

\begin{lemma}\label{l-inj-QFS2}
In the situation of Setting \ref{setting:most-general-foundations-of-log-pure-quasi-F-splitting}, 
let $\Delta$ be an arbitrary $\bQ$-divisor on $X$ and let $D$ be a Weil divisor on $X$.
Then the following are equivalent. 
\begin{enumerate}
    \item $\Phi^S_{X, \Delta, n} \colon \cI_S(\Delta) \to Q^S_{X, \Delta, n}$  splits as an $\MO_X$-module homomorphism.
    \item $\Phi^S_{X, D+\Delta, n} \colon \cI_S(D+\Delta) \to Q^S_{X, D+\Delta, n}$  splits as an $\MO_X$-module homomorphism.
\end{enumerate}
\end{lemma}

\begin{lemma} \label{lem:cohomological-criterion-for-pure-quasi-F-splitting}
In the situation of Setting  \ref{setting:foundations-of-log-pure-quasi-F-splitting}, let  $B$ be a $\bQ$-divisor on $X$ such that $S \not\subseteq \Supp B$ and $\rdown{B}=0$.  Assume that $R$ is a local ring with maximal ideal $\m$. 

Then $(X,S+B)$ is purely $n$-quasi-$F$-split if and only if
\[
H^d_{\m}(\Phi^S_{X, K_X+S+B, n}) \colon H^d_{\m}(X, \cO_X(K_X+B)) \to H^d_{\m}(X, Q^S_{X,K_X+S+B,n})
\]
is injective.
\end{lemma}

\subsection{Quasi-$F$-stable sections}
In this subsection, we define quasi-$F$-stable sections. The reader is advised to first look at Subsection \ref{ss:intro-to-lifting-quasi-stable-sections}, where a simple definition for line bundles in the non-boundary case is given.

\begin{definition} \label{definition:quasi-F-stable-sections} In the situation of Setting \ref{setting:most-general-foundations-of-log-quasi-F-splitting}, let $L$ and $\Delta$ be $\bQ$-divisors on $X$. The \emph{submodule (subspace) of {\cora $n$-}quasi-$F$-stable sections} $qS^0_n(X,\Delta;L)$ is defined by 
\[
qS^0_n(X,\Delta;L) := \image\big( 
H^0(X,Q^{\vee}_{X,\Delta-L,n}) \xrightarrow{\Phi_{X, \Delta -L, n}^{\vee}} H^0(X, \cO_X(\lceil L-\Delta\rceil))\big).
\]
\end{definition}
\noindent The above map $\Phi_{X, \Delta -L, n}^{\vee}$ 
is induced by applying $\mathcal Hom_{\MO_X}(-, \MO_X)$ and $H^0(X, -)$ to 
$\Phi_{X, \Delta -L, n} \colon \cO_X(\Delta-L) \to Q_{X,\Delta-L,n}$.\\

When $\rdown{\Delta}=0$ and $L$ is a Weil divisor, we have 
\[
qS^0_n(X,\Delta; L) \subseteq H^0(X, \cO_X(\lceil L-\Delta\rceil))  = H^0(X,\cO_X(L)).
\]
However, we do need to consider the most general case of the above definition as non-integral $L$ and non-effective $\Delta$ come up naturally in the context of pullbacks and restrictions of Weil divisors.


We are particularly interested in the case {when $\rdown{\Delta}=0$ and  $\Delta \geq \{L\}$}.
For example, this happens for $L$ that comes from restricting a Weil divisor to a plt centre with $\Delta$ being the different  
of a log pair (Lemma \ref{l:special-cases-of-L}). 

\begin{lemma} \label{l:special-cases-of-L}
Under the same notation as in Definition \ref{definition:quasi-F-stable-sections}, 
suppose that $(W,X+B_W)$ is a plt pair, $K_X+\Delta = (K_W+X+B_W)|_X$, and $L=L_W|_X$ for a Weil divisor $L_W$ on $W$. 
Then $\rdown{\Delta}=0$ and $\Delta \geq \{L\}$. 
In particular, $\lceil L - \Delta \rceil = \lfloor L \rfloor$ and 
\[
qS^0_{n}(X,\Delta;L) = qS^0_{n}(X,\Delta-\{L\};\rdown{L}) \subseteq H^0(X,\cO_X(\rdown{L})).
\]
\end{lemma}

\begin{proof}
Since $(W,X+B_W)$ is plt, it holds that $\rdown{\Delta}=0$. 
Let us show $\Delta \geq \{L\}$. 
Fix a prime divisor $\Gamma$ on $X$ with generic point $\xi$. 
Let $(-)_{\xi}$ denote the localisation at $\xi$. Then 
\[
(K_{W}+X)|_{X_{\xi}} = K_{X_{\xi}} + \left(1-\frac{1}{N}\right)\xi \qquad \textrm{and} \qquad\! 
L_W|_{X_{\xi}} = \frac{\ell}{N}\xi
\]
for some $\ell \in \bZ_{\geq 0}$ by Lemma \ref{lem:everything-about-plt-surface-pairs}(2)(3) applied to $(W_{\xi},X_{\xi})$, where $N$ is the determinant of the minimal resolution of $W_{\xi}$. 
Since $(K_W+X)|_{X_{\xi}} \leq (K_W+X+B_W)|_{X_{\xi}} = K_{X_{\xi}} + \Delta_{\xi}$, it follows from $\rdown{\Delta}=0$ that $\Delta_{\xi} = (1-\frac{1}{N}+\lambda)\xi$ for some $0 \leq \lambda < \frac{1}{N}$. Finally, the inequality $\Delta \geq \{ L\}$ holds, because 
\[
\mathrm{coeff}_{\Gamma}(\Delta) = 1-\frac{1}{N}+\lambda \geq \Big\{\frac{\ell}{N}\Big\} = \mathrm{coeff}_{\Gamma}(\{L\}), 
\]
where $\mathrm{coeff}_{\Gamma}(D)$ denotes the coefficient of $\Gamma$ in a $\Q$-divisor $D$. 
\end{proof}

\begin{proposition} \label{prop:qS^0-for-quasi-F-split}
In the situation of Setting \ref{setting:most-general-foundations-of-log-quasi-F-splitting}, let $L$ and $\Delta$ be $\bQ$-divisors on $X$. Suppose that $(X,\{\Delta-L\})$ is $n$-quasi-$F$-split. Then
\[
qS^0_n(X,\Delta;L) = H^0(X,\cO_X(\rup{L-\Delta})).
\]
\end{proposition}

\begin{proof}
By definition, $\Phi_{X, \{\Delta -L \}, n} \colon \cO_X \to Q_{X,\{\Delta-L\},n}$ 
splits, and so
\[
\Phi_{X, \Delta -L, n} :
\cO_X(\rdown{\Delta-L})  \to Q_{X,\Delta-L,n}
\]
splits as well (Lemma \ref{l-inj-QFS1}). 
In particular, its $\MO_X$-dual 
\[
Q^{\vee}_{X,\Delta-L,n} =\cHom_{\MO_X}(Q_{X,\Delta-L,n}, \MO_X) 
\to \cHom_{\MO_X}(\MO_X(\rdown{\Delta -L}), \MO_X) =  \cO_X(\rup{L-\Delta}) 
\]
is a split surjection, thus
$
H^0(X,Q^{\vee}_{X,\Delta-L,n}) \to H^0(X,\cO_X(\rup{L-\Delta}))
$
is also surjective. This concludes the proof of the lemma.
\end{proof}

Analogously, we make an adjoint definition. 

\begin{definition}\label{d-pureQFS-sections}
In the situation of Setting \ref{setting:most-general-foundations-of-log-pure-quasi-F-splitting}, let $L$ and $B$ be $\bQ$-divisors on $X$. The \emph{submodule (subspace) of purely {\cora $n$-}quasi-$F$-stable sections} 
$qS^0_{n,\adj}(X,S+B;L)$ is defined {as follows}
\[
qS^0_{n,\adj}(X,S+B;L) := \image\big(
H^0(X,(Q^S_{X,S+B-L,n})^{\vee}) \xrightarrow{(\Phi^S_{X, S+B-L, n})^{\vee}}
H^0(X, \cO_X(\lceil L-B\rceil))\big).
\]
\end{definition}

\noindent The above map $(\Phi_{X, { S+B} -L, n}^S)^{\vee}$ 
is induced by applying $\mathcal Hom_{\MO_X}(-, \MO_X)$ and $H^0(X, -)$ to 
$\Phi^S_{X,S+B -L, n} \colon \cI_S(S+B- L) \to Q^S_{X,S+B-L,n}$.

\begin{proposition} \label{prop:qS^0-for-pure-quasi-F-split}
In the situation of Setting \ref{setting:most-general-foundations-of-log-pure-quasi-F-splitting},  let $L$ and $B$ be $\bQ$-divisors on $X$. Suppose that $(X,S+\{B-L\})$ is purely $n$-quasi-$F$-split. Then
\[
qS^0_{n,\adj}(X,S+B;L) = H^0(X,\cO_X(\rup{L-B})).
\]
\end{proposition}
\begin{proof}
The proof is completely analogous to that of Proposition \ref{prop:qS^0-for-quasi-F-split}.
\end{proof}

 \begin{lemma}\label{l-qS^0-adj-nonadj}
In the situation of {Setting \ref{setting:most-general-foundations-of-log-pure-quasi-F-splitting}},  
 let $B$ be a $\Q$-divisor on $X$  such that $S \not\subseteq \Supp B$ and $\rdown{B}=0$. 
Let $L$ be a Weil divisor on $X$. 

Then the following inclusions hold {for every $0 \leq t < 1$}: 
\[
qS^0_{n, \adj}(X, S+B; L) \subseteq 
qS^0_n(X, {tS+} B; L) \subseteq H^0(X, \MO_X(L)).
\]
\end{lemma}

\begin{proof}
 We have the following commutative diagram: 
\[
\begin{tikzcd}[row sep=1cm,column sep=2cm]
W_n\MO_X({tS+}B-L) \arrow{r}{ F_1}\arrow[hookrightarrow]{d}{i} & F_*W_n\MO_X(p({tS+}B-L)) \arrow[hookrightarrow]{d}{j}\\
W_n\cI_S(S+B-L)\arrow{r}{F_2} \arrow{d}{R^{n-1}}& F_*W_n\cI_S({ p}(S+B-L))\\
\mathllap{\MO_X({tS+}B-L) =}\  \cI_S(S+B-L),
\end{tikzcd}
\]
where $i$ and $j$ are the natural inclusions, {and $F_1$ and $F_2$ denote the Frobenii}. 
Since $Q_{X, {tS+}B-L, n}$ and $Q^S_{X, S+B-L, n}$ are the pushouts of 
$(F_1, R^{n-1} \circ i)$  and $(F_2, R^{n-1})$ respectively, {we get a map $Q_{X, {tS+}B-L, n} \to Q^S_{X, S+B-L, n}$ sitting inside the following diagram}
\begin{center}
\begin{tikzcd}[column sep = huge]
{\MO_X(-L)} \arrow[bend right = 30]{rr}{\Phi^S_{X, S+B-L, n}} \arrow{r}{\Phi_{X, {tS+}B-L, n}} &  Q_{X, {tS+}B-L, n} \arrow{r} & Q^S_{X, S+B-L, n}.
\end{tikzcd}    
\end{center}

Now, by applying $\Hom_{\cO_X}(-, \cO_X)$, we get a factorisation of 
{$(\Phi^S_{X, S+B-L, n})^{\vee}$} 
\[
H^0(X, (Q^S_{X, S+B-L, n})^{\vee}) \to H^0(X, Q_{X, {tS+}B-L, n}^{\vee}) 
\xrightarrow{{\Phi_{X, {tS+}B-L, n}^{\vee}}} H^0(X, \MO_X(L)), 
\]
which immediately implies the statement of the lemma.
\end{proof}

{\begin{remark} \label{remark:stable-sections-changing-parameter}
In the situation of {Setting \ref{setting:most-general-foundations-of-log-quasi-F-splitting}}, let $L$, $\Delta$, and $\Delta'$ be $\bQ$-divisors on $X$ such that $\Delta \leq \Delta'$. Then an analogous argument to the one above shows that 
$qS^0_n(X, \Delta; L) \subseteq qS^0_n(X, \Delta'; L)$.
\end{remark}}

\subsection{Quasi-$F$-splittings under birational morphisms}
In this subsection, we will show that a log pullback of a quasi-$F$-split pair is quasi-$F$-split provided that some additional assumptions are satisfied. 
We start by proving two auxiliary lemmas.

\begin{lemma} \label{lemma:pushforward-of-Watanabe-F-split}
In the situation of Setting \ref{setting:most-general-foundations-of-log-quasi-F-splitting}, let
{$(X,\Delta)$ be a $p$-compatible log pair} such that $(X,\{\Delta\})$ 
is klt. 
Let $f \colon Y \to X$ be a proper birational morphism 
of integral normal divisorially Cohen-Macaulay 
schemes. 
Let $\Delta_Y$ be the log pullback of $\Delta$, that is, $K_Y + \Delta_Y = f^*(K_X+\Delta)$ and $f_*\Delta_Y=\Delta$. 

Suppose that 
$
R^jf_*\cO_Y(K_Y + \rup{-p^i(K_Y+\Delta_Y)}) = 0
$
for every $i\geq 0$ and $j>0$. Then 
\[
R^jf_*\cO_Y(p^i(K_Y+\Delta_Y)) = 
\begin{cases}
\cO_X(p^i(K_X+\Delta)) \qquad & (j=0)\\
0 \qquad & (j>0)
\end{cases}
\]
for every $i \geq 0$.
\end{lemma}
\noindent 
Note that $K_Y + \rup{-p^i(K_Y+\Delta_Y)} = K_Y + \{p^i\Delta_Y\} -p^i(K_Y+\Delta_Y)$, and so 
the assumption $R^jf_*\cO_Y(K_Y + \rup{-p^i(K_Y+\Delta_Y)}) = 0$ 
is a consequence of Kawamata--Viehweg vanishing whenever it is valid, provided that {we assume that $f$ is a log resolution of $(X,\Delta)$, or more generally that} $(Y,\{p^i\Delta_Y\})$ is klt for every $i \geq 0$.

\begin{proof}
Fix $i \geq 0$. 
Set $B := \{\Delta\}$ and let $B_Y$ be its log pullback, that is, $K_Y + B_Y = f^*(K_X+B)$ and $f_*B_Y = B$.

By $p^i(K_Y+\Delta_Y) = f^*( p^i(K_X+\Delta))$, 
we have that
\[
f_*\cO_Y(p^i(K_Y+\Delta_Y)) = \cO_X(p^i(K_X+\Delta)),
\]
and so we only need to show that $R^jf_*\cO_Y(p^i(K_Y+\Delta_Y))=0$ for $j>0$.

Set $D_Y := -p^i(K_Y+\Delta_Y)$. 
By Lemma \ref{lem:dual-of-vanishing}, it suffices to show that 
\begin{enumerate}
    \item $f_*\cO_Y(K_Y + \rup{D_Y})$ is Cohen--Macaulay and 
    \item $R^jf_*\cO_Y(K_Y + \rup{D_Y}) =0$ for every $j>0$. 
\end{enumerate}
Note that (2) holds by assumption{, and so it} is enough to prove (1).

Set $D := -p^i(K_X+\Delta)$ and recall that
\begin{align*}
K_Y + \rup{D_Y} &= K_Y + \{p^i\Delta_Y\} -p^i(K_Y+\Delta_Y) \quad{\rm and}\\
K_X + \rup{D} &= K_X + \{p^i\Delta\} -p^i(K_X+\Delta).
\end{align*}
{In particular, $K_X+\lceil D \rceil$ is $\bQ$-Cartier.} We claim that $f_*\cO_Y(K_Y+\lceil D_Y \rceil) = \cO_X(K_X + \lceil D \rceil)$. To this end, we calculate 
\begin{align*}
\Gamma &:= K_Y + \lceil D_Y \rceil - f^*(K_X+\lceil D \rceil) \\
&\hphantom{:}= K_Y + \{p^i\Delta_Y\} - f^*(K_X + \{p^i\Delta\}) \\
&\hphantom{:}= f^*B - B_Y + \{p^i\Delta_Y\} - f^*\{p^i\Delta\} \\
&\hphantom{:}\geq -B_Y. 
\end{align*}
Here the third equality follows from the identity $K_Y+B_Y = f^*(K_X+B)$ and the last inequality holds by the fact that $B = \{\Delta\} \geq \{p^i\Delta\}$ as $\Delta$ is $p$-compatible.

Since $(X,B)$ is klt, {we have} $\rup{-B_Y} \geq 0$, and so $\rup{\Gamma} \geq 0$. As $\Gamma$ is $f$-exceptional, we get that
$
f_*\cO_Y(K_Y+\lceil D_Y \rceil) = \cO_X(K_X + \lceil D \rceil)
$
(cf.\ \cite[Lemma 2.31]{BMPSTWW20}), {and} 
as $X$ is divisorially Cohen--Macaulay, $f_*\cO_Y(K_Y+\lceil D_Y \rceil)$ is Cohen--Macaulay. 
Thus, (1) holds. 
\end{proof}

Our key technical result is the following.

\begin{lemma} \label{lemma:pushforward-of-C}
Under exactly the same assumptions as in Lemma \ref{lemma:pushforward-of-Watanabe-F-split}, we have that 
\[
R^jf_* Q_{Y,K_Y+\Delta_Y,n} = 
\begin{cases}
Q_{X,K_X+\Delta,n} \qquad & (j=0)\\
0 \qquad & (j>0)
\end{cases}
\]
for every $n \in \Z_{>0}$. 
Moreover, $f_* Q^{\vee}_{Y,\Delta_Y,n} = Q^{\vee}_{X,\Delta,n}.$
\end{lemma}
\noindent The reader should be aware that $\Delta_Y$ need not be $p$-compatible.

\begin{proof}
First, we claim that 
\begin{equation}\label{e1:pushforward-of-C}
Rf_*\cB_{Y,p^i(K_Y+\Delta_Y{ )}} = \cB_{X,p^i(K_X+\Delta{)}}
\quad 
{\text{ for every } i\geq 0}.
\end{equation}
Consider the short exact sequence {(see (\ref{eq:definition-of-log-B}))}
\[
0 \to \cO_Y(p^{i}(K_Y+\Delta_Y)) \xrightarrow{F} F_* \cO_Y(p^{i+1}(K_Y+\Delta_Y)) \to \cB_{Y,p^i(K_Y+\Delta_Y{ )}} \to 0. 
\]
By applying $R^jf_*(-)$ to this exact sequence and using
Lemma \ref{lemma:pushforward-of-Watanabe-F-split}, we get that 
 $R^jf_*\cB_{Y,p^i(K_Y+\Delta_Y{)}}=0$ for every $j>0$. {Now, consider} 
the following commutative diagram in which each horizontal sequence is exact
\[
\begin{tikzcd}[column sep = small]
0 \arrow{r} & f_*\cO_Y(p^{i}(K_Y+\Delta_Y)) \arrow{r} & F_*f_* \cO_Y(p^{i+1}(K_Y+\Delta_Y)) \arrow{r} & f_*\cB_{Y,p^i(K_Y+\Delta_Y{ )}} \arrow{r} & 0\\
0 \arrow{r} & \cO_X(p^{i}(K_X+\Delta)) \arrow{u}{\cong} \arrow{r} & F_* \cO_X(p^{i+1}(K_X+\Delta)) \arrow{r} \arrow{u}{\cong} & \cB_{X,p^i(K_X+\Delta{)}} \arrow{r} \arrow{u} & 0.
\end{tikzcd}
\]
By the five lemma, 
$f_*\cB_{Y,p^i(K_Y+\Delta_Y{ )}} = \cB_{X,p^i(K_X+\Delta{ )}}$, {and} so (\ref{e1:pushforward-of-C}) holds.

By applying an analogous argument to the one above and using (\ref{eq:Bn-sequence}) and (\ref{e1:pushforward-of-C}), we get 
$Rf_*\cB_{Y,K_Y+\Delta_Y,n} = \cB_{X,K_X+\Delta,n}. 
$
{Then, by} repeating the same procedure {again but}  for (\ref{eq:C-quotient-sequence}), 
we finally obtain 
\begin{equation} \label{eq:first-pushforward-of-C}
Rf_* Q_{Y,K_Y+\Delta_Y,n} = Q_{X,K_X+\Delta,n}.
\end{equation}

To prove the second part of the lemma, we observe that
\begin{align*}
f_*\cHom_{\MO_Y}(Q_{Y,\Delta_Y,n}, \cO_Y) &\stackrel{(1)}{=} f_*\cHom_{\MO_Y}(Q_{Y,K_Y+\Delta_Y,n}, \omega_Y) \\
&\stackrel{(2)}{=} \cH^0Rf_*R\cHom_{\MO_Y}(Q_{Y,K_Y+\Delta_Y,n}, \omega_Y) \\
&\stackrel{(3)}{=} \cH^0R\cHom_{\MO_X}(Rf_*Q_{Y,K_Y+\Delta_Y,n}, \omega_X) \\
&\stackrel{(4)}{=} \cH^0R\cHom_{\MO_X}(Q_{X,K_X+\Delta,n}, \omega_X) \\
&\stackrel{(5)}{=} \cHom_{\MO_X}(Q_{X,\Delta,n}, \cO_X).
\end{align*}
 First, (1) and (5)  
 follow from (\ref{e-C-tensor})
 as these equalities can be checked on the regular loci by Remark \ref{remark:hom-reflexive}.  
Equality (2) holds by $Rf_*R\cHom_{\MO_Y}(Q_{Y,K_Y+\Delta_Y,n}, -) \cong 
R(f_* \circ \cHom_{\MO_Y}(Q_{Y,K_Y+\Delta_Y,n}, -))$ \cite[Ch.\ II, Lemma 7.3.2]{God58}. 
 Equality (3) follows from Grothendieck duality, because $X$ and $Y$ are Cohen--Macaulay. 
 Finally, (4) holds by (\ref{eq:first-pushforward-of-C}). 
\end{proof}

\begin{corollary}\label{Cor:min res preserve quasi-$F$-splitting}
\label{cor:pulling-back-quasi-F-splittings}
In the situation of Setting \ref{setting:most-general-foundations-of-log-quasi-F-splitting}, let $(X,\Delta)$ be a $p$-compatible klt $n$-quasi-$F$-split pair.
Let $f \colon Y \to X$ be a projective birational morphism of
integral normal divisorially Cohen-Macaulay schemes  
such that $\Delta_Y \geq 0$, {where $\Delta_Y$ is the log pullback of $\Delta$,  that is,}
$K_Y + \Delta_Y = f^*(K_X+\Delta)$ and $f_*\Delta_Y=\Delta$.

Assume that
\[
R^jf_*\cO_Y(K_Y + \rup{-p^i(K_Y+\Delta_Y)}) =0
\]
for all $i\geq 0$ and $j>0$. Then $(Y,\Delta_Y)$ is $n$-quasi-$F$-split.
\end{corollary}

\begin{proof}
By Lemma \ref{lemma:pushforward-of-C},  we have  that
\[
\Hom_{\MO_Y}(Q_{Y,\Delta_Y,n}, \cO_Y)= \Hom_{\MO_X}(Q_{X,\Delta,n}, \cO_X).
\]
Let $\alpha \colon Q_{X,\Delta,n} \to \cO_X$ be a splitting of $\Phi_{X, \Delta, n} \colon \cO_X \to Q_{X,\Delta,n}$.

The above equation shows that $\alpha$ extends to an $\MO_Y$-module homomorphism 
\[
\alpha_Y \colon Q_{Y,\Delta_Y,n} \to \cO_Y. 
\]
{In particular}, the following diagram is commutative:
\[
\begin{tikzcd}[column sep = large]
\arrow[d,dash,shift left=.1em] \arrow[d,dash,shift right=.1em] f_*\MO_Y \arrow{r}{f_*\Phi_{Y, \Delta_Y, n}} &  f_* Q_{Y,\Delta_Y,n} \arrow{r}{f_* \alpha_Y} & \arrow[d,dash,shift left=.1em] \arrow[d,dash,shift right=.1em] f_*\MO_Y\\
\ \MO_X \arrow{r}{\Phi_{X, \Delta, n}} & \ Q_{X,\Delta,n} {\arrow{u}{f^*}} \arrow{r}{\alpha} & \ \cO_Y.
\end{tikzcd}
\]
Here we use that $\Delta_Y \geq 0$, so that 
the domain of $\Phi_{Y, \Delta_Y, n}$ is indeed $\cO_Y$.
Since the lower horizontal composite arrow is the identity, so is the upper one: 
$f_* \alpha_Y \circ f_*\Phi_{Y, \Delta_Y, n} = {\rm id}$. 
Therefore, we get $ \alpha_Y \circ \Phi_{Y, \Delta_Y, n} (1_{\MO_Y(Y)}) = 1_{\MO_Y(Y)}$ 
for the identity element $1_{\MO_Y(Y)}  \in \Gamma(Y, \MO_Y)$, 
which implies $ \alpha_Y \circ \Phi_{Y, \Delta_Y, n}  ={\rm id}$. 
\end{proof}

\begin{remark}
In fact, provided that {similar} assumptions to these of Corollary \ref{cor:pulling-back-quasi-F-splittings} are satisfied, one can deduce from Lemma \ref{lemma:pushforward-of-C} that
$qS^0_n(X,\Delta; L) = qS^0_n(Y,\Delta_Y, L_Y)$ where $L$ is a Cartier divisor and $L_Y := f^*L$. 
\end{remark}

We also explain the converse implication.

\begin{proposition}
\label{prop:pushing-forward-quasi-F-splittings}
In the situation of Setting \ref{setting:most-general-foundations-of-log-quasi-F-splitting}, let $\Delta$ be an effective $\bQ$-divisor on $X$ such that $(X,\Delta)$ is $n$-quasi-$F$-split. Let $g \colon X \to Y$ be a projective birational morphism 
to an integral normal excellent scheme $Y$. Then $(Y,g_*\Delta)$ is $n$-quasi-$F$-split.
\end{proposition}
\begin{proof}
Let $Y' \subseteq Y$ be an open subset such that ${\rm codim}_Y(Y\, \setminus\, Y') \geq 2$ and $g|_{X'} \colon X' \to Y'$ is an isomorphism for $X' := g^{-1}(Y')$. 
Since $(X, \Delta)$ is $n$-quasi-$F$-split, 
also $(X', \Delta|_{X'}) = (Y', g_*\Delta|_{Y'})$ is $n$-quasi-$F$-split. 
It follows from Remark \ref{r-big-open-QFS} that {$(Y, g_*\Delta)$} is 
$n$-quasi-$F$-split. 
\end{proof}

Combining the previous results, we can show the following.
\begin{corollary} \label{corollary:log-pullback-quasi-F-split-dimension-two}
Let $f \colon Y \to X$ be a projective birational morphism 
of two-dimensional integral normal Noetherian $F$-finite $\mathbb{F}_p$-schemes admitting dualising complexes. 
Let $\Delta$ be a $\Q$-divisor on $X$ such that
$(X, \Delta)$ is klt, {$p$-compatible,} and {$n$-}quasi-$F$-split. Write $K_Y+\Delta_Y = f^*(K_X+\Delta)$ and suppose that $\Delta_Y \geq 0$.  

Then
$(Y, \Delta_Y)$ is {$n$-}quasi-$F$-split. 
\end{corollary}

\begin{proof}
Let $h \colon Z \to Y$ be the minimal resolution of $Y$ and let $g \colon Z \xrightarrow{h} Y \xrightarrow{f} X$ be the composition of $h$ and $f$. Write $K_Z + \Delta_Z = g^*(K_X+\Delta)$. Since $h^*K_Y-K_Z$ is effective by minimality of $h$, we get that $\Delta_Z \geq 0$. By Proposition \ref{prop:pushing-forward-quasi-F-splittings}, it is enough to show that $(Z,\Delta_Z)$ is $n$-quasi-$F$-split. Thus, we may replace $(Y,\Delta_Y)$ by $(Z,\Delta_Z)$ and assume from now on that $Y$ is regular. Then, 
\[
R^jf_*\cO_Y(K_Y + \rup{-p^i(K_Y+\Delta_Y)}) =  
R^jf_*\cO_Y(K_Y + \{p^i\Delta_Y\}-p^i(K_Y+\Delta_Y))= 0
\]
for $i\geq 0$ and $j>0$ by 
\cite[Proposition 3.2]{tanaka16_excellent}. 
Moreover,
$X$ is divisorially Cohen--Macaulay (Remark \ref{remark:divisorially-Cohen-Macaulay}).
Thus the assumptions of Corollary \ref{cor:pulling-back-quasi-F-splittings} are satisfied, and so $(Y,\Delta_Y)$ is $n$-quasi-$F$-split.
\end{proof}

We warn the reader that $(Y, \{p^i\Delta_Y\})$ in the above proof need not be klt, but, fortunately, \cite[Proposition 3.2]{tanaka16_excellent} only requires that $Y$ is regular and the coefficients of the boundary $\bQ$-divisor are smaller than $1$.

%% file: section4.tex
\section{Inversion of adjunction}
\subsection{Introduction} \label{ss:intro-to-lifting-quasi-stable-sections}
We explain the idea behind lifting quasi-$F$-stable sections (Theorem \ref{thm:intro-quasi-F-split-lifting}). 
For simplicity, we assume that $X$ is a smooth projective integral scheme 
over a perfect field of characteristic $p>0$ and $S$ is a smooth prime divisor on $X$.

Let $\cL$ be a line bundle on $X$. We defined the space of {\cora $n$-}quasi-$F$-stable sections as
\[
qS^0_n(X, \cL) = \image( H^0(X, Q^\vee_{X,n} \otimes \cL) \to H^0(X,\cL)),
\]
and its adjoint variant as
\[
qS^0_{n,\adj}(X, S; \cL) = \image(H^0(X, (Q^S_{X,S,n})^\vee \otimes \cL) \to H^0(X,\cL)),
\]
where $Q^S_{X,S,n} \cong Q^S_{X,n} \otimes_{\MO_X} \MO_X(S)$ holds by (\ref{e-C-tensor2}) 
and $Q^S_{X,n}$ is the pushout in
\begin{center}
\begin{tikzcd}
W_n\cI_S \arrow{r}{F} \arrow{d}{R^{{\cora n-1}}} & F_* W_n\cI_S \arrow{d}  \\
\cI_S \arrow{r}{\Phi^S_{X, n}} & Q^S_{X,n} \arrow[lu, phantom, "\usebox\pushoutdr" , very near start, yshift=0em, xshift=0.6em, color=black]. 
\end{tikzcd}
\end{center}
One can easily check (we provide details later) that there exists the following commutative diagram 
in which each horizontal sequence is exact: 
\begin{center}
\begin{tikzcd}
0 \arrow{r} & Q^S_{X,n} \arrow{r} & Q_{X,n}  \arrow{r} & Q_{S,n} \arrow{r} &   0 \\
0 \arrow{r} & \cO_X(-S) \arrow{u}{} \arrow{r} & \cO_X \arrow{u}{} \arrow{r} & \cO_S \arrow{u}{} \arrow{r} & 0.
\end{tikzcd}
\end{center}
A simple calculation with $\cExt$ shows that applying $\cHom_{\MO_X}(- , \cO_X(-S))$ yields:
\begin{center}
\begin{tikzcd}
0 \arrow{r} & (Q_{X,n})^\vee \otimes \MO_X(-S) \arrow{d}{} \arrow{r} & (Q^S_{X,S,n})^\vee \arrow{d}{}  \arrow{r} & (Q_{S,n})^\vee \arrow{d}{} \arrow{r} &   0 \\
0 \arrow{r} & \cO_X(-S)  \arrow{r} & \cO_X  \arrow{r} & \cO_S  \arrow{r} & 0.
\end{tikzcd}
\end{center}
Finally tensoring by $\cL$ gives the map
\[
qS^0_{n, \adj}(X, S; \cL) \to qS^0_n(S, \cL|_S)
\]
featured in Theorem \ref{thm:intro-quasi-F-split-lifting}. We are now ready to provide the following. 

\begin{proof}[Sketch of the proof of Theorem \ref{thm:intro-quasi-F-split-lifting}]
Set $d := \dim X$, and $Q_{X,n}(-A) := Q_{X,n} \otimes \cO_X(-A)$. 
By the above diagram, it is enough to show that 
\[
H^1(X, (Q_{X,n})^{\vee} \otimes \MO_X(L) \otimes \MO_X(-S))=0. 
\]
By Serre duality, this is equivalent to 
$H^{d-1}(X, Q_{X,n}(-A))=0$.
{\cora 
Recall that we assume that 
\begin{enumerate}
    \item $H^1(X, \cO_X(K_X + p^iA))=0$ and 
    \item $qS^0_{n-i}(X, \cO_X(K_X + p^iA)) = H^0(X, \cO_X(K_X+p^iA))$. 
\end{enumerate}
for all $i \geq 1$.} 
We have the following implications:
\begin{align*}
H^{d-1}(X, {Q_{X,n}(-A)})\!=\! 0\! 
\xLeftarrow{\text{ {Assumpt.}\ \!(1) and } (\ref{eq:intro-C-restriction-sequence}) }
 &H^{d-1}(X, (F_*(Q_{X,n-1}/\cO_X)) \!\otimes\! \MO_X({-A}))\!=\! 0  \\
\xLeftarrow{\small \text{ { Assumpt.}\ \!(2) and  (\ref{eq:obvious-sequence-for-C})}}
&H^{d-1}(X, Q_{X,n-1} \otimes \MO_X(-pA))\!=\!0 \\
\xLeftarrow{\text{repeat the above } n-1 \text{ times}}
&H^{d-1}(X, Q_{X,{1}} \otimes \MO_X(-p^{n-1}A))\!=\!0,
\end{align*}
where the last vanishing {follows from Assumption (1) as $Q_{X,1}=F_*\cO_X$}. 
The second implication used the long exact sequence of cohomology associated to
\begin{equation} \label{eq:obvious-sequence-for-C}
0 \to \cO_X \to Q_{X,n-1} \to Q_{X,n-1}/\cO_X \to 0     
\end{equation}
and the fact that $H^d(X,\cO_X(-pA)) \to H^d(X,Q_{X,n-1}(-pA))$ is injective. This injectivity follows by applying Grothendieck duality to the map \[H^0(X,Q^{\vee}_{X,n-1}(K_X+pA)) \to H^0(X,\cO_X(K_X+pA)),\] which is surjective by Assumption (2).
\end{proof}

Theorem \ref{thm:intro-inversion-of-adjunction} follows by repeatedly applying Theorem \ref{thm:intro-quasi-F-split-lifting} and using the fact that $-kS$ is sufficiently ample for $k\gg 0$ so that it can be used to kill cohomology via the Serre vanishing theorem. Of course, in practice, {the argument} is much more technically involved.

\subsection{The construction of the restriction maps for $qS^0$}
We start by showing that the restriction of $Q_{X,\Delta,n}$ to $S$ is equal to $Q_{S,\Delta|_S,n}$.
\begin{proposition} \label{prop:restrictiong-ses-for-C}
In the situation of Setting \ref{setting:most-general-foundations-of-log-quasi-F-splitting}, suppose that $X$ is divisorially Cohen-Macaulay. Let $S$ be a normal prime divisor on $X$ and 
let $\Delta$ be a {$\bQ$-Cartier} $\bQ$-divisor on $X$ such that 
$S \not\subseteq \Supp\, \Delta$ and
$(X,S+\{p^i\Delta\})$ is plt for every $i \geq 0$.

Then restricting $Q_{X,\Delta, n}$ to $S$ yields the following short exact sequence:
\[
0 \to Q^{S}_{X,\Delta,n} \to Q_{X,\Delta,n} \to Q_{S,\Delta_S,n} \to 0,
\]
where $\Delta_S = \Delta|_S$.
\end{proposition}
\noindent Here,  $\Delta|_S$ is defined as in Definition \ref{definition:restricting-divisor}. Also note that we do not assume that $\rdown{\Delta}=0$ as we will often set $\Delta = K_X+S+B - L$ for a plt pair $(X,S+B)$ and $L$ being the Weil divisor whose quasi-$F$-stable sections we want to study.

\begin{proof}
First, we construct the short exact sequence
\begin{equation} \label{eq:witt-restriction}
0 \to W_n \cI_S(\Delta) \to W_n\cO_X(\Delta) \xrightarrow{\res} W_n\cO_S(\Delta_S) \to 0.
\end{equation}
The $W_n\cO_X$-module homomorphism $\res \colon W_n\cO_X(\Delta) \to W_n\cO_S(\Delta_S)$ is defined by the formula 
\[
(a_0, \ldots, a_{n-1}) \mapsto (a_0|_S, \ldots, a_{n-1}|_S),
\]
where $a_i \in \Gamma(U, F^i_*\cO_X(p^i\Delta))$ for some open subset $U \subseteq X$, and $a_i|_S := \res(a_i)$ for {the} surjective homomorphism $\res \colon F^i_*\cO_X(p^i\Delta) \to F^i_*\cO_S(p^i\Delta_S)$ induced from Proposition \ref{prop:basic_restriction}. 
In particular, $\res \colon W_n\cO_X(\Delta) \to W_n\cO_S(\Delta_S)$
is well-defined and surjective.
Last, $\kernel(\res) = W_n \cI_S(\Delta)$, because $a_i|_S = 0$ if and only if $a_i \in \Gamma(U, F^i_*\cO_X(-S+p^i\Delta))$.

By (\ref{eq:witt-restriction}), 
we get the following commutative diagram in which each horizontal sequence is exact
\begin{center}
\begin{tikzcd}
0 \arrow{r} & W_{n-1} \cI_S({p}\Delta) \arrow{d}{FV} \arrow{r} & W_{n-1}\cO_X({ p}\Delta) \arrow{d}{FV} \arrow{r}{\res} & W_{n-1}\cO_S({ p}\Delta_S) \arrow{d}{FV} \arrow{r} & 0 \\
0 \arrow{r} & W_n\cI_S({ p}\Delta) \arrow{r} & W_n\cO_X({ p}\Delta) \arrow{r}{\res} & W_n\MO_S({p}\Delta_S) \arrow{r} & 0.
\end{tikzcd}
\end{center} 
Since all the vertical arrows are injective, the snake lemma yields the required exact sequence (cf.\ (\ref{e-big-BC-diagram2})).
\end{proof}

We will now construct the restriction map for $qS^0$ that is featured in Theorem \ref{thm:intro-quasi-F-split-lifting}. With exactly the same assumptions as in Proposition \ref{prop:restrictiong-ses-for-C},  
consider the following commutative diagram 
in which each horizontal sequence is exact:
\[
\begin{tikzcd}
0 \arrow{r} & Q^{S}_{X,\Delta,n}  \arrow{r} & Q_{X,\Delta,n} \arrow{r} & Q_{S,\Delta_S,n} \arrow{r} & 0 \\
0 \arrow{r} & \cO_X(\Delta-S) \arrow{r} 
\arrow{u}{{\Phi^S_{X, \Delta, n}}} & \cO_X(\Delta) \arrow{r} \arrow{u}{{\Phi_{X, \Delta, n}}} & \cO_S(\Delta_S) \arrow{r} \arrow{u}{{\Phi_{S, \Delta_S, n}}} & 0.
\end{tikzcd}
\]
As we shall verify below, applying $\cHom_{{\MO_X}}(- , \omega_X)$ yields:
\begin{equation}\label{e-C-adj-diag}
\begin{tikzcd}[column sep = small]
0 \arrow{r} & (Q_{X,\Delta-K_X,n})^\vee \arrow{d}{} \arrow{r} & (Q^S_{X,\Delta-K_X, n})^\vee \arrow{d}{}  \arrow{r} & (Q_{S,\Delta_S-K_S, n})^\vee \arrow{d}{} \arrow{r} &   (\star) \\
0 \arrow{r} & \cO_X(K_X\!-\!\rdown{\Delta})  \arrow{r} & \cO_X(K_X\!+\!S\!-\!\rdown{\Delta})  \arrow{r} & \cO_S(K_S\!-\!\rdown{\Delta_S})  \arrow{r} & 0,
\end{tikzcd}
\end{equation}
with $(\star)$ being zero if $(X,\{p^i\Delta\})$ is naively keenly $F$-pure for every $i\geq 0$.

\begin{proof}[Verification] 
We check that the upper {horizontal sequence} dualises in this fashion; 
the case of the lower {horizontal sequence} is completely analogous, standard, and left to the reader 
{\cora (note that the surjectivity of 
$\cO_X(K_X + S - \rdown{\Delta})  
\to \cO_S(K_S - \rdown{\Delta_S})$  
follows from the fact that $X$ is divisorially Cohen-Macaulay)}. First,  $\cHom_{{\MO_X}}(\cF, \cO_X)$ and $\cHom_{{\MO_X}}(\cF, \omega_X)$ are reflexive for every coherent sheaf $\cF$, and so we immediately get
\begin{align*}
\cHom_{{\MO_X}}(Q_{X,\Delta,n}, \omega_X) &\cong (Q_{X,\Delta-K_X,n})^{\vee} \\
\cHom_{{\MO_X}}(Q^S_{X,\Delta,n}, \omega_X) &\cong (Q^S_{X,\Delta-K_X,n})^{\vee},
\end{align*}
{because} these {isomorphisms} can be checked on the regular locus of $X$ {by Remark \ref{remark:hom-reflexive}} {(cf.\ (\ref{e-C-tensor}), (\ref{e-C-tensor2}))}. For the induced closed immersion $i \colon S \to X$, 
we have $\cHom_{\MO_X}(i_*Q_{S,\Delta_S,n}, \omega_X)=0$, 
as $i_*Q_{S,\Delta_S,n}$ is torsion. {Set} 
$d := \dim X$. {Then}
\begin{align*}
\cExt^1_{{\MO_X}}(i_*Q_{S,\Delta_S,n}, \omega_X) &=
\cH^{-d+1}R\cHom_{{\MO_X}}(i_*Q_{S,\Delta_S,n}, \omega^{\mydot}_X) \\
&= i_*\cH^{-d+1}R\cHom_{{\MO_S}}(Q_{S,\Delta_S,n}, \omega^{\mydot}_S)\\
&= i_*\cHom_{{\MO_S}}(Q_{S,\Delta_S,n}, \omega_S) \\
&= (Q_{S,\Delta_S-K_S,n})^{\vee}.
\end{align*}
by Grothendieck duality, where $\omega_S = \cH^{-d+1}(\omega^{\mydot}_S)$ is the dualising sheaf of $S$. 
Here, the first and third equalities hold since $X$ and $S$ are Cohen-Macaulay.
The last equality follows from Remark \ref{remark:hom-reflexive}.

Finally, assuming that $(X,\{p^i\Delta\})$ is {naively keenly $F$-pure} for every $i\geq 0$, we get that 
\[
\cExt^1_{{\MO_X}}(Q_{X,\Delta,n}, \omega_X) = 0
\]
as $Q_{X,\Delta,n}$ is Cohen-Macaulay by Proposition \ref{prop:C-cohen-macaulay}. 
\end{proof}

We are ready to define the restriction map for $qS^0$.

\begin{definition} \label{definitin:restriction-for-qS^0} 
In the situation of {Setting \ref{setting:most-general-foundations-of-log-quasi-F-splitting}},
let $(X,S+B)$ be a  divisorially Cohen-Macaulay plt pair, 
where
$S$ is a normal prime divisor and 
$B$ is a $\Q$-divisor such that $\rdown{B}=0$.
{\cora In particular, we have $S \not \subseteq \Supp B$.} 
Suppose that $(X,S+\{p^iB\})$ is plt for every $i\geq 0$ and write $K_S + B_S = (K_X+S+B)|_S$. 

Let $L$ be a $\bQ$-Cartier  Weil divisor {such that $S \not\subseteq \Supp L$}. Recall that 
\begin{align*}
qS^0_{n,\adj}(X,S+B;L) &= \image\big(
H^0(X,(Q^S_{X,S+B-L,n})^{\vee}) \to  
H^0(X, \cO_X(L)\big)\ {\rm and } \\
qS^0_n(S, B_S;L|_S) &= \image\big(
H^0(S,(Q_{S,B_S -L|_S,n})^{\vee}) \to H^0(S, \cO_X(\lceil L|_S-B_S\rceil))\big).
\end{align*}
Then we {obtain} a map
\[
qS^0_{n, \adj}(X, S+B; L) \to qS^0_n(S, B_S; L|_S)
\]
induced by plugging $\Delta=K_X+S+B-L$ into the above diagram (\ref{e-C-adj-diag}) to get
\begin{equation} \label{eq:diagram-res-C}
\begin{tikzcd}[column sep = small]
0 \arrow{r} & (Q_{X,S+B-L,n})^\vee \arrow{d}{} \arrow{r} & (Q^S_{X,S+B-L, n})^\vee \arrow{d}{}  \arrow{r} & (Q_{S,B_S-L|_S, n})^\vee \arrow{d}{} \arrow{r} &   (\star) \\
0 \arrow{r} & \cO_X(L-S)  \arrow{r} & \cO_X(L)  \arrow{r} & \cO_S(\rup{L|_S - B_S})  \arrow{r} & 0.
\end{tikzcd}
\end{equation}
Note that $(\star)$ is zero when $(X,\{p^iB\})$ is {naively keenly $F$-pure}
for every $i \geq 0$. 
\end{definition}
Under our assumptions, $\rup{L|_S- B_S} = \rdown{L|_S}$ (see Lemma  \ref{l:special-cases-of-L}), and so the lower horizontal row in Diagram (\ref{eq:diagram-res-C}) does not contradict Proposition \ref{prop:basic_restriction}.

\begin{remark} \label{remark:implicit-replacement-of-L}
In the above definition, in order to make $(K_X+S+B)|_S$ well-defined, we always, implicitly, pick a canonical divisor $K_X$ so that $S \not \subseteq \Supp\,(K_X+S+B)$. Specifically, this is achieved by using prime avoidance to pick $K_X=-S+D$ for a Weil divisor $D$ such that $S \not \subseteq \Supp D$.

More importantly, we will apply Definition \ref{definitin:restriction-for-qS^0} even when $S \subseteq \Supp L$. In that case, we implicitly find $L' \sim L$ such that $S \not \subseteq \Supp L'$ and then replace $L$ by $L'$. 
As making such replacements in the body of the text would make the proofs and  notation impenetrable, we elected \emph{not} to do it. 
Instead, we added footnotes and  remarks that explicate and justify each such replacement.  
We also emphasise that $\{L\}=\{L'\}$, and so all the sheaves in (\ref{eq:diagram-res-C}) are independent of the choice of $L'$ up to a (compatible) isomorphism.
\end{remark}

\begin{remark}
In the spirit of the above argument, we set up key exact sequences describing $Q_{X,\Delta,n}^\vee$. Specifically, with the same assumptions as in Proposition \ref{prop:restrictiong-ses-for-C}, we suppose additionally that $(X,\{p^i\Delta\})$ 
is 
{naively keenly $F$-pure} for every $i\geq 0$. Then we apply $\cHom_{{\MO_X}}(-, \omega_X)$ to (\ref{eq:C-quotient-sequence}) and (\ref{eq:C-restriction-sequence}), in order to get 
\begin{align} \label{eq:dual-C-quotient-sequence}
&0 \to \cB^{\vee}_{X,\Delta-K_X,n}  \to Q^{\vee}_{X,\Delta-K_X,n} \to \cO_X(K_X - \rdown{\Delta}) \to 0 \\
\label{eq:dual-C-restriction-sequence}
&0 \to  F_* \cO_X(K_X - \rdown{p\Delta}) \to Q^{\vee}_{X,\Delta-K_X,n} \to  F_*\cB^{\vee}_{X,p\Delta-K_X,n-1} \to 0.
\end{align}
The above sequences are exact as $\cB_{X,\Delta, n}$ and $F_*\cO_X(p\Delta)$ are Cohen-Macaulay (see Proposition \ref{prop:C-cohen-macaulay}), and so $\cExt^1_{{\MO_X}}(B_{X,\Delta,n}, \omega_X)=0$  and $\cExt^1_{{\MO_X}}(F_*\cO_X(p\Delta), \omega_X)=0$.
\end{remark}

\subsection{Lifting sections}

We are ready to prove Theorem \ref{thm:intro-quasi-F-split-lifting} in full generality.
\begin{theorem} \label{thm:lifting-sections-basic}
In the situation of Setting \ref{setting:foundations-of-log-quasi-F-splitting},
let $(X,S+B)$ be a divisorially Cohen-Macaulay  plt pair, where 
$S$ is a normal prime divisor and 
$B$ is a $\Q$-divisor such that 
$\rdown{B}=0$. Suppose that $(X,S+\{p^iB\})$ is plt and $(X,\{p^iB\})$ is {naively keenly $F$-pure} for every $i\geq 0$. 
Let $L$ be a $\bQ$-Cartier Weil divisor such that
\begin{enumerate}
    \item $H^1(X, \cO_X(K_X + \lceil p^{i}A \rceil))=0$ and
    \item $qS_{n-i}^0(X,\{p^iB\}; K_X + \lceil p^iA \rceil) =H^0(X,\cO_X(K_X + \lceil p^iA \rceil))$
\end{enumerate}
for all $i \geq 1$, where $A := L - (K_X+S+B)$.
Then the map 
\[
qS^0_{n,\adj}(X,S+B;L) \to qS^0_n(S,B_S; L|_S)
\]
from Definition \ref{definitin:restriction-for-qS^0} is surjective (cf.\ Remark \ref{remark:implicit-replacement-of-L}).
\end{theorem}

\noindent Note that $K_X + \lceil p^iA \rceil = K_X + \{p^iB\} + p^iA$. We do not assume that $A$ is ample necessarily, though that will usually be the case in applications.

In alignment with Remark \ref{remark:implicit-replacement-of-L}, what we formally claim in Theorem \ref{thm:lifting-sections-basic} is that  $qS^0_{n,\adj}(X,S+B;L') \to qS^0_n(S,B_S; L'|_S)$ is surjective for every (equivalently, one) {Weil} divisor $L' \sim L$ such that $S \not \subseteq \Supp L'$.

\begin{proof}
Up to replacing $L$ by a linearly equivalent {Weil} divisor, we may assume that $S \not \subseteq \Supp L$. This reduction process does not affect Assumptions (1) and (2).

By the diagram {\cora \eqref{eq:diagram-res-C}}, 
it is enough to show that  
\[
 H^1(X, Q^{\vee}_{X,S+B-L,n}) = 0.
\]
Recall that $A = L - (K_X+S+B)$, and so $S+B-L = -K_X-A$. We will prove by descending induction on $i \geq 0$ that
\begin{equation} \label{eq:basic-lifting-inductive-step}
H^{1}(X, Q^{\vee}_{X,-K_X-p^iA,n-i}) = 0.
\end{equation}
This is clear for the base of the induction {$i=n-1$} as then 
\[
Q^{\vee}_{X,-K_X-p^{n-1}A,1} = \cHom_{\cO_X}(F_*\cO_X(-p(K_X+p^{n-1}A)), \cO_X) \cong F_*\cO_X(K_X + \rup{p^{n}A}) 
\]
by Grothendieck duality, and so the statement follows by Assumption (1).
Thus we may assume that 
\begin{equation} \label{eq:easy-lifting-sections-induction}
H^{1}(X, Q^{\vee}_{X,-K_X-p^{i+1}A,n-i-1}) = 0
\end{equation}
for $0 \leq i \leq {n-2}$ and aim for showing \eqref{eq:basic-lifting-inductive-step}.

The short exact sequence \eqref{eq:dual-C-restriction-sequence} for $\Delta=-p^iA$ specialises to:
\[
0 \to  F_* \cO_X(K_X + \rup{p^{i+1}A}) \to Q^{\vee}_{X,-K_X-p^iA,n-i} \to  F_*\cB^{\vee}_{X,-K_X-p^{i+1}A,n-i-1} \to 0.
\]
and so by {\cora our} assumptions, it is enough to show that
\[
H^1(X,\cB^{\vee}_{X,-K_X-p^{i+1}A,n-i-1})=0.
\]
The short exact sequence \eqref{eq:dual-C-quotient-sequence} for $\Delta=-p^{i+1}A$ specialises to:
\[
0 \to \cB^{\vee}_{X,-K_X-p^{i+1}A,n-i-1}  \to Q^{\vee}_{X,-K_X-p^{i+1}A,n-i-1} \to \cO_X(K_X + \rup{p^{i+1}A}) \to 0.
\]
Thus $H^1(X, \cB^{\vee}_{X,-K_X-p^{i+1}A,n-i-1}) =0$ is valid, because 
\begin{itemize}
    \item $H^1(X, Q^{\vee}_{X,-K_X-p^{i+1}A,n-i-1})=0$ by {the induction hypothesis} \eqref{eq:easy-lifting-sections-induction},
    \item \vspace{0.2em} $H^0(X,Q^{\vee}_{X,-K_X-p^{i+1}A,n-i-1}) \to H^0(X, \cO_X(K_X + \rup{p^{i+1}A}))$ is surjective as  
    \[
    qS_{n-i-1}^0(X,\{p^{i+1}B\}; K_X + \lceil p^{i+1}A \rceil) =H^0(X,\cO_X(K_X + \lceil p^{i+1}A \rceil)).
    \]
\end{itemize}
The second observation follows automatically from Definition \ref{definition:quasi-F-stable-sections} after we observe that $-K_X-p^{i+1}A = \{p^{i+1}B\} - (K_X+\lceil p^{i+1}A \rceil)$. The proof is concluded.
\end{proof}

\subsection{Inversion of adjunction in the log Calabi--Yau case}

We apply the above theorem to show that inversion of adjunction holds 
in the log Calabi--Yau case when $X$ is projective over a field.

\begin{theorem} \label{thm:log-trivial-inversion-of-adjunction}
Let $(X,S+B)$ be a divisorially Cohen-Macaulay 
{$p$-compatible} plt pair which is projective over an {$F$-finite} field $k$ of characteristic $p>0$, 
 where $S$ is a normal prime divisor and 
$B$ is a $\Q$-divisor {satisfying} 
$\rdown{B}=0$. 
Suppose that
\begin{enumerate}
    \item $K_X+S+B \equiv 0$,
    \item $(X,B)$ is {naively keenly $F$-pure},
    \item $(S,B_S)$ is $n$-quasi-$F$-split, where $K_S+B_S = (K_X+S+B)|_S$, and
    \item $H^1(X, \cO_X(K_X + \rup{-p^i(K_X+S+B)}))=0$ for all $i \geq 1$.
\end{enumerate}
Then $(X,S+B)$ is purely $n$-quasi-$F$-split.
\end{theorem}

\begin{proof}
This follows from  Theorem \ref{thm:lifting-sections-basic} for $L=0$ {(cf.\ Lemma \ref{lem:dual-criterion-for-pure-quasi-F-splitting})} if we verify that
\[
qS^0_{n-i}(X, \{p^iB\}; K_X + \rup{-p^i(K_X+S+B)}) = H^0(X, \cO_X(K_X + \rup{-p^i(K_X+S+B)})) 
\]
{for all $i \geq 1$. }
But this is clear, because 
\[
K_X + \rup{-p^i(K_X+S+B)} 
\equiv K_X + \{p^iB\} \leq K_X+B \equiv -S, 
\]
and so 
\[
H^0(X, \cO_X(K_X + \rup{-p^i(K_X+S+B)}))=0. \qedhere
\]
\end{proof}

\begin{corollary} \label{cor:cubic-surface}
Let $k$ be an 
algebraically closed field 
of characteristic $p>0$ and 
let $X$ be a smooth {del Pezzo} surface over $k$. 
Then 
a general member $C \in |-K_X|$ is a smooth elliptic 
 curve 
and $(X,C)$ is purely $2$-quasi-$F$-split. 
\end{corollary}

\begin{proof}
It follows from \cite[Theorem 1.4]{Kawakami-Nagaoka} that a general member $C$ is a smooth elliptic curve. 
Then we have $K_X+C \sim 0$ and $H^1(X, \MO_X(K_X))=0$. 
Since $C$ is $2$-quasi-$F$-split (see {Remark \ref{remark:elliptic-curves-quasi-F-split}}), 
$(X, C)$ is purely $2$-quasi-$F$-split by Theorem \ref{thm:log-trivial-inversion-of-adjunction}. 
\end{proof}

\begin{remark}
If $k$ is of characteristic $p>3$, then in Corollary \ref{cor:cubic-surface} it is enough to assume that $-K_X$ is nef and big (as opposed to ample). Indeed, 
then a general member $C \in |-K_X|$ is smooth by \cite[Theorem 1.4]{Kawakami-Nagaoka}, and so the same argument as above shows that  
$(X, C)$ is purely $2$-quasi-$F$-split.
\end{remark}

\subsection{Inversion of adjunction for  anti-semi-ample divisors}

In this subsection, we prove Theorem \ref{thm:intro-inversion-of-adjunction} in full generality. We start by showing the following lemma.

\begin{lemma} \label{lem:auxiliary-inversion-of-adjunction-anti-ample}
In the situation of Setting \ref{setting:foundations-of-log-quasi-F-splitting}, suppose that $R$ is a local ring and $H^0(X, \MO_X)=R$. 
Let $(X,S+B)$ be a divisorially Cohen-Macaulay  
log pair, where 
$S$ is a normal prime divisor and 
$B$ is a $\Q$-divisor such that $\rdown{B}=0$ and $(X,S+\{p^iB\})$ is plt for every $i\geq 0$. Set $K_S + B_S = (K_X+S+B)|_S$ and let $L$ be a $\bQ$-Cartier Weil divisor on $X$.
Suppose that the following hold. 
\begin{enumerate}
\item $-S$ is $\Q$-Cartier and semiample.
\item \vspace{0.05em}$(S,B_S)$ is $m$-quasi-$F$-split.
\item \vspace{0.05em}
$qS^0_{m,\adj}(X,S+{B}; L-kS) \!\to\! qS^0_m(S,{B}; (L-kS)|_S)$ 
is surjective for all $k \geq 0$. 
\end{enumerate}
Then
\[
qS^0_{m,\adj}(X,S+{B}; L) = H^0(X,\cO_X(L)).
\]
\end{lemma}

\noindent
Here, Assumption (3) should be formally read as follows (cf.\ Remark \ref{remark:implicit-replacement-of-L}): for every $k \geq 0$ and for {every} (equivalently, for one) Weil {divisor} $D_{k} \sim L-kS$ such that $S \not \subseteq \Supp D_{k}$, we have that $qS^0_{m, \adj}(X,S+{B}; D_{k}) \to qS^0_{m}(S, {B_{S}}; D_{k}|_S)$ is surjective.

\begin{proof}
Let $\m$ be the maximal ideal of $R$. 
We have $\m \in \pi(S)$. 
\\

\noindent \textbf{Reduction to \eqref{e1-l-anti-ample}.} 
We may assume that 
$\Supp L$ contains no components of $\Supp(S+B)$. 
We claim that {\cora it suffices to show that} 
\begin{equation}\label{e1-l-anti-ample}
\mathfrak{m}H^0(X,\cO_X(L)) + qS^0_{m,\adj}(X,S+B; L) = H^0(X,\cO_X(L)).
\end{equation}
Assuming {\eqref{e1-l-anti-ample}}, we get $\mathfrak{m}M = M$, where 
\[
M := H^0(X,\cO_X(L))/qS^0_{m,\adj}(X,S+B; L) 
\]
is a finitely generated $R$-module. Hence $M=0$ by Nakayama's lemma (see \cite[Proposition 2.6]{AM69} or \cite[\href{https://stacks.math.columbia.edu/tag/00DV}{Tag 00DV}]{stacks-project}. This concludes the proof of the lemma under Assumption \eqref{e1-l-anti-ample}.\\

\noindent \textbf{Reduction to \eqref{eq:second-induction}.} 
We are left with showing {(\ref{e1-l-anti-ample})}. First, we observe that
\begin{equation} \label{eq:integral-closure}
 H^0(X,\cO_X(L-kS)) \subseteq \mathfrak{m}H^0(X,\cO_X(L)),
\end{equation}
for $k\gg 0$. 
Indeed, 
there exist an integer $l>0$ and a projective morphism $f \colon X \to X'$ 
such that $X'$ is a projective $R$-scheme and $\MO_X(-lS) \cong f^*\cH$, where $\cH$ is a very ample invertible sheaf on $X'$. 
{Note that} $H^0(X,\cO_X(-lS)) \subseteq \m H^0(X,\cO_X)$ {as} $H^0(X, \MO_X) =R$. 
Let $s_1,\ldots, s_r \in H^0(X', \cH)$ 
be a {collection} of sections which globally generate $\cH$. 
Let 
\[
0 \to K \to  \MO_{X'}^{\oplus r}  
\to \cH \to 0
\]
be a short exact sequence of coherent locally free sheaves, where $\cO_{X'}^{\oplus r} \to \cH$ is a map on local sections given by $(f_1,\ldots,f_r) \mapsto {\sum_{i=1}^r f_is_i}$ and $K := \kernel({\MO_{X'}^{\oplus r} \to \cH})$. The map $\cO_{X'}^{\oplus r} \to \cH$ is surjective, because $s_1,\ldots, s_r$ globally generate $\cH$. 

Pick any integer $k \gg 0$ and write $k = lq +r$, {where $q$ and $r$ are integers} satisfying $q \gg 0$ and $0 \leq r <l$. 
By {tensoring the above sequence} with 
$f_*\MO_X(L-(k-l)S) = f_*\MO_X(L {- (r-l)}S) \otimes \cH^{\otimes q}$ and {applying} $H^0(X', -)$, {we get the following homorphism}
\[
{H^0(X, \cO_X(L-(k-l)S))^{\oplus r}}\to H^0(X, \cO_X(L-kS)), \qquad 
{ 
(g_1, ..., g_r) \mapsto \sum_{i=1}^r g_is_i. }
\]
{This homomorphism} is surjective for $k \gg l$
as $H^1(X', K \otimes f_*\MO_X(L {- (r-l)}S) \otimes \cH^{\otimes q})=0$ by Serre's vanishing. This immediately implies (\ref{eq:integral-closure}).

Now, in view of (\ref{eq:integral-closure}), {in order to prove (\ref{e1-l-anti-ample})}, it is enough to show that
 \begin{equation} \label{eq:second-induction}
 H^0(X,\cO_X(L-kS)) + qS^0_{m,\adj}(X,S+B; L) = H^0(X,\cO_X(L))
 \end{equation}
for all $k \geq 0$. \\

\noindent \textbf{Proof of  (\ref{eq:second-induction}).} 
We shall prove (\ref{eq:second-induction}) by increasing induction on $k$. Since the base case of the induction $k=0$ is trivial, we may assume that  (\ref{eq:second-induction})
is true for some $k\geq 0$ and aim for showing that
\begin{equation} \label{eq:l-IOA-antiample-k+1}
H^0(X,\cO_X(L-(k+1)S)) + qS^0_{m,\adj}(X,S+B; L) = H^0(X,\cO_X(L)).
\end{equation}
We start by making two observations. First, note that
\begin{equation} \label{eq:inclusion_of_im_psi-new}
qS^0_{m,\adj}(X,S+B; L) \supseteq qS^0_{m,\adj}(X,S+B; L-kS)    
\end{equation}
 Second, observe that 
Proposition \ref{prop:qS^0-for-quasi-F-split} yields
\begin{equation}\label{e-IOAlem-QFS-S}
qS^0_m(S, {B_{S}}; (L-kS)|_S)
 = H^0(S, \cO_S(\rup{(L-kS)|_S - {B_{S}}})). 
\end{equation}
{\cora Here, formally, we state that $qS^0_m(S, B_{S}; D_k|_S)
 = H^0(S, \cO_S(\rup{D_k|_S - B_{S}}))$ for every (equivalently, one) Weil divisor $D_k \sim L-kS$ such that $S \not \subseteq \Supp D_k$.}
Here we used that $(S, \{{B_{S}} - (L-kS)|_S\})$ is $m$-quasi-$F$-split, which is true by Assumption (2) and the following claim.

\begin{claim}
$\{{B_{S}} - (L-kS)|_S\} \leq B_S$.
\end{claim}
\begin{proof}
Pick an arbitrary prime divisor $\Gamma \subseteq S$ and let $\xi$ be its generic point. To prove the claim, it is enough to show the following inequality of the coefficients
\begin{equation} \label{eq:inequality-on-coeffs}
\mathrm{coeff}_{\Gamma}(B_S) \geq \mathrm{coeff}_{\Gamma}(\{{B_{S}} - (L-kS)|_S\}).
\end{equation}
Let $(-)_{\xi}$ denote the localisation at $\xi$. Then $(X_{\xi}, S_{\xi} + B_{\xi})$ is a two-dimensional plt pair. Denote by  $N$ the determinant of the intersection matrix of the minimal resolution of $X_{\xi}$. 

After replacing $X$ by a small open neighbourhood of $\xi \in X$, 
we get that
\begin{align*}
(K_X+S)|_S &= K_S + \Big(1-\frac{1}{N}\Big)\Gamma \qquad \text{ and }\qquad (L-kS)|_S = \frac{\ell_k}{N}\Gamma
\end{align*}
for some 
$\ell_k \in \bZ$ by Lemma \ref{lem:everything-about-plt-surface-pairs}(2)(3) applied to the plt pair $(X_{\xi},S_{\xi})$. 
By replacing $X$ by a small neighbourhood again, 
we get $B_S = (1-\frac{1}{N}+\lambda)\Gamma$ for some $\lambda \in \Q$. 
We have that $0 \leq \lambda <\frac{1}{N}$ because $\rdown{B_S}=0$ and $(K_X+S)|_S \leq (K_X+S+B)|_S = K_S + B_S$.
Then
\[
\mathrm{coeff}_{\Gamma}(B_S)
= 1 - \frac{1}{N} + \lambda 
\geq 
\left\{ 1 - \frac{1}{N} + \lambda -\frac{\ell_k}{N}\right\} 
 = \mathrm{coeff}_{\Gamma}(\{{B_{S}} - (L-kS)|_S\}). 
\]
Thus Assertion (\ref{eq:inequality-on-coeffs}) holds and the proof of the claim is concluded. \qedhere

\end{proof}

Now pick a section $\gamma \in H^0(X,\cO_X(L))$. 
By (\ref{eq:second-induction}), we have $\gamma= \gamma_1 + \gamma_2$ 
for some 
$\gamma_1 \in H^0(X,\cO_X(L-kS))$ and 
$\gamma_2 \in qS^0_{m, \adj}(X, S+B;L)$. 

Let $L' \sim L$ be a Weil divisor such that $S \not \subseteq \Supp(L'-kS)$. Via the induced isomorphism $\cO_X(L-kS) \cong \cO_X(L'-kS)$, we denote the image of $\gamma_1$ by $\gamma'_1 \in H^0(X,\cO_X(L'-kS))$. Then 
\[
(\gamma'_1)|_S \in H^0(S, \cO_S(\rup{(L'-kS)|_S-B_{S}}))
\]
(cf.\ the lower {horizontal sequence} {in the diagram (\ref{eq:diagram-res-C})}). 
By (\ref{e-IOAlem-QFS-S}) and Assumption (3), 
there exists a section 
\[
\sigma' \in qS^0_{m,\adj}(X,S+{B}; L'-kS) \subseteq H^0(X,\cO_X(L'-kS)) 
\]
such that $\sigma'|_S = \gamma'_1|_S$. Hence, 
\[
\gamma_1' = (\gamma_1'-\sigma') + \sigma',
\]
where $\gamma_1'-\sigma' \in H^0(X,\cO_X(L'-(k+1)S))$ and 
$\sigma' \in qS^0_{m,\adj}(X,S+{B}; L')$ 
by (\ref{eq:inclusion_of_im_psi-new}).

 In particular, if $\sigma \in qS^0_{m,\adj}(X,S+{B}; L)$ denotes the section corresponding to $\sigma'$ via the isomorphism $\cO_X(L') \cong \cO_X(L)$, then
\[
\gamma_1 = (\gamma_1-\sigma) + \sigma
\]
for $\gamma_1-\sigma \in H^0(X,\cO_X(L-(k+1)S))$. This concludes the proof of the claim 
{(\ref{eq:l-IOA-antiample-k+1})}.
\end{proof}

\begin{theorem} \label{thm:inversion-of-adjunction-anti-ample}
In the situation of Setting \ref{setting:foundations-of-log-quasi-F-splitting}, suppose that $R$ is a local ring and $H^0(X, \MO_X)=R$. 
Let $(X,S+B)$ be a divisorially Cohen-Macaulay {$p$-compatible}  
plt pair, where 
$S$ is a normal prime divisor and 
$B$ is a $\Q$-divisor such that $\rdown{B}=0$. 
Let $L$ be a $\bQ$-Cartier Weil divisor {on $X$ for which} $A := L-(K_X+S+B)$ is ample. {Suppose the following.} 
\begin{enumerate}
\item $-S$ is $\Q$-Cartier and semiample. 
\item $(S,B_S)$ is $n$-quasi-$F$-split,  where $K_S + B_S = (K_X+S+B)|_S$. 
\item $H^1(X, \cO_X(K_X + \lceil p^iA - kS \rceil))=0$
for all integers $i \geq 1$ and $k\geq 0$. 
\item $(X, B)$ is naively keenly $F$-pure. 
\end{enumerate}
Then, for {every} $m \gg n$, {we have that}
\[
 qS^0_{m, \adj}(X,S+B;L) = H^0(X,\cO_X(L)).
\]
\end{theorem}

\begin{proof}
Since $\lfloor B \rfloor =0$, we get that $L = K_X + S + \lceil A\rceil$. 
As $S \not\subseteq \Supp B = \Supp \{ A\}$, 
we may assume that $S \not\subseteq \Supp A$. 

Fix $m \gg 0$ {and}
set $r := m-n$.
By descending induction on $i$, we will show that the map
\begin{equation} \label{eq:main-induction}
qS^0_{m-i, \adj}(X,S+B_i; L_{i}-kS) \to qS^0_{m-i}(S,B_{S,i}; (L_{i}-kS)|_S)
\end{equation}
from Definition \ref{definitin:restriction-for-qS^0} is surjective\footnote{Formally, (\ref{eq:main-induction}) should be read as follows: for every $k \geq 0$ and $i \in \{0, 1, ..., r\}$ and for all (equivalently, for one) Weil divisors $D_{i,k} \sim L_i-kS$ such that $S \not \subseteq \Supp D_{i,k}$, we have that $qS^0_{m-i, \adj}(X,S+B_i; D_{i,k}) \to qS^0_{m-i}(S,B_{S,i}; D_{i,k}|_S)$ is surjective.} 
for all {$k \geq 0$ and $i \in \{0, 1, ..., r\}$}. 
Here $B_i := \{p^iB\}$, $L_{i} := K_X + S + \rup{p^iA}$, and $(K_X+S+B_i)|_S = K_S + B_{S,i}$. 
 Set $A_{i,k} := p^iA-kS$ and note that 
 \[
 L_{i}-kS = {K_X + S + \rup{A_{i,k}} = K_X + S + B_{i} + A_{i,k}}. 
\]
The reader should be aware that $B_i$ is not necessarily $p$-compatible. Plugging $i=0$ {into} (\ref{eq:main-induction}) yields the theorem by Lemma \ref{lem:auxiliary-inversion-of-adjunction-anti-ample}. \\
\vspace{0.0em}

\noindent \textbf{Base of induction.} We show that (\ref{eq:main-induction})   is surjective  for $i=r$. 
To this end, by the diagram {(\ref{eq:diagram-res-C})}  in Definition \ref{definitin:restriction-for-qS^0}, 
it is enough to show that 
\[
H^1(X, Q^{\vee}_{X,S+B_r-(L_{r}-kS),n}) = 
{H^1(X, \cHom_{\MO_X}(Q_{X, -K_X+B_r,n}, \MO_X(\rup{A_{r, k}}))) }
\]
is equal to zero. Since there are only finitely many possibilities 
for $B_r$, this holds by Fujita's vanishing, {see \cite[Theorem 1.5]{keeler03} and \cite{keeler03errata}} (we emphasise that $-S$ is nef and $A$ is ample).\\

\noindent \textbf{Inductive step.} 
Fix an integer $i$ with $0 \leq i < r = m-n$. 
We assume that
\[
qS^0_{m-i', \adj}(X,S+B_{i'}; L_{i'} - kS) \to qS^0_{m-i'}(S,B_{S,i'}; (L_{i'}-kS)|_S)
\]
is surjective for all $i' > i$ and $k\geq 0$,  and aim for showing (\ref{eq:main-induction}). {\cora Again, this should hold up to replacing $L_{i'} - kS$ by a Weil divisor $D_{i',k} \sim L_{i'} - kS$ such that $S \not \subseteq \Supp D_{i',k}$.}

By Theorem \ref{thm:lifting-sections-basic},  it is enough to prove that
\begin{align*}
&H^1(X, \cO_X(K_X + \lceil p^jA_{i,k} \rceil))=0 \text{ and }\\
&qS^0_{m-i-j}(X,\{p^jB_i\}; K_X + \lceil p^jA_{i,k}\rceil) = H^0(X,\cO_X(K_X + \lceil p^jA_{i,k}\rceil))
\end{align*}
for all $j\geq 1$. The former condition holds by Assumption (3). 
The latter condition, in view of  $\{p^jB_i\} = \{p^j\{p^iB\}\} = \{p^{i+j}B\} = B_{i+j}$, can be rewritten as
\begin{equation} \label{eq:key-thing-to-prove}
qS^0_{m-i'}(X,B_{i'}; L_{i'}-k'S) = H^0(X,\cO_X(L_{i'}-k'S))
\end{equation}
for $i' := i+j$ and $k' := {kp^j}+1$. 
By the induction hypothesis and Lemma \ref{lem:auxiliary-inversion-of-adjunction-anti-ample},  we get
\begin{equation*} 
qS^0_{m-i', \adj}(X, S+B_{i'},  L_{i'}-k'S) = H^0(X,\cO_X(L_{i'}-k'S)). 
\end{equation*}

{Here the assumptions of Lemma \ref{lem:auxiliary-inversion-of-adjunction-anti-ample} are satisfied as $(X,S+B)$ is $p$-compatible and 
$\{p^lB_{i'}\} = \{p^{l+i'}B\}$ for every $l \in \bZ_{\geq 0}$}.
In view of Lemma \ref{l-qS^0-adj-nonadj}, 
this implies (\ref{eq:key-thing-to-prove}). \qedhere
\end{proof}

\begin{corollary} \label{cor:inversion-of-adjunction-anti-ample}
In the situation of Setting \ref{setting:foundations-of-log-quasi-F-splitting},  suppose that $R$ is a local ring and $H^0(X, \MO_X)=R$. 
Let $(X,S+B)$ be a divisorially Cohen-Macaulay {$p$-compatible} 
plt pair, where 
$S$ is a normal prime divisor and 
$B$ is a $\Q$-divisor such that $\rdown{B}=0$.
{\cora Suppose that} $-(K_X+S+B)$ is ample and the following hold. 
\begin{enumerate}
\item $-S$ is $\Q$-Cartier and semiample.  
\item $(S,B_S)$ is $n$-quasi-$F$-split,  where $K_S + B_S = (K_X+S+B)|_S$. 
\item  $H^1(X, \cO_X(K_X + \lceil -p^i(K_X+S+B) - kS \rceil))=0$
for all  $i \geq 1$ and $k\geq 0$. 
\item $(X, B)$ is naively keenly $F$-pure. 
\end{enumerate}
Then $(X, S+B)$ is purely $m$-quasi-$F$-split  for some $m \geq n$.
\end{corollary}

\begin{proof}
By applying Theorem \ref{thm:inversion-of-adjunction-anti-ample} for $L=0$, 
the induced map 
\[
H^0(X,(Q^S_{X,S+B,m})^{\vee}) \to  
H^0(X, \cO_X)
\]
is surjective for some $m \gg n$ (cf.\ Definition \ref{d-pureQFS-sections}). 
By Lemma \ref{lem:dual-criterion-for-pure-quasi-F-splitting}, 
$(X,S+B)$ is purely $m$-quasi-$F$-split. 
\end{proof}

%% file: section5.tex
\section{Higher log Cartier operator and criterion for quasi-$F$-splitting}

\subsection{Introduction} \label{ss:sketch-of-proof-of-higher-Cartier-criterion}

We start by explaining the idea of the proof of Theorem \ref{thm:intro-higher-Cartier-criterion-for-quasi-F-split}. 
{\cora 
Let $X$ be a $d$-dimensional smooth projective Fano variety over a perfect field of characteristic $p>0$. Suppose that
\begin{enumerate}
    \item  $H^{d-2}(X,\Omega^1_X \otimes \omega_X) = 0$ and
    \item $H^{d-2}(X,B_1\Omega^2_X \otimes \omega^{p^k}_X) = 0$ for every $k \geq 0$.
\end{enumerate}
}
\noindent 
By {Proposition \ref{prop:intro-definition-via-splitting}}, 
{\eqref{eq:intro-C-quotient-sequence},} and Lemma \ref{lem:intro-cohomological-quasi-F-spliteness}(1), it is enough to show that
\[
H^{d-1}(X, (F_*W_n\cO_X/W_n\cO_X) \otimes \omega_X) = 0.
\]
Note that Serre's isomorphism (\ref{eq:Serre-map}) stipulates that $F_*W_n\cO_X/W_n\cO_X \cong B_n\Omega^1_X$. {By Serre's vanishing, we can pick $n>0$ such that $H^{d-1}(X, \Omega^1_X \otimes \omega^{p^n}_X)=0$. Then}
\begin{align*}
H^{d-1}(X, B_n\Omega^1_X \otimes \omega_X) = 0 &
\xLeftarrow{\hspace{0.08em}\hphantom{lal} \text{ Assumption (1) and } \eqref{eq:higher-Cartier-operator} \hphantom{lal}\hspace{0.08em}}
H^{d-1}(X,Z_n\Omega^1_X \otimes \omega_X)=0    \\
&\xLeftarrow{\small \text{ Assumption (2) and Lemma } \ref{lemma:key-ses-for-Zn}}
H^{d-1}(X, Z_{n-1}\Omega^1_X \otimes \omega^p_X) =0 \\
&\xLeftarrow{\hphantom{l} \text{repeating the above line } n \text{ times} \hphantom{l}}
H^{d-1}(X, \Omega^1_X \otimes \omega^{p^n}_X)=0.
\end{align*}

Using Theorem \ref{thm:intro-higher-Cartier-criterion-for-quasi-F-split}, we can give an alternative proof of the following result, which is a special case of  Corollary \ref{cor:cubic-surface}.

\begin{corollary}
Let $X$ be a smooth del Pezzo surface over a perfect field of characteristic $p > 0$. Then $X$ is quasi-$F$-split.
\end{corollary}

\begin{proof}
By Corollary \ref{c-descent2} (or \cite[Corollary 2.18]{KTY22}), we may assume that $k$ is algebraically closed.
Now note that $H^0(X, \Omega^1_X \otimes \omega_X)=0$ by Akizuki-Nakano vanishing (\cite[Coroilaire 2.8]{deligneillusie}) as $X$ is $W(k)$-liftable (indeed $X$ is a blowup of $\mathbb{P}^2$ {\cora or $\mathbb{P}^1\times \mathbb{P}^1$}at a finite number of points).
Moreover, we have \[H^0(X,B_1\Omega_X^2 \otimes \omega^{p^k}_X) \subseteq H^0(X,F_* \cO_X((p^{k+1}+1)K_X))=0\] as $X$ is Fano. Thus the statement follows from Theorem \ref{thm:intro-higher-Cartier-criterion-for-quasi-F-split}.
\end{proof}

In this section, we develop the logarithmic variant of Theorem \ref{thm:intro-higher-Cartier-criterion-for-quasi-F-split}.

\subsection{Higher Cartier operator with $\Q$-boundary}

Until the end of this subsection, we work under the assumptions of Setting \ref{setting:higher-log-Cartier-operator}.

\begin{setting} \label{setting:higher-log-Cartier-operator}
We assume the following: 
\begin{enumerate}
\setcounter{enumi}{-1}
\item $k$ is a perfect field $k$ of characteristic $p>0$, 
    \item $X$ is a regular integral separated scheme essentially of finite type  over $k$, 
    \item $E$ is a simple normal crossing divisor on $X$, 
and 
\item $\Delta$ is a $\bQ$-divisor on $X$ satisfying $\Supp(\{\Delta\})\subset \Supp E$. 
\end{enumerate}
\end{setting}

Let $E = \sum_{j=1}^s E_j$ be the irreducible decomposition.
For a $\Q$-divisor $\Gamma$, we set $\Omega_X^i(\log E)(\Gamma) := \Omega_X^i(\log E) \otimes \MO_X(\Gamma)$. 

\begin{remark}\label{r-local-imm}
Recall that there exists a factorisation 
\[
\pi \colon X \xrightarrow{j} X' \xrightarrow{\pi'} \Spec k
\]
where $\pi$ denotes the induced morphism, 
$j$ is a localising immersion, and $\pi'$ is separated and of finite type (see Subsection \ref{ss:notation}(\ref{it:essential})). 
After replacing $X'$ by an open subset of itself,  
we may assume that 
\begin{enumerate}
    \item[(1)'] $X'$ is a smooth variety over $k$, 
    \item[(2)'] there is a simple normal crossing divisor $E'$ on $X'$ with $E'|_X = E$, and 
\item[(3)'] there {is} a $\Q$-divisor $\Delta'$ on $X'$ such that $\Supp(\{\Delta'\})\subset\Supp E'$ and $\Delta'|_X = \Delta.$
\end{enumerate}
Here, we may assume (1)' and (2)' as open neighbourhoods of regular points and simple normal crossing divisors are regular and simple normal crossing, respectively (cf.\  \cite[\href{https://stacks.math.columbia.edu/tag/00OE}{Tag 00OE}]{stacks-project} and  \cite[\href{https://stacks.math.columbia.edu/tag/0BIA}{Tag 0BIA}]{stacks-project}).
{By} comparing $X$ and $X'$, {we can reduce many problems about $X$} to the case when $\pi \colon X \to \Spec k$ is of finite type. 
\end{remark}

As in Subsection \ref{ss:iterated-Cartier-operator}, one can define a logarithmic Cartier operator\footnote{{The isomorphism (\ref{e-logE-Cartier-op}) 
holds if $X \to \Spec k$ is of finite type. 
For the general case, we {invoke} Remark \ref{r-local-imm}. {Specifically,
we apply $j^*$} to the isomorphism  
$C^{-1} \colon \Omega^i_{X'}(\log E') \xrightarrow{\cong} \cH^i(F_*\Omega^{\mydot}_{X'}(\log E'))$,  
and use {that} $j^*$ is exact, $F_*j^* = j^*F_*$, {and} $j^*\Omega_{X'}^i(\log E') = \Omega_X^i(\log E)$. 
}} yielding an isomorphism of $\MO_X$-modules
\begin{equation}\label{e-logE-Cartier-op}
    C^{-1} \colon \Omega^i_X(\log E) \xrightarrow{\cong} \cH^i(F_*\Omega^{\mydot}_X(\log E)).
\end{equation}

{\cora The following observation (\cite{deligneillusie}, \cite{hara98a}) allows for extending this isomorphism to study the properties of log pairs: the inclusion 
\begin{equation}\label{e-Hara-obs}
F_*\Omega^{\mydot}_X(\log E) \hookrightarrow F_*(\Omega^{\mydot}_X(\log E) (D))
\end{equation}
is a quasi-isomorphism 
for every Weil divisor $D= \sum_{j=1}^s r_j E_j$ with $0 \leq r_j < p$, 
where $\Omega^{\mydot}_X(\log E) (D) := \Omega^{\mydot}_X(\log E) \otimes \MO_X(D)$ (see \cite[Section 4.2.1]{deligneillusie} or \cite[Lemma 3.3]{hara98a}}). 

\begin{lemma} \label{lemma:Hara}
{In the situation of Setting \ref{setting:higher-log-Cartier-operator}}, the inclusion
\[
(F_*\Omega^{\mydot}_X(\log E)) \otimes \cO_X(\Delta) \hookrightarrow F_*(\Omega^{\mydot}_X(\log E)( p\Delta))
\]
is a quasi-isomorphism.
\end{lemma}
{We remind the reader that $\Delta$ is a (possibly non-effective) $\Q$-divisor}.
\begin{proof}
Apply Hara's observation (\ref{e-Hara-obs})
for $D = \lfloor p\{\Delta\} \rfloor$ and tensor both sides by $\cO_X(\lfloor \Delta \rfloor)$. Here we use that $\rdown{ p\{\Delta\}} + p\rdown{\Delta} = \rdown{p\Delta}$.
\end{proof}

In what follows, we will study the complex
$F_*\Omega^{\mydot}_X(\log E)({p\Delta})$. 
{\cora Here, $F_*\Omega^{\mydot}_X(\log E)({p\Delta})$ should be read as 
$F_*(\Omega^{\mydot}_X(\log E)({p\Delta}))$ (see Remark \ref{remark:vigilant-about-parentheses}).}
We shall denote its cycles and boundaries by $B_1\Omega^i_{X}(\log E)(p\Delta)$ and $Z_1\Omega^i_{X}(\log E)(p\Delta)$, respectively: 
\begin{align}
\label{eq:definition-of-log-B1-and-Z1} B_1\Omega^i_{X}(\log E)(p\Delta) &:= {\rm Im}(d\colon 
F_*\Omega^{i-1}_X(\log E)({p\Delta})
\to 
F_*\Omega^{i}_X(\log E)({p\Delta})
) \\
Z_1\Omega^i_{X}(\log E)(p\Delta) &:= {\rm Ker}(d\colon  
F_*\Omega^{i}_X(\log E)({p\Delta})
\to 
F_*\Omega^{i+1}_X(\log E)({p\Delta}) \nonumber
).
\end{align}
In particular, $B_1\Omega^i_{X}(\log E)(p\Delta)$ and $Z_1\Omega^i_{X}(\log E)(p\Delta)$ are coherent $\MO_X$-modules.\\

Combined with the Cartier isomorphism (\ref{e-logE-Cartier-op}) {\cora and} Lemma \ref{lemma:Hara} yields an isomorphism 
\[
C_{\Delta}^{-1} \colon \Omega^i_X(\log E)({\Delta}) \xrightarrow{\cong} \cH^i(F_*\Omega^{\mydot}_X(\log E)({p\Delta})), 
\]
which leads to the following short exact sequence {{of} $\MO_X$-module homomorphisms}
\begin{equation} 
0\to B_1\Omega^{i}_X(\log E)(p\Delta)\to Z_1\Omega^{i}_X(\log E)(p\Delta) \xrightarrow{C_{\Delta}} \Omega^{i}_X(\log E)(\Delta) \to 0.\label{eq:ses-log-Cartier-operator}
\end{equation}
Define an $\MO_X$-module homomorphism
\[
C_{\Delta, n} := F^n_*C_{p^{n}\Delta} \colon F^n_*Z_1\Omega^i_X(\log E)(p^{n+1}\Delta) \to F^n_*\Omega^i_X(\log E)(p^{n}\Delta).
\]

\begin{remark}
In the spirit of Remark \ref{remark:iterated-partial-map}, one can think of $C_{\Delta}$ as a partial map \[
F_*\Omega^i_X(\log E)(p\Delta) \rightharpoonup \Omega^i_X(\log E)(\Delta),
\]
{\cora i.e., this is defined not on all of the domain $F_*\Omega^i_X(\log E)(p\Delta)$ but on an $\MO_X$-submodule of it.} 
The higher Cartier operator $C^n_{\Delta}$ will be the composition
\[
F^n_*\Omega^i_X(\log E)(p^n\Delta) \xrightharpoonup{C_{\Delta,n-1}} F^{n-1}_*\Omega^i_X(\log E)(p^{n-1}\Delta) \xrightharpoonup{C_{\Delta,n-2}} \ldots \xrightharpoonup{C_{\Delta}} \Omega^i_X(\log E)(\Delta),
\]
{the sheaf} $Z_n\Omega^i_X(\log E)(p^n\Delta)$ will be the domain of definition of $C^n_{\Delta}$, 
and {the sheaf} $B_n\Omega^i_X(\log E)(p^n\Delta)$ will be the kernel of $C^n_{\Delta}$.
\end{remark}

\begin{definition}
We inductively define
$Z_n\Omega_X^i(\log E)(p^n\Delta)$ and $B_n\Omega_X^i(\log E)(p^n\Delta)$, {which are} coherent $\MO_X$-modules, 
{as follows}
\begin{align*}
Z_n\Omega_X^i(\log E)(p^n\Delta) &:= (C_{\Delta,n-1})^{-1}(Z_{n-1}\Omega_X^i(\log E)(p^{n-1}\Delta)) \subseteq F^{n}_*\Omega^i_X(\log E)(p^n\Delta),\\
B_n\Omega_X^i(\log E)(p^n\Delta) &:= (C_{\Delta,n-1})^{-1}(B_{n-1}\Omega_X^i(\log E)(p^{n-1}\Delta)) \subseteq F^{n}_*\Omega^i_X(\log E)(p^n\Delta),
\end{align*}
where $Z_0\Omega_X^i(\log E)(\Delta):=\Omega_X^i(\log E)(\Delta)$ and $B_0\Omega_X^i(\log E)(\Delta):=0$.
\end{definition}

In particular, we get surjective {$\MO_X$-module homomorphisms} 
\begin{align}\label{e-Cartier-Zn-Zn-1}
C_{\Delta,n-1}|_{Z_n\Omega_X^i(\log E)(p^n\Delta)} \colon &Z_n\Omega_X^i(\log E)(p^n\Delta) \to Z_{n-1}\Omega_X^i(\log E)(p^{n-1}\Delta) \textrm{ {and} } \\
\label{e-Cartier-Bn-Bn-1} C_{\Delta,n-1}|_{B_n\Omega_X^i(\log E)(p^n\Delta)} \colon &B_n\Omega_X^i(\log E)(p^n\Delta) \to B_{n-1}\Omega_X^i(\log E)(p^{n-1}\Delta)
\end{align}
featured inside the following diagram (\ref{eq:Bn-1_Zn-1}):
\usetagform{hidden}
{\small\begin{equation} \label{eq:Bn-1_Zn-1}
\begin{tikzcd}[column sep = small]
\!0 \!\ar{r} & \!\! F^{n-1}_{*}B_1 \Omega^i_X(\log E)(p^n\Delta)  \arrow[d,dash,shift left=.1em] \arrow[d,dash,shift right=.1em] \!\ar{r} & \!\!  F^{n-1}_*Z_1\Omega_X^i(\log E)(p^n\Delta) \!\ar{r}{(\star)} & \! F^{n-1}_* \Omega_X^i(\log E)(p^{n-1}\Delta) \!\ar{r} & \!  0\\
\! 0 \!\ar{r} & \!\! F^{n-1}_{*}B_1 \Omega^i_X(\log E)(p^n\Delta) \arrow[d,dash,shift left=.1em] \arrow[d,dash,shift right=.1em] \!\ar{r} & \! Z_n\Omega_X^i(\log E)(p^n\Delta) \arrow[ru, phantom, "\usebox\pullbackdl" , very near start, yshift=0.2em, xshift=-2.7em, color=black] \arrow[hook]{u} \!\ar{r} & \!  Z_{n-1}\Omega^i_X(\log E)(p^{n-1}\Delta) \arrow[hook]{u} \!\ar{r} & \!  0 \\
\!0 \!\ar{r} & \!\! F^{n-1}_{*}B_1 \Omega^i_X(\log E)(p^n\Delta) \!\ar{r} & \! B_n\Omega_X^i(\log E)(p^n\Delta) \arrow[ru, phantom, "\usebox\pullbackdl" , very near start, yshift=0.2em, xshift=-2.7em, color=black] \arrow[hook]{u} \!\ar{r} & \!  B_{n-1}\Omega^i_X(\log E)(p^{n-1}\Delta) \arrow[hook]{u} \!\ar{r} & \! 0,
\end{tikzcd}
\end{equation}}
\!\!{\noindent where $(\star)  = C_{\Delta,n-1}$. Here each horizontal sequence is exact: the top one by (\ref{eq:ses-log-Cartier-operator}) and the other two by definition of a pullback.}
\usetagform{default}

\begin{definition} We define the \emph{higher Cartier operator} 
\[
C^{n}_{\Delta} \colon Z_n \Omega^i_X(\log E)(p^n\Delta) \to \Omega^i_X{(\log E)}(\Delta)
\]
as the composition of the above  {$\MO_X$-module homomorphisms (\ref{e-Cartier-Zn-Zn-1})} 
\[
Z_n\Omega_X^i(\log E)(p^n\Delta) \to Z_{n-1}\Omega_X^i(\log E)(p^{n-1}\Delta) \to \cdots \to \Omega_X^i(\log E)(\Delta).
\]
Specifically,  $C^{n}_{\Delta} := (C_{\Delta,0}|_{Z_1\Omega_X^i(\log E)(p\Delta)})\circ \cdots \circ (C_{\Delta,n-1}|_{Z_n\Omega_X^i(\log E)(p^n\Delta)})$.\\
\end{definition}

As in Subsection \ref{ss:iterated-Cartier-operator}, this name is justified by the short exact sequence
\begin{equation}
\label{eq:ses-log-higher-Cartier-operator}
0 \to B_n \Omega^i_X(\log E)(p^{n}\Delta) \to Z_n \Omega^i_X(\log E)(p^{n}\Delta) \xrightarrow{{C^n_{\Delta}}} \Omega^i_X(\log E)(\Delta) 
\to 0.
\end{equation}

Moreover, we get the following.
\begin{lemma} \label{lemma:key-ses-for-log-Zn}
{In the situation of Setting \ref{setting:higher-log-Cartier-operator}}, {there exists an exact sequence}
\[
0 \to Z_n\Omega^i_X(\log E)(p^{n}\Delta) \to F_*Z_{n-1}\Omega^i_X(\log E)(p^{n}\Delta) \to B_1\Omega^{i+1}_X(\log E)(p\Delta) \to 0
\]
of $\MO_X$-module homomorphisms.
\end{lemma}
\begin{proof}
The construction of the above sequence is analogous to that in Lemma \ref{lemma:key-ses-for-Zn}: take the pull-back of the sequence {(see (\ref{eq:definition-of-log-B1-and-Z1}))}
\begin{equation}
 0 \to Z_1\Omega^i_X(\log E)(p\Delta) \to F_*\Omega^i_X(\log E)(p\Delta) \xrightarrow{F_*d} B_1\Omega^{i+1}_X(\log E)(p\Delta) \to 0, \label{exaxt sequence from differential}   
\end{equation}
({note that} this is the $n=1$ case of the statement of the lemma) 
by 
\[
F_*C^{n-1}_{p\Delta} \colon F_*Z_{n-1}\Omega^i_X(\log E)(p^{n}\Delta) \to F_*\Omega^i_X(\log E)(p\Delta). \qedhere
\]
\end{proof}

We state the log variant of Serre's isomorphism {(\ref{eq:Serre-map})}.
\begin{lemma}\label{lem:log-Serre's map}
{In the situation of Setting \ref{setting:higher-log-Cartier-operator}}, 
\begin{equation}\label{e-log-Serre}
0\to W_n\cO_X(\Delta) \xrightarrow{F} F_{*}W_n\cO_X(p\Delta)\xrightarrow{s} B_n\Omega_X^{1}(\log E)(p^n\Delta)\to 0
\end{equation}
is an exact sequence  of $W_n\MO_X$-module homomorphisms, where 
\[
s(F_*(f_0,\ldots,f_{n-1}))\coloneqq 
F_*^n(f_{0}^{p^{n-1}-1}df_{0} + f_1^{p^{n-2}-1}df_1+\cdots+df_{n-1}).
\]
\end{lemma}

\noindent In particular, $B_{X,\Delta,n} \cong B_n\Omega_X^{1}(\log E)(p^n\Delta)$ (cf.\ (\ref{eq:definition-of-log-B})).
\begin{proof}
By Remark \ref{r-local-imm}, we may assume that $X$ is a smooth variety over $k$. We first check that $s$ is a $W_n\MO_X$-module homomorphism. 
It is easy to see that $s$ is an additive homomorphism.
Pick an open subset $U$ of $X$ and an element $F_*(f_0,\ldots,f_{n-1}) \in \Gamma(U, F_{*}W_n\cO_X(p\Delta))$. 
Then $s$ is a $W_n\MO_X$-module homomorphism, because 
the following equations hold for $\zeta \in {\MO_X(U)}$:
\begin{align*}
s( V[\zeta] \cdot F_*(f_0,\ldots,f_{n-1}) ) 
&= s( F_*( FV[\zeta] \cdot (f_0,\ldots,f_{n-1}) ){ )} \\
&=s( F_*( p \cdot [\zeta] \cdot (f_0,\ldots,f_{n-1})){ )} \\
&\overset{(\star)}{=}p \cdot s( F_*( [\zeta] \cdot (f_0,\ldots,f_{n-1})){ )} =0\\
     s( [\zeta] \cdot F_*(f_0,\ldots,f_{n-1}) ) 
 &= s( F_*( F[\zeta] \cdot (f_0,\ldots,f_{n-1}) )) \\
 &=
  s( F_*(\zeta^pf_0, \zeta^{p^2}f_1, \ldots, \zeta^{p^n}f_{n-1}) ) \\
 &= F_*^n( (\zeta^pf_{0})^{p^{n-1}-1}d(\zeta^pf_{0}) + (\zeta^{p^2}f_1)^{p^{n-2}-1}d(\zeta^{p^2}f_1)+\cdots+d(\zeta^{p^n}f_{n-1}))\\
&=[\zeta] \cdot F_*^n( f_{0}^{p^{n-1}-1}df_{0} + f_1^{p^{n-2}-1}df_1+\cdots+df_{n-1}), 
\end{align*}
where $[\zeta] \in W_n\cO_X(U)$ denotes the Teichm\"uler lift of $\zeta$ and $(\star)$ follows from the fact that $s$ is additive. Therefore, the sequence (\ref{e-log-Serre}) consists of $W_n\MO_X$-module homomorphisms. 

Let us prove that (\ref{e-log-Serre}) is exact at $F_*W_n\MO_X(p\Delta)$, 
{that is}, ${\rm Im}(F) = {\rm Ker} (s)$. 
Since the problem is local, we may assume that $X$ is affine. 
{By} taking the tensor product with the invertible $W_n\MO_X$-module $W_n\MO_X(-H)$ for an effective Cartier divisor $H$ {satisfying} $\Delta -H \leq 0$, we may assume that $\Delta \leq 0$. 
We then get the following commutative diagram of $W_n\MO_X$-module homomorphisms in which all the vertical arrows are the natural inclusions: 
\[\begin{tikzcd}
	0 \ar{r} & W_n\cO_X \ar{r} & F_{*}W_n\cO_X \ar{r} & B_n\Omega_X^{1} \ar{r} & 0 \\
	0 \ar{r} & W_n\cO_X(\Delta) \ar[hook]{u} \ar{r}{F}  & F_{*}W_n\cO_X(p\Delta) \ar[hook]{u} \ar{r}{s} & B_n\Omega_X^{1}(\log E)(p^n\Delta) \ar{r} \ar[hook]{u}  & 0. 
\end{tikzcd}\]
The upper horizontal sequence is exact by 
\cite[Lemme 2 in Section 7]{serre58} or 
\cite[Ch.\ I, Remarques 3.12(a)]{illusie_de_rham_witt}. 
Inside $F_*W_n\MO_X$, we have $W_n\MO_X \cap F_*W_n\MO_X(p\Delta) = W_n\MO_X(\Delta)$, 
which in turn implies ${\rm Im}(F) = \Ker(s)$.

Thus, it is enough to show that $s \colon F_{*}W_n\cO_X(p\Delta) \to  B_n\Omega_X^{1}(\log E)(p^n\Delta)$ is surjective. 
By inducting on $n$ and
using the exact sequence (see (\ref{eq:Bn-1_Zn-1})) {\small\[
0 \to F_*^{n-1}B_1\Omega_X^1(\log E)(p^{n}\Delta) \xrightarrow{v} B_{n}\Omega_X^1(\log E)(p^{n}\Delta) \xrightarrow{{C_{\Delta,n-1}}}
B_{n-1}\Omega_X^1(\log E)(p^{n-1}\Delta) \to  0,
\]}
\!\!
any local section of $B_n\Omega_X^{1}(\log E)(p^n\Delta)$ can be written as 
\begin{equation} \label{eq:desciption-of-B}
f_{0}^{p^{n-1}-1}df_{0}+ f_1^{p^{n-2}-1}df_1 + \cdots + df_{n-1} \textrm{ {for} some } f_{i}\in \mathcal{O}_X(p^{i+1}\Delta).
\end{equation}
To verify this, we first note that $v$ is a natural inclusion within $F^n_*\Omega_X$, and so 
\[
\mathrm{Im}(v) = \big\{df_{n-1} \mid f_{n-1} \in \cO_X(p^n\Delta)\big\} \subseteq F^n_*\Omega^1_X.
\]
Moreover, since $C_{\Delta,n-1}(f^{p-1}df) = df$, we have that
\begin{alignat*}{2}
\big\{C_{\Delta,n-1}&(f_{0}^{p^{n-1}-1}df_{0}+ f_1^{p^{n-2}-1}df_1 + \ldots + f_{n-2}^{p-1}df_{n-2}) &&\mid f_{i}\in \mathcal{O}_X(p^{i+1}\Delta)\textrm{ for } 0 \,{\leq}\, i \,{\leq}\, n{-}2 \big\} \\ 
=\ &\big\{f_{0}^{p^{n-2}-1}df_{0}+ f_1^{p^{n-3}-1}df_1 + \ldots + df_{n-2} &&\mid f_{i}\in \mathcal{O}_X(p^{i+1}\Delta)\textrm{ for } 0 \,{\leq}\, i \,{\leq}\, n{-}2 \big\} \\
=\ &B_{n-1}\Omega_X^1(\log E)(p^{n-1}\Delta),&& 
\end{alignat*}
where the last equality follows by induction. This concludes the verification of (\ref{eq:desciption-of-B}).
Therefore, $F_{*}W_n\cO_X(p\Delta)\xrightarrow{s} B_n\Omega_X^{1}(\log E)(p^n\Delta)$ is surjective. 
\qedhere
\end{proof}

\begin{lemma}
In the situation of Setting \ref{setting:higher-log-Cartier-operator}, we have that $B_n \Omega^i_X(\log E)(p^{n}\Delta)$ and $Z_n \Omega^i_X(\log E)(p^{n}\Delta)$ are locally free $\MO_X$-modules for all $i\geq 0$ and $n\geq  0$. 
\end{lemma}
\begin{proof}
We note that if the last term of a short exact sequence is locally free, then the sequence splits locally, and if, in addition, the first term or the middle term is locally free, then the 
other one is also locally free.
By {(\ref{eq:Bn-1_Zn-1})}, the following  sequences
{\small\begin{align*}
&0\to F^{n-1}_{*}B_1 \Omega^i_X(\log E)(p^{n}\Delta)\to B_n \Omega^i_X(\log E)(p^{n}\Delta)\xrightarrow{C_{\Delta,n-1}} B_{n-1} \Omega^i_X(\log E)(p^{n-1}\Delta)\to 0,\\
&0\to F^{n-1}_{*}B_1 \Omega^i_X(\log E)(p^{n}\Delta)\to Z_n \Omega^i_X(\log E)(p^{n}\Delta)\xrightarrow{C_{\Delta,n-1}} Z_{n-1} \Omega^i_X(\log E)(p^{n-1}\Delta)\to 0,
\end{align*}}
\!{are exact,} and thus we may assume that $n=1$. 
First, $Z_1\Omega^d_X(\log E)(p\Delta)=F_{*}\omega_X(E+p{\Delta})$ is locally free.
If $Z_1\Omega^{i}_X(\log E)(p\Delta)$ is locally free,
then so is $B_1\Omega^i_X(\log E)(p\Delta)$ by (\ref{eq:ses-log-Cartier-operator}), and if $B_1\Omega^{i}_X(\log E)(p\Delta)$ is locally free, then so is $Z_1\Omega^{i-1}_X(\log E)(p\Delta)$ by (\ref{exaxt sequence from differential}).
This proves the assertion of the lemma.
\end{proof}
 
\begin{remark} \label{remark:vigilant-about-parentheses}
We warn the reader to be {careful with} how the parentheses should be parsed. For example, $F_*\Omega^{i}_X(\log E)(p\Delta)$ should be read as $F_*(\Omega^{i}_X(\log E)(p\Delta))$. Similarly, $Z_n\Omega^{i}_X(\log E)(p^n\Delta)$ should be read as $Z_n(\Omega^{i}_X(\log E)(p^n\Delta))$. 

For this reason, we have {that}
\begin{align*}
B_n\Omega_X^{i}(\log E)(p^n\Delta)\otimes \cO_X(H)&=B_n\Omega_X^{i}(\log E)(p^n\Delta+p^nH), \text{ and }\\ Z_n\Omega_X^{i}(\log E)(p^n\Delta)\otimes \cO_X(H)&=Z_n\Omega_X^{i}(\log E)(p^n\Delta+p^nH)
\end{align*}
for a Cartier divisor $H$.
\end{remark}

\begin{remark}
The results of this subsection should hold for all regular integral Noetherian $F$-finite $\F_p$-schemes $X$. Indeed, the Cartier isomorphism is still valid in this generality at least in the non-log case
(cf.\ \cite[Proposition 2.5]{Shi07}, \cite[the first two paragraphs of Section 2]{KM21}). 
The authors decided not to pursue such a generalisation to avoid dealing with multiple technical details.
\end{remark}

\subsection{Criterion for quasi-$F$-splitting}
Since we use the Cartier operator, we assume that the schemes we consider are essentially of finite type over a perfect field of characteristic $p>0$.

\begin{theorem} \label{thm:higher-Cartier-criterion-for-quasi-F-split}
Let $(R, \m)$ be a local domain essentially of finite type over 
a perfect field of characteristic $p>0$.  
Let $(X,{\Delta)}$ be a $d$-dimensional 
log pair which is projective over $R$,  where $\rdown{\Delta}=0$. Let $f \colon Y \to X$ be a log resolution of $(X,\Delta)$ and let 
$B_Y$ be a $\bQ$-divisor such that  $\rdown{B_Y}\leq 0$, $-(K_Y+B_Y)$ is ample, and $f_*B_Y = \Delta$. 
Set $E \coloneqq \mathrm{Supp}(B_Y)$. Suppose that
\begin{enumerate}
    \item  $H^{d-2}_{\m}(Y,\Omega_Y^1(\log E)(K_Y + B_Y)) = 0$ and
    \item $H^{d-2}_{\m}(Y,B_1\Omega^2_Y(\log E)(p^k(K_Y + B_Y))) = 0$ for every $k \geq 1$.
\end{enumerate}
Then $(X,\Delta)$ is quasi-$F$-split.
\end{theorem}

\noindent {We refer to Notation \ref{notation:global-local-cohomology} for the definition of $H^i_{\m}$.} When applying the above theorem, we will often take $(X,\Delta)$ to be klt and $B_Y$ {to be} a small perturbation of 
the log pullback of ${\Delta}$. Also, to avoid confusion, we emphasise that 
\[
B_1\Omega^2_Y(\log E)(p^k(K_Y + B_Y)) \cong  B_1(\Omega^2_Y(\log E)(p^kB_Y)) \otimes \omega^{p^{k-1}}_Y.
\]
\begin{proof}
Write $K_Y+\Delta_Y = f^*(K_X+\Delta)$. Since $-(K_Y+B_Y)$ is ample and $f_*B_Y=\Delta$, the negativity lemma implies that $K_Y+B_Y \geq K_Y+\Delta_Y$, and so $f_*\cO_Y(K_Y+B_Y) = \cO_X(K_X+\Delta)$. 
Moreover, we obtain a natural map $f^* \colon Q_{X,K_X+\Delta,n} \to f_*Q_{Y,K_Y+B_Y,n}$ featured inside the following commutative diagram:
\[
\begin{tikzcd}[column sep = huge]
H^d_{\m}(\cO_X(K_X+\Delta)) \ar{r}{\Phi_{X,K_X+\Delta,n}} \arrow[d,dash,shift left=.1em] \arrow[d,dash,shift right=.1em]  & H^d_{\m}(Q_{X,K_X+\Delta,n}) \ar{d}
 \ar{d}{f^*}\\
H^d_{\m}(\cO_Y(K_Y+B_Y)) \ar{r}{\Phi_{X,K_Y+B_Y,n}} & H^d_{\m}(Q_{Y,K_Y+B_Y,n}).
\end{tikzcd}
\]
Note that the injectivity of the upper horizontal arrow is equivalent to $(X,\Delta)$ being $n$-quasi-$F$-split  by Lemma \ref{lem:cohomological-criterion-for-log-quasi-F-splitting}. 
Thus, to prove the theorem, it is enough to show that the lower horizontal arrow
\begin{equation} \label{Cartier-criterion-local-cohomology-injectivity}
H^d_{\m}(Y, \mathcal{O}_Y(K_Y+B_Y)) 
\to H^d_{\m}(Q_{Y,K_Y+ B_Y,n}) \ \textrm{is  injective for some } n>0.
\end{equation}
To this end, 
we can combine {\eqref{eq:C-quotient-sequence}, \eqref{eq:definition-of-log-B},} 
 and Lemma \ref{lem:log-Serre's map}, to get the short exact sequence:
\[
0 \to \mathcal{O}_Y({B_Y}) \to Q_{Y,{B_Y},n} \to B_n\Omega^1_Y(\log E)(p^n{B_Y}) \to 0.
\]
Thus, to show {Claim (\ref{Cartier-criterion-local-cohomology-injectivity})}, 
it is enough to prove that 
\[
H^{d-1}_{\m}(Y,B_n\Omega^1_Y(\log E)(p^n {B_Y}) \otimes {\omega_Y}) = 0.
\]
Now we use (\ref{eq:ses-log-higher-Cartier-operator}):
\[
0 \to B_n\Omega^1_Y(\log E)(p^n {B_Y}) \to Z_n\Omega^1_Y(\log E)(p^n { B_Y}) \xrightarrow{{C^n_{{B_Y}}}}  \Omega^1_Y(\log E)({B_Y}) \to 0,
\]
By tensoring by {$\omega_Y$} again, the problem is reduced to showing  that 
\begin{itemize}
    \item $H^{d-1}_{\m}(Y,Z_n\Omega^1_Y(\log E)(p^n{B_Y}) \otimes { \omega_Y}) = 0$ and
    \item $H^{d-2}_{\m}(Y,\Omega^1_Y(\log E)({B_Y}) \otimes { \omega_Y})=0$.
\end{itemize}

The latter assertion is nothing but our Assumption (1).
Thus, it is sufficient to prove the former assertion. 

Now we use the short exact sequence:
{\small \[
0 \!\shortrightarrow\! 
F^k_*Z_{n{-}k} \Omega^1_Y(\log E)(p^n{B_Y}) \!\shortrightarrow\! 
F^{k{+}1}_*Z_{n{-}k{-}1} \Omega^1_Y(\log E)(p^n{B_Y}) \!\shortrightarrow\! 
F^k_*B_1\Omega^2_Y(\log E)(p^{k{+}1}{B_Y}) \!\shortrightarrow\! 0
\]}
\!\!from Lemma \ref{lemma:key-ses-for-log-Zn}. 
By tensoring by {$\omega_Y$} and repeatedly applying Assumption (2), 
{we get the following injections
\begin{align*}
H^{d-1}_{\m}(Y,\, Z_n \Omega^1_Y(\log E)(p^n{B_Y}) \otimes {\omega_Y}) 
\xhookrightarrow{\hphantom{aa}}\ 
&H^{d-1}_{\m}(Y, F_*Z_{n-1} \Omega^1_Y(\log E)(p^n{B_Y}) \otimes { \omega_Y}) \\
\xhookrightarrow{\hphantom{aa}}\  &\cdots \\
\xhookrightarrow{\hphantom{aa}}\
&H^{d-1}_{\m}(Y, F_*^n\Omega^1_Y(\log E)(p^n{B_Y}) \otimes { \omega_Y}). 
\end{align*}
}
Hence, it suffices to show that 
\[
H^{d-1}_{\m}(Y,F^{n}_*\Omega^1_Y(\log E)(p^n{B_Y}) \otimes {\omega_Y}) = H^{d-1}_{\m}(Y, \Omega^1_Y(\log E)(p^n(K_Y + {B_Y})))
\]
is zero for $n \gg 0$. 
By Matlis duality (Proposition \ref{prop:Matlis-duality}) and exactness of the derived $\mathfrak m$-completion in our setting \cite[\href{https://stacks.math.columbia.edu/tag/0A06}{Tag 0A06}]{stacks-project}, it is enough to verify that
\begin{align*}
H^{-d+1}&R\Hom(R\pi_*\Omega^1_Y(\log E)(p^n(K_Y + {B_Y})), \omega^{\mydot}_R) \\
= H^{1}&R\Hom(\Omega^1_Y(\log E)(p^n(K_Y + {B_Y})), \omega_Y) \\
= H^1&(Y, \Omega^1_Y(\log E)^{\vee} \otimes \cO_Y(K_Y-p^n(K_Y + {B_Y})))
\end{align*}
is equal to zero, where $\pi \colon Y \to \Spec R$ is the natural projection. This, in turn, is a consequence of Serre's vanishing (\cite[\href{https://stacks.math.columbia.edu/tag/0B5U}{Tag 0B5U}]{stacks-project} or \cite[Theorem III.8.8]{hartshorne77}).
\end{proof}

\begin{remark} \label{remark:log-criterion-no-B}
As indicated after Theorem \ref{thm:intro-higher-Cartier-criterion-for-quasi-F-split}, 
Assumption (2) of Theorem \ref{thm:higher-Cartier-criterion-for-quasi-F-split} is valid if the following holds: 
\[
H^{d-i-1}_{\m}(Y,\Omega^i_Y(\log E)(p^{ k-1}(K_Y+B_Y)))=H^{d-i}_{\m}(Y,\Omega^i_Y(\log E)(p^{ k}(K_Y+B_Y))) = 0
\]
for every $k \geq {1}$ and $i \geq 2$.

To prove this, we consider the short exact sequences {(see \eqref{eq:definition-of-log-B1-and-Z1} and \eqref{eq:ses-log-Cartier-operator}): 
\begin{align*}
&0 \to B_1\Omega^i_Y(\log E)(p\Delta) \to Z_1\Omega^i_Y(\log E)(p\Delta) \xrightarrow{C_{\Delta}} \Omega^i_Y(\log E)(\Delta) \to 0, \\
&0 \to Z_1\Omega^i_Y(\log E)(p\Delta) \to F_*\Omega^i_Y(\log E)(p\Delta) \xrightarrow{F_{*}d} B_1\Omega^{i+1}_Y(\log E)(p\Delta) \to 0
\end{align*}}
for $\Delta = p^{k-1}(K_Y+B_Y)$.
Then, for $k\geq 1$ and $i\geq 2$, we have the following implications:
\begin{align*}
H^{d-i}_{\m}(Y&, B_1\Omega^i_Y(\log E)(p^{k}(K_Y+B_Y)) ) = 0\\
&\xLeftarrow{H^{d-i-1}_{\m}(Y, \Omega^i_Y(\log E)(p^{k-1}(K_Y+B_Y)))=0} 
H^{d-i}_{\m}(Y, Z_1\Omega^i_Y(\log E)(p^{k}(K_Y+B_Y))) = 0\\ &\xLeftarrow{H^{d-i}_{\m}(Y, \Omega^i_Y(\log E)(p^{k}(K_Y+B_Y)) )=0}
H^{d-i-1}_{\m}(Y, B_1\Omega^{i+1}_Y(\log E)(p^{k}(K_Y+B_Y))) = 0.
\end{align*}
In particular, we get that $H^{d-2}(Y, B_1\Omega^i_Y(\log E)(p^{k}(K_Y+B_Y)))=0$ 
by starting with $i=2$ and repeating the above $d-1$ times.
Therefore, Assumption (2) of Theorem \ref{thm:higher-Cartier-criterion-for-quasi-F-split} is valid. 
\end{remark}

\begin{remark} \label{remark:log-criterion: splitting of resolution}
In the setting of Theorem \ref{thm:higher-Cartier-criterion-for-quasi-F-split}, if we assume that 
{$\lfloor B_Y\rfloor=0$ instead of $\lfloor B_Y\rfloor\leq 0$}, then \eqref{Cartier-criterion-local-cohomology-injectivity} and Lemma \ref{lem:cohomological-criterion-for-log-quasi-F-splitting} show that $(Y,B_Y)$ is quasi-$F$-split.
\end{remark}
 The above theorem implies that one-dimensional log Fano pairs are quasi-$F$-split.

\begin{corollary} 
Let $\kappa$ be a field {which is} finitely generated over a perfect field $k$ of characteristic $p>0$. 
Let $X$ be a regular projective curve over $\kappa$ and 
let $\Delta$ be an effective $\bQ$-divisor on $X$ such that $(X,\Delta)$ is klt and  $-(K_X+\Delta)$ is ample. Then $(X,\Delta)$ is quasi-$F$-split.
\end{corollary}
\begin{proof}
Note that $\kappa$ is essentially of finite type over $k$. 
Thus, we may apply Theorem \ref{thm:higher-Cartier-criterion-for-quasi-F-split} 
for ${\rm id} \colon Y =X \to X$.
\end{proof}

\begin{corollary}\label{c-QFS-2dim-rel}
Let $R$ be a ring essentially of finite type over a perfect field 
of characteristic $p>0$. 
Let $X$ be a two-dimensional integral normal scheme and 
let $\pi \colon X\to \Spec R$ be a projective morphism with $\dim \pi(X) \geq 1$. 
Let $\Delta$ be an effective $\bQ$-divisor on $X$ such that $(X,\Delta)$ is klt and $-(K_X+\Delta)$ is ample. Then $(X,\Delta)$ is quasi-$F$-split.
\end{corollary}

\noindent In particular, two-dimensional klt pairs 
$(X,\Delta)$ essentially
of finite type over a perfect field of characteristic $p>0$ are always quasi-$F$-{pure}.

\begin{proof} 
{By} taking the Stein factorisation of $\pi \colon X \to \Spec R$ and a localisation of $R$, 
we may assume that $H^0(X, \MO_X)=R$ and $R$ is a local domain. 
Let $f \colon Y \to X$ be a log resolution of $(X,\Delta)$ {sitting inside the following diagram} 
\[
g \colon Y \xrightarrow{f} X \xrightarrow{\pi} \Spec R. 
\]
Set $K_Y+\Delta_Y = f^*(K_X+\Delta)$. 
As $X$ is $\Q$-factorial, 
there exists an $f$-exceptional effective $\Q$-divisor $F$ such that $-F$ is $f$-ample. 
Set $B_Y := \Delta_Y + \epsilon F$ for some $0 < \epsilon \ll 1$. 
We may assume that $\rdown{B_Y} \leq 0$, $-(K_Y+B_Y)$ is ample, and $f_*B_Y = \Delta$. 

It is enough to verify Assumption (1) and {Assumption} (2) of Theorem  \ref{thm:higher-Cartier-criterion-for-quasi-F-split}. 
{Recall that}
\[
H^0 R\Gamma_{\m} Rg_* \cF  = H^0_{\m}(g_*\mathcal{F}) =0
\]
for a torsion-free coherent $\MO_Y$-module $\cF$, where the first equality holds by $R\Gamma_{\m} Rg_* = R(\Gamma_{\m} \circ g_*)$ and the second one follows from the fact that  $g_*\mathcal{F}$ is also torsion-free. 
Hence, Assumptions (1) and (2) 
of Theorem  \ref{thm:higher-Cartier-criterion-for-quasi-F-split}
hold automatically.
\qedhere
\end{proof}

%% file: section6.tex
\section{Log del Pezzo surfaces} \label{section-log-del-Pezzo}
 In this section, we study when  log del Pezzo pairs are quasi-$F$-split.
\begin{remark}
In what follows, $k$ is always a perfect field of characteristic $p>0$. We will often {take} base change  to the algebraic closure $\bar{k}$ of $k$. We point out that this operation preserves all natural classes of singularities (such as klt or lc), but it does \emph{not} preserve $\bQ$-factoriality 
(see \cite[Remark 2.7(1)]{GNT16} for details). 
\end{remark}

\begin{remark} \label{remark:properties-del-Pezzo-type}
Let $k$ be a perfect field of characteristic $p>0$ and let $X$ be a surface of del Pezzo type over $k$. Let $\Delta$ be a $\bQ$-divisor such that $(X,\Delta)$ is log del Pezzo. Then
\begin{enumerate}
    \item for every $\bQ$-divisor $D$, we can run a $D$-MMP which terminates with a minimal model or a Mori fibre space,
    \item every nef $\bQ$-divisor on $X$ is semiample, {and}
    \item if $f \colon X \to Y$ is a birational morphism and 
    $Y$ is a normal projective surface, then $Y$ is of del Pezzo type. Specifically, 
    $(Y, f_*\Delta)$ is log del Pezzo.
\end{enumerate}
Assertion (1) holds as we can run an MMP for 
\begin{equation} \label{eq:del-Pezzo}
\epsilon D =  K_X + \Delta - (K_X+\Delta) + \epsilon D
\end{equation}
and any $0 < \epsilon \ll 1$ by \cite[Theorem 1.1]{tanaka16_excellent}. For a nef $\bQ$-divisor $D$, Assertion (2) holds by the base point free theorem  \cite[Theorem 1.2]{tanaka12} applied to (\ref{eq:del-Pezzo}). 
As for Assertion (3), we refer to \cite[Lemma 2.8]{Tan20} (cf.\ \cite[Lemma 2.4]{Bir16}). 
\end{remark}

We start with the case when the boundary is empty.

\begin{theorem}\label{thm:qFs of klt dPs}
Let $X$ be a klt del Pezzo surface over {an algebraically closed field} of characteristic $p>0$.
Then the following are equivalent.
\begin{enumerate}
    \item There exists a log resolution $f\colon Y\to X$ such that $(Y,\Exc(f))$ lifts to $W(k)$.
    \item $X$ is quasi-$F$-split.
\end{enumerate}
In particular, $X$ is quasi-$F$-split if $p>5$.
\end{theorem}

\begin{proof}
We first prove that (1) {implies} (2). Suppose that there exists a log resolution $f\colon Y\to X$ such that $(Y,\Exc(f))$ lifts to $W(k)$. Set $\Delta_Y \coloneqq f^{*}K_X-K_Y$.
{In particular,}  ${\lfloor} \Delta_Y \rfloor \leq 0$. 
By increasing the coefficients of $\Delta_Y$, 
we can find an $f$-exceptional $\Q$-divisor $B_Y$ such that $B_Y \geq \Delta_Y$, 
$\rdown{B_Y} \leq 0$, and $-(K_Y+B_Y)$ is ample.
{We can conclude the proof of the theorem by applying Theorem} \ref{thm:higher-Cartier-criterion-for-quasi-F-split} 
provided that we verify its assumptions. 
{Assumption} (1) of Theorem \ref{thm:higher-Cartier-criterion-for-quasi-F-split} holds by \cite[Theorem 2.4]{Kawakami-Nagaoka}.
{Assumption} (2) of Theorem \ref{thm:higher-Cartier-criterion-for-quasi-F-split} is also satisfied, since we have {the following inclusions}
\begin{align*}
H^0(Y, B_1\Omega_Y^2(\log\,E)({\lfloor} p^k(K_Y+B_Y){\rfloor}))&\xhookrightarrow{\hphantom{aa}}
H^0(Y, \mathcal{O}_Y(K_Y+E+{\lfloor} p^k(K_Y+B_Y){\rfloor}))\\
&\xhookrightarrow{\hphantom{aa}} H^0(X, \mathcal{O}_X((1+p^k)K_X)) =0
\end{align*}
for all $k\in\mathbb{Z}_{>0}$, where $E := \Exc(f)$.

We next prove that (2) {implies} (1). 
Let $\pi\colon Z\to X$ be the minimal resolution {of $X$ and let $K_Z+\Delta_Z = \pi^* K_X$}.
Then $Z$ is quasi-$F$-split by 
Corollary \ref{corollary:log-pullback-quasi-F-split-dimension-two}.  In turn, Theorem \ref{thm:p^2-lift of pairs} and Corollary \ref{cor:p^2-lift of snc pairs} show that $(Z, \Exc(\pi))$ lifts to $W_2(k)$. Then $(Z, \Exc(\pi))$ lifts to $W(k)$ by \cite[Proposition 2.8 (1)]{Kawakami-Nagaoka}.

Finally, if $p>5$, 
then a log resolution as in (1) exists by \cite[Theorem 1.2]{ABL}. 
\end{proof}

\begin{corollary} \label{cor:del-Pezo-type-surfaces-are-quasi-F-split}
Let $Y$ be a
surface of del Pezzo type over a perfect field $k$ of characteristic $p>5$. Then $Y$ is quasi-$F$-split.
\end{corollary}

\begin{proof}
By Proposition \ref{prop:descent}, we may replace $k$ by its algebraic closure, and so, from now on, we can assume that $k$ is algebraically closed. Let $\Delta_Y$ be an effective $\bQ$-divisor such that $(Y,\Delta_Y)$ is a log del Pezzo pair. Let $f\colon Y \to X$ be the composition of 
\begin{enumerate}
    \item a $(-K_Y)$-MMP $g' \colon Y \to X'$ (see Remark \ref{remark:properties-del-Pezzo-type}(1)) and
\item the $(-K_{X'})$-semiample fibration $h \colon X' \to X$ (see Remark \ref{remark:properties-del-Pezzo-type}(2)).
\end{enumerate}
By construction,  $X$ is a klt del Pezzo surface (see Remark \ref{remark:properties-del-Pezzo-type}(3)) and {$f^{*}K_{X}-K_{Y}$} is effective.
Since $p>5$, the surface $X$ is quasi-$F$-split by Theorem \ref{thm:qFs of klt dPs}. 
Therefore, we can conclude that $Y$ is quasi-$F$-split by Corollary  \ref{corollary:log-pullback-quasi-F-split-dimension-two}. \qedhere
\end{proof}

\subsection{Log del Pezzo surfaces in large characteristic}

Our goal is to show that all log del Pezzo pairs $(X,\Delta)$ with standard coefficients are $2$-quasi-$F$-split if the characteristic $p$ is large enough (Theorem \ref{thm:qFsplit-del-Pezzo-in-large-characteristic}). 

\begin{remark} For the convenience of the reader, we provide an outline of the proof of Theorem \ref{thm:qFsplit-del-Pezzo-in-large-characteristic}. We start by fixing a positive integer $N \gg 0$ independent of $(X,\Delta)$ (see Setting \ref{setting:del-pezzo-large-characteristic} for the precise definition).  
Using the BAB conjecture, one can show that if $(X,\Delta)$ is $\frac{1}{N}$-klt, then it is globally $F$-regular provided that 
$p \gg N$ (see \cite[Remark 3.3]{CTW15b}).

Therefore, the key difficulty is to show that if $(X,\Delta)$ is \emph{not} $\frac{1}{N}$-klt, then it is quasi-$F$-split (this is exactly the content of Theorem \ref{thm:qFsplit-del-Pezzo}). To this end, we extract a prime divisor $C$ via $f \colon Z \to X$ {(possibly, $f=\mathrm{id}$)} such that $1-\frac{1}{N} < b < 1$ for $K_Z+bC + B = f^*(K_{X}+\Delta)$, where $C \not \subseteq \Supp B$ and $B$ is an effective $\bQ$-divisor. Note that $(Z,aC+B)$ is log del Pezzo if $a := b + \epsilon$ for $0 < \epsilon \ll 1$.

By definition of $N\gg 0$ and by running some MMPs, 
the problem is reduced to the case when $(Z,C+B)$ is log canonical and $-(K_Z+C+B)$ is semiample. Set $\kappa := \kappa(Z, -(K_Z+C+B))$. Then thanks to different variants of inversion of adjunction:
\begin{enumerate}
    \item $(Z,C+B)$ is sharply $F$-split if $\kappa= 2, \kappa = 1$, or both $\kappa =0$ and $(Z,C+B)$ is \emph{not} plt (see Lemma \ref{lem:F-split-inversion-of-adjunction-kappa-two}, Lemma \ref{lem:F-split-inversion-of-adjunction-kappa-one}, and \cite[Lemma 2.9]{CTW15b}, respectively),
    \item $(Z,C+B)$ is purely $2$-quasi-$F$-split if it is plt and $\kappa=0$.
\end{enumerate}
\end{remark}

Most of the {above} argument follows the idea of \cite{CTW15b} 
except for Case (2), in which $(Z,C+B)$ need not be $F$-split in general (cf.\  \cite{CTW15a}), but we can {now} show that it is always purely $2$-quasi-$F$-split.

For the convenience of the reader, we decided to rewrite the argument of \cite{CTW15b} referring to some of the lemmas therein but not the proofs themselves. In doing so, we also streamlined the argument and made it, in a way, more general as we now cover the case of standard coefficients, too.

We emphasise that, for us, the set of standard coefficients contains $1$.
\begin{setting} \label{setting:del-pezzo-large-characteristic}
We pick integers {$N_1,N_2>5$} which satisfy the following properties:
\begin{enumerate} 
    \item given a  two-dimensional {projective lc} pair $(X,\Delta)$ with standard coefficients over $k$ such that $X$ is klt and $K_X+\Delta \equiv 0$, we have that $N_1(K_X+\Delta)$ is Cartier, 
    \item  
    given a projective lc pair $(X, \Delta)$ over $k$ and a morphism $\pi \colon X \to T$ {such that}
    \begin{enumerate}
    \vspace{0.2em}
    \begin{minipage}{0.4\linewidth}
        \item $\dim X \leq 2$, 
        \item $K_X+\Delta \equiv_{\pi} 0$,
        \item $\pi_*\cO_X = \cO_T$,
    \end{minipage}
    \begin{minipage}{0.6\linewidth}
        \item {$\mathrm{coeff}(\Delta) \subseteq \{ 1 - \frac{1}{n} \mid n \geq 1\} \cup (1-\frac{1}{N_2},1]$},
        \item $T$ is a projective normal variety,
        \item {$\rho(X/T)=1$} 
    \end{minipage}
    \vspace{0.05em}
    \end{enumerate}
    we have that the coefficients of the $\pi$-horizontal divisors in $\Delta$ belong to the set $\{1-\frac{1}{n} \mid 1 \leq n \leq {N_2}\}\cup \{1\}$. Here $\mathrm{coeff}(\Delta)$ denotes the set consisting of the coefficients of $\Delta$.
\end{enumerate}
Finally, set $N:=N_1\cdot N_2$. {Here we say that a prime divisor $C$ on $X$ is {\em $\pi$-horizontal} if $\pi(C) = T$.} 
The existence of such {integers} 
{$N_1$ and $N_2$} follows from \cite[Proposition 2.5 and Lemma 2.3]{CTW15b}. 
As for {Case} {(2)},
one should expect ${N_2}=42$ (cf.\ Remark \ref{rem:main-theorem-for-42}), but we do not know any small explicit bound on $N$ {in} {Case} {(1)}.  Note that we assume that $N>5$ for convenience, so that the condition $p>N$ implies that $p>5$ for a prime number $p$.
\end{setting}

\begin{lemma} \label{lem:lc-is-sharply-F-pure}
Let $(X,\Delta)$ be a two-dimensional {lc} pair of finite type over a perfect field 
$k$ of characteristic $p>5$. Suppose that $(X, \Delta)$ is klt outside $\Supp \rdown{\Delta}$ and $\Delta$ has standard coefficients. Then $(X,\Delta)$ is sharply $F$-pure.
\end{lemma}

\begin{proof}
By taking the base change to the algebraic closure of $k$, 
the problem is reduced to the case when $k$ is algebraically closed. 
The question is local, so we {can} fix a closed point $x \in X$ around which we will work. 
We may assume that 
$X \setminus \{x\}$ and $\Supp\, (\Delta) \setminus \{x\}$ 
are smooth and that the irreducible components of $\Delta$ pass through~$x$. 

If $(X,\Delta)$ is klt at $x$, then it is strongly $F$-regular by \cite[Theorem 3.1]{hx13}. 
If $(X,\Delta)$ is plt and not klt at $x$, then
$(X,\Delta)$ is purely $F$-regular (hence sharply $F$-pure) by inversion of $F$-adjunction \cite[Theorem A]{Das15} 
and the fact that every klt pair on $\bP^1$ is strongly $F$-regular.

Finally, suppose that $(X,\Delta)$ is not plt at $x$. Let $\Gamma \subseteq \Supp \Delta$ be a {prime divisor} passing through $x$ {such that $\mathrm{coeff}_{\Gamma}(\Delta) = 1-\frac{1}{n}$ for $n \in \bZ_{>0}$}.
\begin{claim} We have that $n\leq 6$.
\end{claim}
\begin{proof}
Assume by contradiction that $n > 6$, and take a dlt blowup $f \colon Y \to X$ of $(X,\Delta)$ {(see \cite[Theorem 4.7]{tanaka16_excellent}) for its definition and the proof of existence). Specifically, we require} that $f^*(K_X+\Delta) = K_Y + E + f_*^{-1}\Delta$, where  $E$ denotes the reduced $f$-exceptional divisor.  
Pick an $f$-exceptional prime divisor $F \subseteq \Supp E$ intersecting $f_*^{-1} \Gamma$ 
and  write $K_F + \Delta_F = f^*(K_X+\Delta)|_F$. 

Write $\Delta_F = \sum a_i P_i$ for rational numbers $0<a_i\leq 1$ and distinct points $P_i$. Without loss of generality, we may assume that $P_1 \in F \cap f_*^{-1}\Gamma$. Since $n>2$, the pair $(Y, E + f_*^{-1}\Delta)$ is plt at $P_1$ by Lemma \ref{lem:surface-lc-pair-with-high-coeffs-is-plt}, and so $1 > a_1 \geq 1-\frac{1}{n}$ by Lemma \ref{lem:everything-about-plt-surface-pairs}(5) and (\ref{eq:coefficients-increase-in-adjunction}). In particular, $a_1 > \frac{5}{6}$ as $n>6$. Since $\Delta_F$ has standard coefficients and $\deg \Delta_F = 2$, this contradicts Remark \ref{remark:ACC-P^1}.
\qedhere
\end{proof}

Since $(X,\Delta)$ is klt outside $\Supp \rdown{\Delta}$ and $(X, \Delta)$ is lc, 
$(X, (1-\frac{1}{m})\rdown{\Delta}+ \{ \Delta\})$ is a klt pair with standard coefficients for every integer {$m>0$}, and so 
 it is strongly $F$-regular again by \cite[Theorem 3.1]{hx13}. 
Then the limit pair $(X,\Delta)$ is $F$-pure by \cite[Theorem 3.2]{CTW15a} (cf.\ \cite[Proposition 2.6]{hara06}, \cite[Remark 5.5]{schwede08}). 
Since $p>5$ and the coefficients of $\Delta$ {are contained in {the set} $\{\frac{1}{2}, \frac{2}{3}, \frac{3}{4}, \frac{4}{5}, \frac{5}{6}, 1\}$}
by the above claim, we get that $(X,\Delta)$ is sharply $F$-pure.
\end{proof}

\begin{lemma} \label{lem:F-split-inversion-of-adjunction-kappa-two}
Let $(Z,C+B)$ be a two-dimensional projective lc  pair over an algebraically closed field $k$ of characteristic $p>5$, where $C$ is a prime divisor and $B$ is an effective $\bQ$-divisor with standard coefficients. Suppose that
\begin{enumerate}
    \item $(Z,aC+B)$ is log del Pezzo for some $a\geq 0$ and
    \item $-(K_Z+C+B)$ is nef and big.
\end{enumerate}
\noindent Then $(Z,C+B)$ is sharply $F$-split.
\end{lemma}
\begin{proof}
We first reduce the problem to the case when $-(K_Z+C+B)$ is ample, as this is needed for the inversion of adjunction to work.
\begin{claim}
We may assume that 
\begin{enumerate}
\item[(2)']    $-(K_Z+C+B)$ is ample.
\end{enumerate}
\end{claim}

\begin{proof}
Since $(Z, aC+ B)$ is log del Pezzo, the nef $\bQ$-divisor $-(K_Z+C+B)$ is semiample {(see Remark \ref{remark:properties-del-Pezzo-type}(2))}. Thus, there exists a birational morphism $f \colon Z \to Z'$ 
to a projective normal surface $Z'$ 
such that $K_Z+C+B \sim_{\Q} f^*(K_{Z'} +C' + B')$ 
and 
\[
-(K_{Z'}+C'+B') \textrm{ is an ample } \Q\textrm{-Cartier } \Q\textrm{-divisor}
\]
for $C' := f_*C$ and $B' := f_*B$. 
In particular, $(Z', C'+B')$ is lc. 
Moreover, {the push-forward $C'$ is still a prime divisor, given that} the $\Q$-divisor $(1-a)C =  (K_Z+C+B) -(K_Z+ aC+B) $ is $f$-ample. 
As $(Z, aC+B)$ is log del Pezzo, also $(Z', aC'+B')$ is log del Pezzo {(see Remark \ref{remark:properties-del-Pezzo-type}(3))}. 

{Since we may assume} that {the lemma} is valid when {(2)'} holds, $(Z', C'+B')$ is sharply $F$-split. 
Therefore, $(Z, C+B)$ is sharply $F$-split by  \cite[Lemma 3.3]{GT16} in view of the 
{\cora $\Q$-linear equivalence} 
$K_Z + C + B \sim_{\bQ} f^*(K_{Z'} + C' + B')$. \qedhere
\end{proof}

From now on, we assume that {(2)' holds.} Since the Cartier index of $(Z+C+B)$ may be divisible by $p$, we might need to perturb $B$.

\begin{claim} There exists a $\bQ$-divisor $B' \geq B$ such that 
 \begin{enumerate}
     \item[(a)] $(Z,C+B')$ is log canonical, 
     \item[(b)] $-(K_Z+C+B')$ is ample, and 
     \item[(c)] $(p^e-1)(K_Z+C+B')$ is Cartier for some $e>0$.
 \end{enumerate}
\end{claim}

\begin{proof}
Note that the pair $(Z,C+B)$ is sharply $F$-pure by Lemma \ref{lem:lc-is-sharply-F-pure}, and so by \cite[Lemma 2.11]{CTW15b}, we get an effective $\bQ$-divisor $H$ such that
$(Z,C+B+H)$ is sharply $F$-pure (in particular, log canonical) {and $(p^e-1)(K_Z+C+B+H)$ is Cartier for some $e>0$. Then (a) and (b) are satisfied for $B' = B + \epsilon H$ and any $0 < \epsilon \ll 1$. 
By a standard argument, we can pick $0 < \epsilon \ll 1$ so that (c) is also satisfied\footnote{Specifically, we pick $m \gg 0$ such that the Cartier index of $p^mH$ is not divisible by $p$, and then we can take, for example, $\epsilon =\frac{1}{p^m+1}$. The Cartier index of $K_Z+C+B + \epsilon H = K_Z+C+B+H - \frac{p^m}{p^m+1}H$ is not divisible by $p$, and so there exists $e>0$ such that $(p^e-1)(K_Z+C+B+\epsilon H)$ is Cartier.}}. 
\end{proof}

We can now move on {to} the main part of the proof. Since 
\[
(K_Z+C) \cdot C \leq (K_Z+C+B) \cdot C < 0,
\]
we get that $C \cong \bP^1$ (see \cite[Theorem 3.19]{tanaka12}). As $C$ is normal, we can write $K_C + B_C = (K_Z+C+B)|_C$ and $K_C + B'_C = (K_Z+C+B')|_C$. We claim that $(C,B'_C)$ is $F$-split. Indeed,
\begin{itemize}
    \item if $(C,B_C)$ is klt, then $(C,B_C)$ is globally $F$-regular by \cite{watanabe91} (as $B_C$ has standard coefficients), and so is $(C,B'_C)$ given that $\epsilon$ in the proof of the above claim can be chosen sufficiently small; 
    \item if $(C,B_C)$ is not klt, then the log canonical pair $(C,B'_C)$ is not klt either, and so it is  $F$-split by inversion of $F$-adjunction (see \cite[Lemma 2.9]{CTW15b}). 
\end{itemize} 
Finally, by {(c)}, we can apply the inversion of $F$-adjunction (see \cite[Lemma 2.7]{CTW15b}) to get that $(Z, C+B')$ is  $F$-split. {In particular, it is sharply  $F$-split, and thus so is $(Z,C+B)$.}
\end{proof}

\begin{lemma} \label{lem:F-split-inversion-of-adjunction-kappa-one} 
Let $(Z,C+B)$ be a two-dimensional projective 
{lc} pair over {an algebraically closed field} {$k$} 
of characteristic $p>{N_2}$, where $C$ is a prime divisor and $B$ is an effective $\bQ$-divisor with standard coefficients. 
Suppose that
\begin{enumerate}
    \item $(Z,aC+B)$ is log del Pezzo
    for some $a \geq 0$ and
    \item $-(K_Z+C+B)$ is nef, and $\kappa(Z, -(K_Z+C+B))=1$.
\end{enumerate}
Then $(Z,C+B)$ is sharply $F$-split. 
\end{lemma}

\begin{proof}
First, we reduce the problem to the case when the following additional assumption holds.
\begin{assumption} \label{ass:kappa=1}
There exists  a projective morphism $\pi \colon Z \to T$ to  a smooth projective curve $T$ such that $\pi_*\cO_Z = \cO_T$, $\rho(Z/T)=1$, and 
    $-(K_Z+C+B) \sim_{\Q} \pi^*A$ for some ample $\Q$-divisor $A$ on $T$.    
\end{assumption}
\noindent To this end, we note that by (2), there exists a morphism $\pi \colon Z \to T$ to 
a smooth projective curve $T$ such that $\pi_*\cO_Z = \cO_T$ and $-(K_Z+C+B) \sim_{\Q} \pi^*A$ for some ample $\Q$-divisor $A$ on $T$.
By running a $K_Z$-MMP over $T$, we get the following factorisation 
\[
\pi \colon Z \xrightarrow{h} Z' \xrightarrow{\pi'} T,
\]
where $\pi' \colon Z' \to T$ is the last step of this MMP. Set $C':=h_*C$ and $B':=h_*B$. Then $(Z',C'+B')$ satisfies all the assumptions of the lemma as well as Assumption \ref{ass:kappa=1}. If we show that $(Z',C'+B')$ is sharply $F$-split, then in view of the equality $K_Z+C+B=h^*(K_{Z'}+C'+B')$ and \cite[Lemma 3.3]{GT16}, we will get that the pair $(Z,C+B)$ is also sharply $F$-split. Therefore, by replacing $(Z,C+B)$ by $(Z',C'+B')$ we may assume  for now on that Assumption \ref{ass:kappa=1} holds.\\

Since $(Z, B)$ is klt,
$(Z,C+B)$ is sharply $F$-pure by Lemma \ref{lem:lc-is-sharply-F-pure}. Also,  note that $C$ is $\pi$-horizontal as $K_Z+C+B \equiv_{{\pi}} 0$ and 
$-(K_Z+aC+B)$ is ample.
Moreover, $C \cong \bP^1$ by \cite[Theorem 3.19]{tanaka12} as 
\[
(K_Z+C) \cdot C \leq (K_Z+C+B) \cdot C < 0.
\]

\begin{claim} There exists a $\bQ$-divisor $B' \geq B$ such $(Z,C+B')$ is log canonical, $(p^e-1)(K_Z+C+B')$ is Cartier for some $e>0$, and $K_Z+C+B' \equiv_{\pi} 0$.
\end{claim}
\begin{proof}[Proof of the claim]
Write $B = B_{\mathrm{hor}} + B_{\mathrm{ver}}$ as the sum of the horizontal and vertical parts. 
By (2) in Setting \ref{setting:del-pezzo-large-characteristic},  the coefficients of $B_{\mathrm{hor}}$ are contained in the set $\{1-\frac{1}{n} \mid 1 \leq n \leq {N_2} \}$, and so their denominators are not divisible by $p$. Write $B_{\mathrm{ver}} = \sum a_iF_i$ for ${0<a_i <1}$ and distinct {prime} divisors $F_i$. Since $\rho(Z/T)=1$, fibres of $\pi$ are irreducible, and so $F_i \equiv_{\pi} 0$.

We will modify $B_{\mathrm{ver}}$ so that the denominators of its coefficients are not divisible by $p$. To this end, define $a'_i \in \bQ$ as follows:
\begin{alignat*}{2}
    a'_i &:= a_i \qquad &&{\text{if}\quad 0 < a_i \leq {\textstyle 1 - \frac{1}{N}} }\\
    \vspace{0.2em} a'_i &:= a_i + \epsilon_i  \qquad &&{\text{if}\quad {\textstyle 1 - \frac{1}{N}} < a_i  <1}
\end{alignat*}
where we pick rational {numbers} $0 \leq \epsilon_i \ll 1$ such that the denominator{s} of $a'_i=a_i + \epsilon_i$ {are} not divisible by $p$. Set $B' =  B_{\mathrm{hor}} + \sum a'_iF_i$. As desired, $(p^e-1)B'$ is {a Weil divisor} for some $e>0$ and $K_Z+C+B' \equiv_{\pi} 0$ given that $F_i \equiv_{\pi} 0$. 

Before moving forward, we make two observations:
\begin{enumerate}
    \item[(i)] $(Z, C+B')$ is sharply $F$-pure and
    \item[(ii)] $-(K_Z+C+B') \sim_{\bQ} \pi^*A'$ for an ample $\bQ$-divisor $A'$ on $T$
\end{enumerate}
Observation {(i)} can be checked locally along every fibre. {Pick an integer $i>0$. If $a_i \leq 1-\frac{1}{N}$, then $a'_i = a_i$, and so $(Z, C+B')$ is sharply $F$-pure along $F_i$ as so is $(Z,C+B)$.} If $a_i > 1 - \frac{1}{N}$, then $(Z, C+B_{\mathrm{hor}} + F_i)$ is lc along $F_i$ by Lemma \ref{lem:kollar-acc} as $N>5$, and so sharply $F$-pure by Lemma \ref{lem:lc-is-sharply-F-pure}. 
Observation {(ii)} 
is automatic as $\epsilon_i$ are chosen sufficiently small and $K_Z+C+B'$ descends to $T$ by {the base point free theorem} \cite[Theorem 4.4(3)]{tanaka16_excellent},    
which can be applied as $Z$ is of del Pezzo type.

Therefore, we can invoke \cite[Lemma 2.12]{CTW15b} to get the existence of an effective $\bQ$-divisor $E$ whose support is contained in some fibres of $\pi$ such that $(Z,C+B'+E)$ is log canonical and $(p^e-1)(K_Z+C+B'+E)$ is Cartier for some $e>0$. Replacing $B'$ by $B'+E$ concludes the proof of the claim.
\qedhere
\end{proof}

As in the proof of the claim, $K_Z+C+B'$ descends to $T$. Now, by a standard argument, we can replace $B'$ by $(1-t)B + tB'$ for some $0 < t \ll 1$, so the statement of the claim is still satisfied\footnote{{Specifically, we pick $m \gg 0$ such that the Cartier index of $p^m(K_Z+C+B)$ is not divisible by $p$, and then we can take, for example, $t =\frac{1}{p^m+1}$. For this choice of $t$, the Cartier index of $K_Z+C+(1-t)B + tB' = \frac{p^m}{p^m+1}(K_Z+C+B) + \frac{1}{p^m+1}(K_Z+C+B')$ is not divisible by $p$.}}. In particular, $K_Z+C+B' \sim_{\bQ} -\pi^*A'$ for some ample $\bQ$-divisor $A'$ on $T$.

 Write $K_{C} + B_C = (K_Z+C+B)|_C$ and  $K_{C} + B'_{C} = (K_Z+C+B')|_{C}$. We claim that $(C,B'_C)$ is  $F$-split. To this end, we consider two cases. First, if $(C,B_C)$ is klt, then it is globally $F$-regular by \cite[Theorem 4.2]{watanabe91} as its coefficients are standard, and so $(C',B'_C)$ is globally $F$-regular as well as $t$ was chosen sufficiently small. Second, if $(C,B_C)$ is not klt, then $(C,B'_C)$ is not klt either (but still log canonical as $(Z,C+B')$ is log canonical), and so it is $F$-split by \cite[Lemma 2.7]{CTW15b} 
 {(here, $-(K_C+B'_C)$ is ample)}. {This lemma}  requires that $(p^e-1)(K_C+B'_C)$ is Cartier for some $e>0$, which holds by the claim.

Therefore, $(Z,C+B')$ is sharply $F$-split by {a subtle variant of} an inversion of $F$-adjunction \cite[Lemma 2.8]{CTW15b}, and thus $(Z,C+B)$ is  sharply $F$-split.
\qedhere
\end{proof}

\begin{lemma} \label{lem:log-trivial-case-plt}
Let $(Z,C+B)$ be a two-dimensional projective plt pair over a perfect field $k$ of characteristic $p>5$, where $C$ is a prime divisor and $B$ is an effective $\bQ$-divisor with standard coefficients. Suppose that
\begin{enumerate}
    \item $(Z,aC+B)$ is log del Pezzo for some $a \geq 0$ and
    \item $K_Z+C+B \equiv 0$.
\end{enumerate}
Then $(Z,C+B)$ is purely $2$-quasi-$F$-split. 
\end{lemma}

\begin{proof}
{In what follows, we verify that the assumptions of Theorem \ref{thm:log-trivial-inversion-of-adjunction} hold for $(Z,C+B)$. Assumption (1) is automatic. Assumption (2) follows from the fact that $(Z,B)$ is a two-dimensional klt pair with standard coefficients in characteristic $p>5$, and so it is strongly $F$-regular. As for Assumption (3), we have that} $(C,B_C)$ is $2$-quasi-$F$-split, 
where $K_C + B_C = (K_Z+C+B)|_C \equiv 0$. Indeed, to check this, we may assume by Corollary  \ref{c-descent2} that $k$ is algebraically closed. Here, $C$ is normal, as $(Z,C+B)$ is a two-dimensional plt pair. Since $K_C+B_C \equiv 0$, the classification of curves show that $C \cong \bP^1$ or $C$ is an elliptic curve. Thus $(C,B_C)$ is $2$-quasi-$F$-split by Proposition \ref{prop:p1-log-calabi-yau} and Remark \ref{remark:elliptic-curves-quasi-F-split}, respectively.
Finally, Assumption (4) is a consequence of the following equalities
\begin{align*}
H^1(Z,\,&\cO_Z(K_{Z} + \rup{-p^i(K_Z+C+B)})) \\
&= H^1(Z, \cO_Z(-\rup{-p^i(K_Z+C+B)})) \\
&= H^1(Z, \cO_Z(-\{p^iB\} +p^i(K_Z+C+B)))  \\
&= H^1(Z, \cO_Z(K_Z+aC + (B-\{p^iB\}) +p^i(K_Z+C+B) - (K_Z+aC+B)))  \\
&= 0.
\end{align*}
Here the first equality follows from Serre duality,
and the last one holds {by Kawamata-Viehweg vanishing \cite[Theorem 1.1]{ABL} given that $p^i(K_Z+C+B)-(K_Z+aC+B)$ is ample and $(Z, aC + (B-\{p^iB\})$ is klt {in view of} $\{p^iB\} \leq B$. This concludes the verification of the assumptions of Theorem \ref{thm:log-trivial-inversion-of-adjunction}, and so of the proof that} $(Z,C+B)$ is purely $2$-quasi-$F$-split.
\end{proof}

\begin{theorem}  \label{thm:qFsplit-del-Pezzo}
Let $(X,\Delta)$ be a log del Pezzo pair over a perfect field of characteristic $p > N$ such that $\Delta$ has standard coefficients. Suppose that $(X,\Delta)$ is not $\frac{1}{N}$-klt. Then $(X,\Delta)$ is $2$-quasi-$F$-split.
\end{theorem}

\begin{proof}
{By taking the base change to the algebraic closure of $k$, we may assume that $k$ is algebraically closed.} Since $(X, \Delta)$ is not $\frac{1}{N}$-klt, there exists 
a birational morphism $f \colon Y \to X$ {(possibly $f=\mathrm{id}$)} 
from a normal projective surface $Y$ 
and a prime divisor $C$ on $Y$ such that
$\Exc(f) \subseteq C$ 
and the log discrepancy of $C$ is at most  $\frac{1}{N}$, that is, 
$K_Y  + bC +B_Y = f^*(K_X+\Delta)$, 
{where $C \not \subseteq \Supp B_Y$, $B_Y \geq 0$, and}  $1 - \frac{1}{N} < b <1$. 
As before, we get that $(Y,C+B_Y)$ is log canonical  (see Lemma \ref{lem:kollar-acc}). 
Note that $-(K_Y + c C + B_Y)$ is ample for some $c \in \Q$ {satisfying} $b  <c <1$. 

Since $Y$ is of del Pezzo type, we may run a $-(K_Y+C+B_Y)$-MMP {for $Y_0 := Y$}
\begin{equation}\label{e1-qFsplit-del-Pezzo}
g\colon Y_0 \xrightarrow{g_0} Y_1 \xrightarrow{g_1}  \cdots \xrightarrow{g_{n-1}} Y_n. 
\end{equation}
Set $Z := Y_n$. Note that $C_Z := g_* C$ is still a prime divisor. Indeed, 
this MMP 
is $(K_Y+C+B_Y)$-positive and $(K_Y + c C + B_Y)$-negative. In particular, it is  
$C$-positive, that is, 
$C_{\ell} \cdot \Exc(g_{\ell}) >0$ { for all $0 \leq  {\ell} < n$, where $C_{\ell}$ is the image of} $C$ on $Y_{\ell}$. 

Set $B_Z := g_*B$. Note that $(Z, cC_Z + B_Z)$ is a del Pezzo pair and $(Z,C_Z+B_Z)$ is log canonical by Lemma \ref{lem:kollar-acc} again.  Now, we shall show that the MMP ends up with a minimal model, in other words, that $-(K_Z+C_Z+B_Z)$ is nef. 
By contradiction, suppose that this is not the case, and so $Z$ admits a Mori fibre space $h \colon Z \to T$. { Specifically, $h$ is a projective morphism to a normal projective variety $T$ such that $\dim T \leq 1$, $h_*\cO_{Z}=\cO_T$,} $\rho(Z/T)=1$, and $K_Z+C_Z+B_Z$ {is} $h$-ample. Since $-(K_Z + cC_Z + B_Z)$ is $h$-ample, there exists $t \in \Q$ such that 
$c < t< 1$ and $K_Z+tC_Z + B_Z \equiv_h 0$. This contradicts Condition (2) in Setting \ref{setting:del-pezzo-large-characteristic} as 
${1 - \frac{1}{N} < b < c < t < 1}$.

\begin{claim}\label{c-qFsplit-del-Pezzo}
$(Z,C_Z+B_Z)$ is purely $2$-quasi-$F$-split. 
\end{claim}

\begin{proof}[Proof of Claim \ref{c-qFsplit-del-Pezzo}]
If $\kappa(Z, -(K_Z+C_Z+B_Z)) =2$ or $1$, then $(Z, C_Z+ B_Z)$ is sharply $F$-split by Lemma \ref{lem:F-split-inversion-of-adjunction-kappa-two} and Lemma \ref{lem:F-split-inversion-of-adjunction-kappa-one}, respectively. 
Thus, for now on, we may assume that $\kappa({Z}, -(K_Z+C_Z+B_Z)) =0$.

By Condition (1) in Setting \ref{setting:del-pezzo-large-characteristic}, we get that $N(K_Z+C_Z+B_Z)$ is Cartier. Since $p>N$, there exists $e>0$ such that $(p^e-1)(K_Z+C_Z+B_Z)$ is Cartier. Thus, if $(Z,C_Z+B_Z)$ is not plt, then it is $F$-split by \cite[Lemma 2.9]{CTW15b}. In particular, it is sharply $F$-split (as $(p^e-1)B_Z$ is Weil), and so purely $1$-quasi-$F$-split {(the verification of the last assertion is similar to the proof of (\ref{e1-Fsplit}) and left to the reader)}.

Hence, we may assume that $(Z,C_Z+B_Z)$ is plt. Then it is $2$-quasi-$F$-split by Lemma \ref{lem:log-trivial-case-plt}.
\end{proof}

We complete the proof of the theorem. Pick $m\gg 0$ coprime with $p$ and write 
\begin{align*}
K_Y+\Phi_m &= g^*(K_{Z}+(1-1/m)C_{Z}+B_{Z}), \\
K_Y+\Phi &= g^*(K_{Z}+C_{Z}+B_{Z}).
\end{align*}Since 
$(Z, C_Z+B_Z)$ is the output of the $-(K_Y+C+B_Y)$-MMP ({see} (\ref{e1-qFsplit-del-Pezzo})), 
the negativity lemma implies that $\Phi \geq C+B_Y$. Moreover, $\Exc(g) \subseteq \Supp\,(\Phi -(C+B_Y))$, because 
{this} MMP is $(\Phi-(C+B_Y))$-negative {given that} $\Phi -(C+B_Y) = (K_Y+\Phi) -(K_Y+C+B_Y)$. 
Therefore, $\Phi_m \geq bC+B_Y$ for $m \gg 0$.

The pair $(Z, (1-\frac{1}{m})C_{Z}+B_{Z})$ is $2$-quasi-$F$-split  
by the claim {(and Lemma \ref{l-qS^0-adj-nonadj})}, and so $(Y,\Phi_m)$ is $2$-quasi-$F$-split by {Corollary} \ref{corollary:log-pullback-quasi-F-split-dimension-two}. In particular, $(Y, bC+B_Y)$ is $2$-quasi-$F$-split. By Proposition \ref{prop:pushing-forward-quasi-F-splittings}, 
$(X,\Delta)$ is $2$-quasi-$F$-split as well, concluding the proof.
\end{proof}

\begin{theorem} \label{thm:qFsplit-del-Pezzo-in-large-characteristic}
There exists a {positive integer} $p_0$ satisfying the following property. Let $(X,\Delta)$ be a {log del Pezzo} pair over a perfect field of characteristic $p{\cora \geq} p_0$ such that $\Delta$ has standard coefficients. Then $(X,\Delta)$ is $2$-quasi-$F$-split. 
\end{theorem}

\begin{proof}
If $(X,\Delta)$ is not $\frac{1}{N}$-klt, then the statement follows from Theorem \ref{thm:qFsplit-del-Pezzo}. Thus, we may assume that $(X,\Delta)$ is $\frac{1}{N}$-klt. In this case, 
we can find $p_0$ independent of $(X,\Delta)$ such that $(X,\Delta)$ is globally $F$-regular if ${\cora p \geq p_0}$ as in \cite[Remark 3.3]{CTW15b}.

{\cora In fact, by \cite[Lemma 3.1]{CTW15b}, there exists a flat family of log del Pezzo surfaces over a reduced quasi-projective scheme $T$ over $\Spec\,\Z$, such that every $\frac{1}{N}$-klt log del Pezzo surfaces over an algebraically closed field of characteristic $p>5$ appears as a fibre.
Then, replacing $T$ with $T\times_{\Spec\,\Z}\Spec\,\Z[1/M]$ for an appropriate $M\in \Z$, we can conclude from \cite[Theorem 1.2]{schwedesmith10} that every fibre is globally $F$-regular.}
\qedhere
\end{proof}

\subsection{Proof of Theorem \ref{thm:intro-3dim-klt}}

Combining  multiple results of this article, we are finally able to prove the following.

\begin{theorem} \label{thm:qFsplit-3dim-klt}
There exists a {positive integer} $p_0$ such that the following holds. 
Let $k$ be a perfect field of characteristic $p {\cora \geq}p_0$ and 
let $(X, \Delta)$ be a three-dimensional $\Q$-factorial affine klt pair of finite type 
over $k$, where $\Delta$ has standard coefficients. 

Then $(X,\Delta)$ is quasi-$F$-split. In particular, $X$ lifts {to} $W_2(k)$.
\end{theorem}

\begin{proof}
In what follows, we assume that $p>5$.  
Pick a {closed} point $x \in X$. After replacing $X$ by an affine open neighbourhood of $x \in X$, 
there exists a projective birational morphism $f \colon Y \to X$ 
such that $x \in f(\Exc(f))$, $E := \Exc(f)$ is a prime divisor, $-(K_Y+E+f^{-1}_*\Delta)$ is ample, and $(Y, E+ f_*^{-1}\Delta)$ 
is a $\Q$-factorial plt pair. 
It is constructed by taking a log resolution $g \colon W \to X$ of $(X,\Delta)$ and running a $(K_W + \Exc(g) + {g}^{-1}_*\Delta)$-MMP 
 {over $X$} \cite[Theorem 1.1]{hnt} 
so that $Y \to X$ is the last step of this {MMP}.

Since $E$ is normal (see, for example, \cite[Theorem 2.11]{GNT16}), 
we can write $K_E + \Delta_E = (K_Y+E+f^{-1}_*\Delta)|_E$. If $f(E)$ is a point, then there exists $p_0$ such that $(E,\Delta_E)$ is quasi-$F$-split by {Theorem} \ref{thm:qFsplit-del-Pezzo-in-large-characteristic}. 
If $f(E)$ is a curve, then 
$(E,\Delta_E)$ is 
quasi-$F$-split by Corollary \ref{c-QFS-2dim-rel}. 

Therefore, $(Y,E+f^{-1}_*\Delta)$ is purely quasi-$F$-split by Corollary \ref{cor:inversion-of-adjunction-anti-ample}; {in what follows, we verify its assumptions. Assumption (1) holds as $\rho(Y/X)=1$, so $-E$ is ample.} 
{Assumption (2) holds by the above paragraph. 
Assumption (3) is valid} as the relative Kawamata--Viehweg vanishing holds by 
\cite[Proposition {\cora 2.3}]{BK20} {(cf.\ \cite[Proof of Theorem 1.1]{HW17})}, where the {required} 
vanishing on $E$ is assured by \cite[Theorem 1.1]{ABL}. 
{Assumption (4) holds} as $(Y,E+f^{-1}_*\Delta)$ is {(locally)} purely $F$-regular by \cite[Lemma {\cora 4.2}]{HW20} applied to the identity morphism; in particular, $(Y,f^{-1}_*\Delta)$ is naively keenly $F$-{pure}. 
{Last, as $Y$ is $\bQ$-factorial and strongly $F$-regular, it is divisorially Cohen-Macaulay (see Remark \ref{remark:divisorially-Cohen-Macaulay}), and since, in addition, $\Delta$ has standard coefficients,  $(Y,E+f^{-1}_*\Delta)$ is $p$-compatible. This concludes the verification of the assumptions of Corollary \ref{cor:inversion-of-adjunction-anti-ample}.}

Finally, $(X,\Delta)$ is quasi-$F$-split by {Lemma \ref{l-qS^0-adj-nonadj} and} Proposition \ref{prop:pushing-forward-quasi-F-splittings}. The last part on the liftability of $X$  follows from \cite[Corollary A.1.2]{achinger-zdanowicz21} (see Subsection \ref{ss:canonical-lift}). 
\end{proof}

Based on \cite{CTW15b} and \cite{lacini20}, Arvidsson--Bernasconi--Lacini showed that 
klt del Pezzo surfaces are log liftable for $p>5$ 
{\cite[Theorem 1.2]{ABL}}. 
Assuming a natural strengthening of their result to log pairs with standard coefficients in the form of the following conjecture, one can show that $p_0= {\cora 43}$ works in the above theorem.
\begin{conjecture} \label{conj:log-liftability}
Let $(X,\Delta)$ be a log del Pezzo pair over a perfect field of characteristic $p>5$ and such that $\Delta$ has standard coefficients. Then there exists {a} 
log resolution $f \colon Y \to X$ of $(X, \Delta)$ such that $(Y, f^{-1}_*\Delta + \Exc(f))$ lifts to {$W(k)$}. 
\end{conjecture}

\begin{remark} \label{rem:main-theorem-for-42}
Assuming the above conjecture, we explain the idea of the proof  that log del Pezzo pairs $(X,\Delta)$ over perfect fields of characteristic 
${\cora p \geq p_0}$ with standard coefficients are quasi-$F$-split, which is a generalisation of {Theorem \ref{thm:qFsplit-del-Pezzo-in-large-characteristic}}. 
This immediately implies Theorem \ref{thm:qFsplit-3dim-klt} for $p_0= {\cora 43}$. {We expand on these ideas in \cite{KTTWYY2}.}

Let $f \colon Y \to X$ be a log resolution of $(X,\Delta)$. Write $K_Y+\Delta_Y = f^*(K_X+\Delta)$ and perturb $\Delta_Y$ by exceptional divisors so that $-(K_Y+\Delta_Y)$ is ample and $f_*\Delta_Y=\Delta$.  We will verify the assumptions of Theorem \ref{thm:higher-Cartier-criterion-for-quasi-F-split} to get that $(Y,\Delta_Y)$ is quasi-$F$-split which immediately yields quasi-$F$-spliteness of $(X,\Delta)$.

By Conjecture \ref{conj:log-liftability}, we have that $(Y,D)$ {lifts to $W(k)$}, where $D = \Supp f^{-1}_*\Delta + \Exc(f)$. Thus 
\[
H^0(Y, \Omega^1_Y(\log D)(\rdown{K_Y+\Delta_Y}))=0
\]
by Akizuki--Nakano vanishing {(see, e.g.,\ \cite[Theorem 2.4]{Kawakami-Nagaoka})}. Moreover,
\[
H^0(X, \cO_X(K_Y + \Supp(\Delta_Y) + \rdown{p^k(K_Y+\Delta_Y)}))=0
\]
for every $k\geq 1$ as can be easily deduced by reducing to a (relative) Picard one case through running a {$K_Y$-}MMP and applying Condition (2) from Setting \ref{setting:del-pezzo-large-characteristic}. That condition is valid for $N=42$ by \cite{kollar94} in characteristic zero, and this immediately extends to positive characteristic as we are assuming log liftability (Conjecture \ref{conj:log-liftability}). This concludes the verification of the assumptions of Theorem \ref{thm:higher-Cartier-criterion-for-quasi-F-split}.
\end{remark}

%% file: section7.tex
\section{Appendix}

The content of this section is independent of the other parts of the article and it only builds on the material from the preliminaries (Section \ref{s-Prelim}).
\subsection{Calculations using local cohomology} \label{ss:calculation-local-cohomology}
We would like to emphasise that, despite a complicated definition, whether a ring is  quasi-$F$-split can be calculated explicitly. The calculation below was one of the key motivations for our work at the onset of the project. We do not need it in our article, but we believe that it is helpful in gaining an insight into how quasi-$F$-splittings work.

\begin{setup*}
Let $k$ be a perfect field of characteristic $p>0$ and let $f \in k[x,y,z]$ be a polynomial of the form $f = z^2 + g$ for $g \in {(x,y)}k[x,y]$. Note that all rational double point singularities are of this form. Set $R := k\llbracket{x,y,z}\rrbracket/(f)$ and assume that $R$ is a domain. We denote the maximal ideal of $R$ by $\m$.
\end{setup*}

Suppose we want to check whether such an $R$ is $2$-quasi-$F$-split. By Lemma \ref{lem:intro-cohomological-quasi-F-spliteness}(2), this is equivalent to verifying whether 
\begin{equation} \label{eq:example-local-cohomology-injective}
    \Phi \colon H^2_{\m}(R) 
    \to H^2_{\m}(Q_{R,2})
\end{equation} 
is injective (note that $R$ is Gorenstein). Since the radical of the ideal $(x,y)$ is equal to $\m$, we can calculate the local cohomologies with respect to $(x,y)$.

\begin{remark}
Recall that given an $R$-module $M$, we have
\[
H^2_{(x,y)}(M) = \cokernel\bigg(M\Big[\frac{1}{x}\Big] \oplus M\Big[\frac{1}{y}\Big] \to M\Big[\frac{1}{xy}\Big]\bigg),
\]
where the map {is given by subtraction.} 
We will denote local cohomology classes in $H^2_{(x,y)}(M)$ by the symbol $\{\frac{m}{x^iy^j}\}$ for $m \in M$ and $i,j \geq  0$.
\end{remark}

\noindent \textbf{Observation 1.} $R$ is $2$-quasi-$F$-split if and only if $\Phi\big\{\frac{z}{xy}\big\} \in H^2_{\m}(Q_{R,2})$ is non-zero.\\

\noindent Indeed, since $R = k\llbracket{x,y}\rrbracket \oplus k\llbracket{x,y}\rrbracket z$ as a $k\llbracket{x,y}\rrbracket$-module, it is free, and so
\[
H^2_{(x,y)}(R) = H^2_{(x,y)}(k\llbracket{x,y}\rrbracket) \otimes_{k\llbracket{x,y}\rrbracket} R = H^2_{(x,y)}(k\llbracket{x,y}\rrbracket) \oplus  H^2_{(x,y)}(k\llbracket{x,y}\rrbracket)z.
\]
{Since $g \in (x,y)k[x,y]$, one can deduce from this that} $\{\frac{z}{xy}\}$ is a multiple of every non-zero element in $H^2_{(x,y)}(R)$ 
{\cora (i.e., this element generates the socle of $H^2_{(x,y)}(R)$).}\footnote{The existence of such a $\{\frac{z}{xy}\}$ is nothing special; {\cora the socle of the highest local cohomology of every Gorenstein local domain is generated by one element}.}
Hence, if some local cohomology class {in} $H^2_{\m}(R)$ lies in the kernel of $\Phi$, then so does $\{\frac{z}{xy}\}$.  \\

\noindent \textbf{Observation 2.} We have that $\Phi\big\{\frac{z}{xy}\big\} = \big\{\frac{[z]^p}{[xy]^p}\big\}$.\\

\noindent Here, we are abusing the notation and denoting the image in $Q_{R,2}$ of a Teichm\"uller lift $[g] \in F_* W_2(R)$ of $g \in R$ by the same symbol. In order to check this equality, consider the following commutative diagram
\begin{center}
\begin{tikzcd}
H^2_{([x],[y])}(W_2(R)) \arrow{r}{F} \arrow{d} & H^2_{([x],[y])}(F_*W_2(R)) \arrow{d} \\
H^2_{(x,y)}(R) \arrow{r}{\Phi} & H^2_{(x,y)}(Q_{R,2}),
\end{tikzcd}
\end{center}
where $[x], [y] \in W_2(R)$ are the Teichm\"uller lifts of $x$ and $y$, respectively. To calculate $\Phi \big\{\frac{z}{xy}\big\}$, we can pick an arbitrary lift of $\big\{\frac{z}{xy}\big\}$ to $H^2_{([x],[y])}(W_2(R))$, then apply $F$, and finally take the image in $H^2_{(x,y)}(Q_{R,2})$. Taking the lift of $\big\{\frac{z}{xy}\big\}$ to be $\big\{\frac{[z]}{[xy]}\big\}$ immediately yields Observation 2.

\begin{example}\label{ex1-RDP}
Let $k$ be a perfect field of characteristic $p=3$ and let $f = z^2+x^3+y^4 \in 
k[x, y, z]$. 
We will show that $R = k\llbracket{x,y,z}\rrbracket/(f)$ is $2$-quasi-$F$-split.

First, we note that by Fedder's criterion, $R$ is not $F$-split, and so is not $1$-quasi-$F$-split. Equivalently, we can check that $F \colon H^2_{(x,y)}(R) \to H^2_{(x,y)}(F_*R)$ is not injective by calculating
\[
F\Big\{\frac{z}{xy}\Big\} = \Big\{\frac{z^3}{(xy)^3}\Big\} = \Big\{\frac{-z(x^3+y^4)}{(xy)^3}\Big\} = -\Big\{\frac{z}{y^3}\Big\} - \Big\{\frac{zy}{x^3}\Big\} = 0.
\]

By Observations 1 and 2, in order to prove that $R$ is $2$-quasi-$F$-split, 
it is enough to verify that 
$\big\{\frac{[z]^3}{[xy]^3}\big\}$ is non-zero in $H^2_{\m}(Q_{R,2})$. Calculate
\begin{align*}
\Big\{\frac{[z]^3}{[xy]^3}\Big\} &= -\Big\{\frac{[z][x^3+y^4]}{[xy]^3}\Big\}\\
&= -\Big\{\frac{[z]([x^3]+[y^4] + V(x^6y^4 + x^3y^8))}{[xy]^3}\Big\}  \\
&= -\Big\{\frac{[z]}{[y]^3}\Big\} - \Big\{\frac{[zy]}{[x]^3}\Big\} - V\Big\{\frac{z^3(x^6y^4 + x^3y^8)}{(xy)^9}\Big\}  \\
&= V\Big\{\frac{z(x^3+y^4)(x^6y^4 + x^3y^8)}{(xy)^9}\Big\} \\
&= -V\Big\{\frac{z}{x^3y}\Big\}.
\end{align*}
Here, 
{\cora the first equality holds by $[z]^3 = [z^3] = [-z(x^3 +y^4)] = -[z][x^3+y^4]$ (note that $[-1]=-1$ is assured by $p>2$),} 
the second equality follows from (\ref{eq:addition-in-W2-using-V}) 
and the last equality is a consequence of $(x^3+y^4)(x^6y^4 + x^3y^8) \equiv 2x^6y^8 \pmod{(x^9,y^9)}$. 
To prove that the above element is non-zero in $H^2_{(x,y)}(Q_{R,2})$, consider the following diagram (see (\ref{eq:key-sequence-for-witt-vectors}) and  (\ref{eq:intro-C-restriction-sequence}))
\begin{center}
\begin{tikzcd}
0 \to H^2_{(x,y)}(F^2_*R) \arrow{d} \arrow{r}{V} & H^2_{([x],[y])}(F_*W_2(R)) \arrow{r} \arrow{d} & H^2_{(x,y)}(F_*R) \arrow{d} \\
0 \to H^2_{(x,y)}(F^2_*R/F_*R) \arrow{r} & H^2_{(x,y)}(Q_{R,2}) \arrow{r}& H^2_{(x,y)}(F_*R),
\end{tikzcd}
\end{center}
where the exactness on the left follows from the vanishing of $H^1_{(x,y)}(R)$ given that $R$ is a normal ring.

The above calculation shows that $\big\{\frac{[z]^3}{[xy]^3}\big\} \in H^2_{([x],[y])}(F_*W_2(R))$ is the image of $\big\{\frac{z}{x^3y}\big\} \in H^2_{(x,y)}(F^2_*R)$ under $V$. By tracing through the above diagram, we see that it is sufficient to show that the image of $\big\{\frac{z}{x^3y}\big\}$ in $H^2_{(x,y)}(F^2_*R/F_*R)$ is non-zero. In turn, by the exact sequence
\[
H^2_{(x,y)}(F_*R) \xrightarrow{F} H^2_{(x,y)}(F^2_*R) \to H^2_{(x,y)}(F^2_*R/F_*R) \to 0, 
\]
it is enough to show that $\big\{\frac{z}{x^3y}\big\}$ does not lie in the image of Frobenius $F \colon H^2_{(x,y)}(F_*R) \to H^2_{(x,y)}(F^2_*R)$. Recall that
\begin{align*}
&H^2_{(x,y)}(R) = H^2_{(x,y)}(k\llbracket{x,y}\rrbracket) \oplus  H^2_{(x,y)}(k\llbracket{x,y}\rrbracket)z \text{ and } \\
&H^2_{(x,y)}(k\llbracket{x,y}\rrbracket) = \bigoplus_{i,j > 0} k\Big\{\frac{1}{x^iy^j}\Big\}.
\end{align*}
These local cohomology groups come with natural multidegree with respect to $(x,y)$ such that the multidegree of $\frac{1}{x^iy^j}$ is $(i,j)$.

\begin{claim*} 
Every homogeneous element {in} the image of $F \colon H^2_{(x,y)}(R) \to H^2_{(x,y)}(F_*R)$ is 
of multidegree $(a,b)$, where $(a \mod 3, b \mod 3) = (0,0) \text{ or } (0,2)$.
\end{claim*}
\noindent Up to replacing $R$ by $F_*R$, this concludes the proof that $R$ is $2$-quasi-$F$-split as the multidegree of  $\big\{\frac{z}{x^3y}\big\}$ is congruent to $(0,1)$ modulo $3$.
\begin{proof}[Proof of the claim]
The above claim is true, because $F\big\{\frac{1}{x^{i}y^{j}}\big\} = \big\{\frac{1}{x^{3i}y^{3j}}\big\}$ is of multidegree congruent to $(0,0)$ modulo $3$, whilst $F\big\{\frac{z}{x^{i}y^{j}}\big\} = \big\{\frac{z^3}{x^{3i}y^{3j}}\big\} = -\big\{\frac{x^3 + y^4}{x^{3i}y^{3j}}\big\}z$ is a sum of classes whose multidegrees are congruent to $(0,0)$ and $(0,2)$ modulo $3$.
\end{proof}
\end{example}

In fact, by a similar local cohomology calculation (together with the definition of pure  {quasi-$F$-splittings} and adjunction, see Subsection \ref{ss:pure-quasi-F-split}), one  can obtain a criterion for when a homogeneous Calabi-Yau hypersurface is $2$-quasi-$F$-split. Note that this result is a corollary of a far more general and powerful result in \cite{KTY22}.
\begin{proposition}[{\cite[Theorem C]{KTY22}}] Let $f \in k[x_1,\ldots, x_n]$ be a homogeneous polynomial of degree $n$ and let $R = k[x_1,\ldots, x_n]/(f)$. Then $R$ is $2$-quasi-$F$-split in a neighbourhood of the origin if and only if 
\[
f^{p-1} \not \in (x_1^p, \ldots,x_n^p) \quad \text{ or } \quad f^{p^2-p-1}\Delta(f)  \not \in (x_1^{p^2}, \ldots, x_n^{p^2}).
\]
Here $\Delta(f) \in k[x_1,\ldots,x_n]$ is defined by the following formula in $W_2({k[x_1,\ldots,x_n]})$:
\[
[f(x_1,\ldots,x_n)]=f([x_1],\ldots,[x_n]) + V\Delta(f),
\]
where in $f([x_1],\ldots, [x_n])$,  we treat $f$ as a polynomial over $W_2(k)$ by lifting the coefficients of monomials via the Teichm\"uller lift $k \to W_2(k)$.  
\end{proposition}

Note that $f^{p-1} \not \in (x_1^p, \ldots,x_n^p)$ is nothing but  Fedder's criterion, and so it verifies whether $R$ is $1$-quasi-$F$-split. The quantity $\Delta(f)$ can be calculated explicitly:
\[
\Delta(f) = \frac{1}{p}\bigg(\Big(\sum_I a_Ix^{I}\Big)^p - \sum_I \big(a^Ix^I\big)^p\bigg),
\]
where $f = \sum_{I}a_Ix^I$ for $a_I \in k$ with $I=(i_1,\ldots, i_n) \in \bZ^{\oplus n}_{\geq 0}$ and $x^I = x_1^{i_1}\ldots x_n^{i_n}$.

The above result (for $2$-quasi-$F$-splittings) can also be deduced from the calculation of logarithms of formal groups of hypersurfaces (see  \cite{stienstra87}), by way of the method of the proof of \cite[Theorem D]{KTY22}.

\subsection{Canonical lifting modulo $p^2$ of  a quasi-$F$-split scheme}
\label{ss:canonical-lift}
We start by recalling the construction from \cite{achinger-zdanowicz21} of a lift associated to a quasi-$F$-splitting (which in turn is motivated by \cite{zdanowicz17}). Let $X$ be a scheme over a perfect field $k$ of characteristic $p>0$. 
Suppose that $X$ is $n$-quasi-$F$-split and let
\[
\sigma \colon F_*W_n\cO_X \to \cO_X
\]
be a {surjective} 
$W_n\cO_X$-module homomorphism  
featured in the definition of quasi-$F$-split schemes (see Diagram \ref{diagram:intro-definition}).
In other words, $\sigma$ is a $W_n\MO_X$-module homomorphism 
satisfying $\sigma({F_*}1)=1$. We call such $\sigma$ a {\em quasi-splitting}. 

Using $\sigma$, {one} can construct a scheme $X_{\sigma}$ which is a lift of  $X$ over $W_2(k)$. 
Specifically, consider the pushout of $V \colon F_*W_n\cO_X \to W_{n+1}\cO_X$ along $\sigma$ sitting in the following diagram of $W_{n+1}\MO_X$-modules
\begin{center}
\begin{tikzcd}
0 \arrow{r} & F_*W_n\cO_X \arrow{r}{V} \arrow{d}{\sigma} & W_{n+1}\cO_X \arrow{r}{R^n} \arrow{d}{{\tau}} & \cO_X \arrow{r} \arrow{d}{=} & 0 \\
0 \arrow{r} & \cO_X \arrow{r} & \arrow[lu, phantom, "\usebox\pushoutdr" , very near start, yshift=0em, xshift=0.3em, color=black] \cO_{\sigma} \arrow{r} & \cO_X \arrow{r} & 0.
\end{tikzcd}
\end{center}

The sheaf $\cO_{\sigma}$ is a sheaf of rings on $X$, because {$\cO_{\sigma} = W_{n+1}\cO_X /\,V(\kernel(\sigma))$ by definition of a pushout} and $V(\kernel(\sigma)) \subset W_{n+1}\cO_X$ is an ideal (cf.\ \cite[Lemma A.1.1]{achinger-zdanowicz21}). The latter is true, {because $V$ and $\sigma$ are naturally  $W_{n+1}\cO_X$-module homomorphisms.}

Since $\cO_{\sigma}$ is an extension of $\cO_X$ by itself, $\cO_{\sigma}$ is an algebra over $W_2(k)$\footnote{Explicitly, {$p$ lies in the image of the $W_{n+1}\cO_X$-homomorphism $\cO_X \to \MO_{\sigma}$} from the exact sequence $0 \to \MO_X \to \MO_{\sigma} \to \MO_X \to 0$; {thus} $p^2 =0$ in $\MO_{\sigma}$, and so the composite ring homomorphism 
$W_{n+1}(k) \to W_{n+1}\MO_X\xrightarrow{\tau} \MO_{\sigma}$ 
factors through $W_{n+1}(k) \to W_{n+1}(k)/p^2 W_{n+1}(k) =W_2(k)$
.}.
 We set $X_{\sigma}:=(|X|, \cO_{\sigma})$ to  be the associated scheme over $W_2(k)$. By the struture of the lower exact sequence in the above diagram, $X_{\sigma}$ is flat over $W_2(k)$, and so it is a lift of $X$ over $W_2(k)$ (see \cite[Corollary A.1.2]{achinger-zdanowicz21} for details).

In what follows, we  show that the above construction allows for lifting of every line bundle and every Cartier divisor modulo $p^2$. The latter result is new and somewhat surprising even for $F$-splittings.

\begin{theorem}\label{thm:p^2-lift of pairs}
Let $X$ be a quasi-$F$-split scheme over a perfect field $k$ of characteristic $p>0$  and let $\sigma \colon F_*W_n\cO_X \to \cO_X$ be a quasi-splitting, that is, 
a $W_n\cO_X$-module homomorphism satisfying $\sigma({F_*1})=1$. Let $X_{\sigma}$ be the associated scheme over $W_2(k)$ constructed above.
Then the following hold:
\begin{enumerate}
\item every invertible sheaf $\cL$ on $X$ admits a canonical lift $\cL_{\sigma}$ to $X_{\sigma}$, and
\item every Cartier divisor $D$ on $X$ admits a canonical lift $D_{\sigma}$ to $X_{\sigma}$.
\end{enumerate}
\end{theorem}

\begin{proof}
By construction, the Teichm\"uller lift induces a multiplicative section
\[
\cO_X \to \cO_{\sigma} \quad \text{ denoted } \quad f \mapsto \tilde{f}
\]
of $\cO_{\sigma} \to \cO_X$ which immediately yields
(1) (just use the above map to lift the transition functions of $\cL$).

For (2), given a local equation $f \in \cO_X$ of $D$, we define $D_{\sigma}$ locally by $\tilde f \in \cO_{\sigma}$. To conclude the proof, it is sufficient to show that 
$\tilde{f}$ is a non-zero-divisor in $\MO_{\sigma}$ and $\cO_{\sigma}/\tilde{f}$ is flat over $W_2(k)$. 

{\cora First, we show that $\tilde{f}$ is a non-zero-divisor in $\MO_{\sigma}$.
Recall that we have an exact sequence
\[
0\to \cO_{\sigma} \xrightarrow{\times p} \cO_{\sigma} \xrightarrow{r} \cO_X \to 0
\]
by the flatness of $\cO_{\sigma}$ over $W_2(k)$,
and an isomorphism $\cO_{X}\xrightarrow{\times p} p\cO_{\sigma}$ of $\cO_{\sigma}$-modules (see \cite[8.7 Properties]{esnault_viehweg}).

Take $\tilde{a}\in \mathcal{O}_{\sigma}$ such that $\tilde{a}\tilde{f}=0$.
We show $\tilde{a}=0$.
Since $f$ is a non-zero-divisor, it follows that $r(\tilde{a})=0$.
Thus we can take $\tilde{b}\in \cO_{\sigma}$ such that $\tilde{a}=p\tilde{b}$.
Since $p\tilde{b}\tilde{f}=0\in p\cO_{\sigma}$, we have $r(\tilde{b})f=0 \in \cO_{\sigma}$. Since $f$ is a non-zero-divisor, we have $r(\tilde{b})=0$, and there exists $\tilde{c}\in \cO_{\sigma}$ such that $\tilde{b}=p\tilde{c}$.
Thus, $\tilde{a}=p\tilde{b}=p^2\tilde{c}=0$, and we conclude.}

Next, $\cO_{\sigma}/\tilde{f}$ is flat over $W_2(k)$, because
from the exact sequence
\[
0 \to \cO_{\sigma} \xrightarrow{\tilde{f}} \cO_{\sigma} \to \cO_{\sigma}/\tilde{f} \to 0,
\]
we can see that $\Tor_1^{W_2(k)}(\cO_{\sigma}/\tilde{f}, k)=0$.
This concludes the proof by 
the local criterion for flatness (\cite[Theorem 22.3 (3')$\to$ (1)]{Matsumura}). 
\end{proof}

\begin{corollary}\label{cor:p^2-lift of snc pairs}
Let $X$ be a smooth variety over a perfect field $k$ of characteristic $p>0$ and let $D$ be a simple normal crossing divisor on $X$.
Assume that $X$ is quasi-F-split and let $\sigma \colon F_*W_n\cO_X \to \cO_X$ be a quasi-splitting, that is, a $W_n\cO_X$-module homomorphism
satisfying $\sigma(F_*1)=1$.

With notation as in Theorem \ref{thm:p^2-lift of pairs},
$(X_{\sigma}, D_{\sigma})$ is log smooth over $W_2(k)$. 
\end{corollary}

\begin{proof}
For any closed point $x \in X$, we have local coordinates $x_1, \dots, x_d \in \cO_X$ at $x \in X$ such that $D$ is defined by $x_1\cdots x_r=0$ for some $r \leq d$.
As in the proof of Theorem \ref{thm:p^2-lift of pairs}, we have lifts $\tilde{x}_1, \dots, \tilde{x}_d \in \cO_{\sigma}$ and $D_{\sigma}$ is defined by $\tilde{x}_1\cdots \tilde{x}_r=0$.
These lifts define the $W_2(k)$-algebra homomorphism 
\[
W_2(k)[X_1, \dots, X_d] \to \cO_{\sigma}, \qquad X_i \mapsto \tilde{x}_i, 
\]
which is \'etale since its mod $p$ reduction is \'etale (cf.\ the local criterion for flatness \cite[Theorem 22.3]{Matsumura}). 
This implies that $(X_{\sigma}, D_{\sigma})$ is log smooth over $W_2(k)$. 
\end{proof}

\begin{corollary}\label{c-QFS-KVV}
Let $X$ be a smooth projective variety over a perfect field $k$ of characteristic $p>0$ and let $D$ be a simple normal crossing divisor on $X$. 
Let $A$ be an ample $\Q$-divisor on $X$. 
Assume that $p \geq \dim X$, 
$X$ is quasi-F-split, and 
$\Supp \{A\} \subseteq \Supp D$. 
Then $H^j(X, \Omega_X^i(\log D) \otimes \MO_X( - \rup{A}))=0$ if $i + j < \dim X$. 
In particular, $H^j(X, \MO_X( - \rup{A}))=0$ for every $j < \dim X$. 
\end{corollary}

\begin{proof}
The assertion follows from Corollary \ref{cor:p^2-lift of snc pairs} and \cite[Theorem 2.4]{Kawakami-Nagaoka}. 
\end{proof}